\documentclass[a4paper,11pt,reqno]{amsart}
\usepackage[english]{babel}
\usepackage[a4paper,margin=2.5cm]{geometry}

\usepackage{amsfonts,amsmath,amsthm,amssymb,amsxtra,calligra,mathrsfs}
\usepackage{graphicx}
\usepackage{tikz}
\usetikzlibrary{automata, arrows.meta, positioning}
\usepackage{tikz-cd} %绘制交换图
\usepackage{xcolor}
\usepackage[all]{xy}
\usepackage{upgreek}
\usepackage{graphicx}
\usepackage{mathtools}
\usepackage{extarrows}
\usepackage{faktor}
\usepackage{diagbox}
\usepackage{tabularx}
\usepackage{comment}
\usepackage{enumitem}
\usepackage[justification=centering]{caption}
\usepackage[pagebackref]{hyperref}
\hypersetup{
    colorlinks=true, 
    citecolor=blue,
    linkcolor=magenta,
	bookmarksnumbered=true,
}

\renewcommand*{\backrefalt}[4]{%
    \ifcase #1 %
        (Not cited.)%
    \or
        $\uparrow$#2.%
    \else
        $\uparrow$#2.%
    \fi
}
\DeclareRobustCommand{\SkipTocEntry}[5]{}

\newtheorem{defn-pro}{Definition-Proposition}
\newtheorem{defn-thm}{Definition-Theorem}
\newtheorem{thm}{Theorem}[section]
\newtheorem{lem}[thm]{Lemma}

\newtheorem{pro}[thm]{Proposition}

\newtheorem{defn}[thm]{Definition}

\newtheorem{question}[thm]{Question}

\newtheorem{conj}[thm]{Conjecture}
\newtheorem{rem}[thm]{Remark}

\newtheorem{Pred}{Prediction}

\theoremstyle{remark}

\newcommand{\bbC}{\mathbb{C}}

\newcommand{\bbF}{\mathbb{F}}

\newcommand{\bG}{\mathbb{G}}

\newcommand{\PP}{\mathbb{P}}

\newcommand{\mcolon}{\mathbin{:}}
\newcommand{\defeq}{\vcentcolon=}

\newcommand{\inverse}[1]{{#1}^{-1}}

\renewcommand{\hom}[3]{\mathrm{Hom}_{#1}({#2},{#3})}
\newcommand{\sheafHom}{\mathscr{H}\text{\kern -3pt {\calligra\large om}}\,}

\newcommand{\mybrace}[2]{\left\{#1 \,\middle|\, #2 \right\}}

\newcommand{\ord}[2]{\mathrm{ord}_{#1}(#2)}

\newcommand{\abs}[1]{\lvert #1 \rvert}

\newcommand{\bZ}{{\mathbb Z}}
\newcommand{\bN}{{\mathbb N}}
\newcommand{\bC}{{\mathbb C}}
\newcommand{\bQ}{{\mathbb Q}}
\newcommand{\bP}{{\mathbb P}}
\newcommand{\bR}{{\mathbb R}}

\newcommand{\bL}{{\mathbb L}}

\newcommand{\q}{/\!\!/}
\newcommand{\cZ}{{\mathcal Z}}
\newcommand{\cO}{{\mathcal O}}

\newcommand{\cP}{\mathcal{P}}
\newcommand{\cB}{\mathcal{B}}

\newcommand{\fm}{{\mathfrak{m}}}

\newcommand{\fX}{{\mathfrak X}}

\newcommand{\rounddown}[1]{\lfloor{#1}\rfloor}

\newcommand\spa{\text{\rm span}}

\newcommand{\II}{I\!I}
\newcommand{\divi}{\mathrm{div}}
\newcommand{\stabop}{\widetilde{O}^+}
\newcommand{\MP}{\rm Mp}
\newcommand{\sstable}[1]{{#1}^\text{ss}}

\newcommand{\parallelsum}{\mathrel{\!/\mkern-5mu/}\!}

\newcommand{\cF}{{\mathcal F}}

\newcommand\cD{{\mathcal{D}}}

\newcommand\cI{{\mathcal{I}}}
\newcommand\cE{{\mathcal{E}}}

\newcommand\BH{\mathbb{H}}

\newcommand\cC{{\mathcal{C}}}

\DeclareMathOperator{\disc}{disc}

\DeclareMathOperator{\proj}{Proj}
\DeclareMathOperator{\Supp}{Supp}
\DeclareMathOperator{\rk}{rk}

\DeclareMathOperator{\stab}{Stab}

\DeclareMathOperator{\Hom}{Hom}

\DeclareMathOperator{\id}{\mathrm{id}}

\DeclareMathOperator{\PGL}{\mathrm{PGL}}
\DeclareMathOperator{\GL}{\mathrm{GL}}
\DeclareMathOperator{\SL}{\mathrm{SL}}

\DeclareMathOperator{\divib}{\mathrm{div}}

\DeclareMathOperator{\Pic}{\mathrm{Pic}}
\DeclareMathOperator{\Cl}{\mathrm{Cl}}
\DeclareMathOperator{\NS}{\mathrm{NS}}
\DeclareMathOperator{\NL}{\mathrm{NL}}

\DeclareMathOperator{\csch}{\mathcal{S}ch}
\DeclareMathOperator{\cset}{\mathcal{S}et}

\DeclareMathOperator{\acup}{Acusp}

\DeclareMathOperator{\pic}{Pic}
\DeclareMathOperator{\aut}{Aut}

\DeclareMathOperator{\im}{\text{\rm Im}}

\DeclareMathOperator{\rank}{{\text{\rm rank}}}

\DeclareMathOperator{\sym}{Sym}
\DeclareMathOperator{\sh}{\mathrm{Sh}}

\renewcommand{\theequation}{\thesection.\arabic{equation}}
\newcommand{\invLattice}{N}

\newcommand{\projection}[2]{\pi_{{#1}^\perp}(#2)}

\newcommand{\setupappendixcounters}[1]{%
    \renewcommand\thesection{#1}
    \refstepcounter{section}
    
    \setcounter{subsection}{0}
    \renewcommand{\thesubsection}{#1.\arabic{subsection}}
    \renewcommand{\theHsubsection}{#1.\arabic{subsection}}
    
    \setcounter{thm}{0}
    \renewcommand{\thethm}{#1.\arabic{thm}}
    \renewcommand{\theHthm}{#1.\arabic{thm}}
    
    \setcounter{equation}{0}
    \renewcommand{\theequation}{#1.\arabic{equation}}
    \renewcommand{\theHequation}{#1.\arabic{equation}}
}

\title{Birational geometry of moduli space of   Del Pezzo pairs}
\author{Long Pan}
\address{Academy of Mathematics and Systems Science, Beijing, China}
\email{panlong@amss.ac.cn}

\author{Fei Si}
\address{The School of Mathematics and Statistics
Xi’an Jiaotong University, 28 West Xianning Road, Xi’an, Shaanxi, P.R.
China, 710049}
\email{sifei@xjtu.edu.cn}

\author{Haoyu Wu}
\address{Yau Mathematical Sciences Center, Jingzhai, Tsinghua University, Haidian District, Beijing, China}
\email{haoyu-wu@tsinghua.edu.cn}

\date{}

\begin{document}
\begin{abstract}
    In this paper, we investigate the geometry of the moduli space $P_d$ of degree $d$ smooth del Pezzo pairs, which consists of a smooth del Pezzo surface $X$ of degree $d$ and a smooth curve $C \sim -2K_X$. 
    More precisely, we study the compactifications of $P_d$ from both Hodge-theoretic and geometric invariant theoretical (GIT) perspectives.  We obtain the class numbers of the Baily-Borel compactification $P_d^\ast$ for $P_d$,  which is an important step toward establishing the Hassett-Keel-Looijenga program for $P_d$. If $d=8$,
    $P_d$ has two connected components. For the component parametrizing del Pezzo pairs $(Bl_p \PP^2, C)$, we propose the Hassett-Keel-Looijenga models $\cF(s)=\proj R(\cF,\Delta(s) )$ via the section rings of certain $\bQ$-line bundles $\Delta(s)$ on the locally symmetric variety $\cF$. These models are expected to connect different birational models of the moduli space arising from K-moduli theory. 
    By constructing an arithmetic stratification on $\cF$ and computing the pullback of $\Delta(s)$ on these strata,  we give arithmetic predictions for the wall-crossing of $\cF(s)$ as $s\in [0,1]$ varies. 
    This work parallels that of Laza-O'Grady \cite{LO19, LaO18}.  %The relation of $\cF(s)$ with the K-moduli spaces of degree $8$ del Pezzo pairs is also discussed.
\end{abstract}

\maketitle
\setcounter{tocdepth}{1}
\tableofcontents

\section{Introduction}

\subsection{Background and motivation} Recent developments in moduli theory provide many compactifications for the moduli space of  (polarized or lattice-polarized) K3 surfaces from various perspectives.
For example,  there are geometric invariant theoretical(GIT)  compactifications by \cite{Shah80} \cite{Shah}.  On the Hodge-theoretic side,  the Torelli theorem for K3 surfaces serves as a bridge to link the moduli space of  K3 surfaces to a locally symmetric space. 
Then it turns out that there are many arithmetic compactifications for  the moduli space of lattice-polarized K3 surfaces, including the Baily-Borel  compactification and Looijenga's semitoric  compactifications \cite{Lo03}. 
Recent breakthroughs in  birational geometry and  K\"ahler-Einstein geometry have led to general solutions for constructing moduli spaces of log Fano varieties (see \cite{Xu21}) and  varieties of general type (see \cite{kollar}).
These works inspire the following results on compactifications of moduli spaces of K3 surfaces:   \cite{ADL22}  uses K-stability theory to provide a compactification for moduli spaces of quartic K3 surfaces. By constructing a recognizable divisor on a K3 surface,  Alexeev-Engel \cite{AE21} views a polarized K3 surface as a KSBA stable pair and provides a modular compactification for moduli spaces of polarized K3 surfaces.

It is a natural question to compare these  various compactifications. For example, Prof. Xu suggested in\cite[Question 10.9]{Xu21}  to study the relation between the compactifications for  moduli spaces of polarized K3 surfaces and the compact K-moduli of Fano $3$-fold pairs. The relation between the GIT compactification and arithmetic compactification has been explored for the moduli space of quartic K3 surfaces in a series of work by Laza-O'Grady \cite{LaO18ii} \cite{LO19}  \cite{LaO18}, where they propose the Hassett-Keel-Looijenga(HKL) program for the moduli space of quartic K3 surfaces, which is a highly nontrivial generalization of the well-known Hassett-Keel program for moduli space of stable curves $\overline{M}_g$. They conjecture that the HKL model defined in \cite{LO19} will interpolate the  GIT compactification and the Baily-Borel compactification for the moduli space of quartic K3 surfaces. This is verified by \cite{ADL22}.

In this paper, we study similar questions for the moduli space of  of K3 surfaces with a typical symmetry. That is, a K3 surface $Y$ with an antisymplectic involution $\tau: Y \rightarrow Y$. We call such a $(Y,\tau)$ a K3 pair.  What makes a K3 pair $(Y,\tau)$ very interesting is that it can be related to a del Pezzo surface pair  $(Y/\tau, C)$ where $C$ is  a  curve on $Y/\tau$  such that $C \sim -2K_{Y/\tau}$ and $C$ is isomorphic to $Y^\tau$ the fixed locus of $\tau$. Let $P_d$ be the moduli space of such del Pezzo surface pairs, then  $P_d$ can also be viewed as the moduli space of  K3 pairs due to the double covering construction. The rich geometric structure of $(Y,\tau)$ then allows us to study the compactification of $P_d$ from all the perspectives mentioned above. More precisely we are concerned with 
\begin{itemize}
    \item  GIT compactification $\overline{P}^{\rm GIT}_d$,
    \item  Baily-Borel compactification $P_d^\ast:=(O(T_d)\setminus \cD_{T_d})^\ast$,
    \item   compactifications $\overline{P}^{K}_{d}(c)$ via  K-stability and compactifications  $\overline{P}^{\mathrm{slc}}_{d}(c)$ via KSBA theory. 
\end{itemize}
and the connections of these compactifications. This paper mainly studies the connections between the first two compactifications. We view the work in this paper as the first step to establish the full HKL program for the moduli space $P_d$.

\subsection{Main results}
In this paper, we first obtain the class numbers of the Baily-Borel  compactifications of the moduli space of del Pezzo surface pairs of degree $d$. 
\begin{thm}\label{mthm1} Let $P_d^\ast$ be the Baily-Borel compactification of the moduli space of del Pezzo pairs of degree $d$, then the class number $\dim_\bQ \mathrm{Cl}(P_d^\ast)$ of $P_d^\ast$ is given by the following table
\begin{table}[ht]
    \centering
    \begin{tabular}{ |c|*{9}{c|} }
    \hline
     $d$ &   $8$ & $7$ & $6$  & $5$  & $4$  & $3$ & $2$ & $1$    \\ \hline
    $\dim_\bQ \mathrm{Cl}(P_d^\ast)$  &$3$ & $4$ & $4$  &  $4$ & $3$ & $3$  &  $3$  & $3$ \\  \hline 
     \end{tabular}
    \caption{Class numbers of $P_d^\ast$}
    \label{tb:Picard }
\end{table}
\end{thm}
\begin{rem}
    If $d=8$, we always refer to the situation where the del Pezzo surface is $Bl_p\bP^2$. 
    The case where the del Pezzo surface is $\bP^1\times \bP^1$ has been studied in the work of Laza-O'Grady \cite{LaO18} and Ascher-DeVleming-Liu  \cite{ADL21}.
\end{rem}

Our method proceeds by first constructing a suitable GIT model $\overline{P}^{\rm GIT}_d$ for $P_d$ and studying  GIT stability of points in the parameter space (see \S\ref{gitcons}). Unlike the case of plane sextic curves, i.e., $d=9$, the construction of $\overline{P}^{\rm GIT}_d$ for $d<9$ is somewhat subtle: for example, to avoid non-reductive GIT, we employ the geometry of the Hilbert scheme of points on $\bP^2$ for $d=7,8$.
Then we  compare $\overline{P}^{\rm GIT}_d$  with the Baily-Borel compactification  $P_d^\ast$.  The class groups of the GIT quotient spaces can be classically computed by \cite{KKV89}. By studying  GIT stability,  we construct a big open subset $U$ shared by $\overline{P}^{\rm GIT}_d$ and  $P_d^\ast$. 
Then we analyze the contribution of boundary divisors in  $P_d^\ast-U$ to the class group  $\Cl_\bQ(P_d^\ast)$ to prove the Theorem \ref{mthm1}. Our method also allows us to determine the generators of $\Cl_\bQ(P_d^\ast)$ (see Remark \ref{rem:generator} for precise details), which is the first step toward establishing the Hassett-Keel-Looijenga program for $P_d$. 

Next,  we focus on the case $d=8$ in more detail. Let $\cF$ be the locally symmetric variety associated to degree $8$ del Pezzo pairs $(Bl_p \bP^2,C)$, i.e., \[\cF = \widetilde{O}^+(T) \backslash \cD_{T},\quad T = U^2\oplus E_8(-1)\oplus E_7(-1)\oplus A_1 .\] There are two natural divisors $H_h$ and $H_u$, which are  called  reduced Heegner divisors (see Definition  \ref{Heegner}).
Let  $\lambda$ be the Hodge line bundle on $\cF$. 
Then we propose to study the space
\begin{equation}
    \cF(s):=\proj \big (\mathop{\bigoplus} \limits_{m \ge 0}H^0(\cF,ml_s\Delta(s)  ) \big ) ,\quad s\in [0,1]\cap \bQ
\end{equation}
where $\Delta(s):=\lambda+s(H_h+25H_u) $
 and $l_s\in \bZ_{>0}$ is the smallest positive integer so that $l_s \Delta(s)$ is Cartier divisor. Following \cite{LaO18},  we call it the HKL model for $\cF$. The first model $\cF(0)=\cF^\ast=P_8^\ast$ is exactly the Baily-Borel compactification for $\cF$.
 Since the GIT model for $P_8$ is morally speaking obtained from the Baily-Borel compactification by contracting the $H_u$ and $H_h$, it's expected that $\cF(1)$ gives back the $\overline{P}_8^{\rm GIT}$.
 To compute walls for $\cF(s)$ when $s\in [0,1] \cap \bQ$ varies, we give an arithmetic stratification for $\cF $ from towers of morphisms in the following diagram.
\subsubsection*{Stratification on the hyperelliptic divisor $H_h$:}
\begin{equation}\label{diah}
    \begin{tikzcd}
  \sh(A_7''')  \arrow[d,dashed] \arrow[dr] & &   &    &      \\
  \sh(A_7'') \arrow[r,dashed] &  \sh(A_6'')  \arrow[d,dashed] \arrow[dr]  &   &    &  \sh(A_1)=\cF    \\
    & \sh(A_6') \arrow[r,dashed]   & \sh(A_5') \arrow[d,dashed] \arrow[dr]      &         &  \sh(A_2)  \arrow[u]    \\
  &  & \sh(A_5) \arrow[r,dashed]  & \sh(A_4) \arrow[r]  & \sh(A_3) \arrow[u]     \\
   &  & & & \sh(D_4) \arrow[u]     \\
  & \sh(E_8) \arrow[r] & \sh(E_7) \arrow[r]   &  \sh(E_6) \arrow[r] &   \sh(D_5) \arrow[u]      \\ 
  &  & &  &   \sh(D_6) \arrow[u]       \\ 
 &  \sh(D_9') \arrow[r]  & \sh(D_8') \arrow[r]  &  \sh(D_7') \arrow[r,dashed]  \arrow[ur] &   \sh(D_7)  \arrow[u,dashed]     \\ 
\end{tikzcd}
\end{equation}
\subsubsection*{Stratification on the unigonal  divisor $H_u$:}
\begin{equation} \label{diau}
    \begin{tikzcd}
& \sh(U_4'') \arrow[r,dashed]  & \sh(U_3') \arrow[r]  & \sh(U_2')  \arrow[r,dashed] & \sh(U_1) \arrow[r]  & \cF
\end{tikzcd}
\end{equation}
In  the above diagram (\ref{diah}) and (\ref{diau}), $\sh(L)$ is the  locally symmetric variety of orthogonal type induced by a lattice $L$. See Definition \ref{defA}, Definition \ref{1stmv} etc for the precise choice of $L$.  For example, if $L=A_n$,  $\sh(A_n)$ is the locally symmetric variety given by the lattice $(E_7\oplus A_n)^\perp_{\II}$.  From the computation of the pullback of the $\bQ$-divisor $\Delta(s)$ along the tower $\sh(A)$ carried out in \S \ref{arithpred}, we predict that the walls are  located at  the image of $\sh(L) \rightarrow \cF$ for $L$ in the Table (\ref{diah}) and (\ref{diau}) above.  Typically, these $L$ are chosen via the principle of  ADE type lattice embedding $A_n \hookrightarrow A_{n+1}$, $D_n \hookrightarrow D_{n+1}$ etc, see \S \ref{arith} for more explanations. The dotted map indicates a modification occurs in the wall-crossing of $\cF(s)$ at this locus.

\begin{thm}\label{mthm2}
Notations as above, the class group of $P_8^\ast$ with $\bQ$-coefficient has a basis $\lambda, H_h,H_u$. Moreover, there is a Borcherds relation 
\[ 76\lambda = H_n+ 2H_h + 57 H_u. \]
Let $\varphi_h: \sh(A_1) \rightarrow P_8^\ast$ and $\varphi_u: \sh(U_1) \rightarrow P_8^\ast$ be the morphism  as in the diagram (\ref{diah}) and (\ref{diau}). Then their images are  $H_h$ and $H_u$ and we have the pullback formulas for the Heegner divisors
\[ \varphi_h^\ast H_h=-\lambda+H_{A_2},\ \ \varphi_u^\ast H_u=-\lambda .\]
Assume the section rings \[ R(s)=\mathop{\bigoplus} \limits_{m \ge 0 } H^0(P_8^\ast,m(\lambda+sH_h+25sH_u))\] are finitely generated for $s\in (0,1) \cap \bQ$ and that the birational contractibility of certain Heegner divisors on the arithmetic stratification holds(see $\widetilde{Z}$ in Theorem \ref{cbf} ).
The wall-crossing of the projective variety $\cF(s)=\proj(R(s))$ is described in Prediction \ref{pre1} and Prediction \ref{pre2}.
\end{thm}
\begin{rem}
    The coefficients for $H_h$ and $H_u$ are chosen  only related to divisorial contractions of K-moduli space wall-crossing,  see  \S \ref{degree8hkl} for more explanation.
\end{rem}
Based on the arithmetic strategy,  the  above arithmetic stratification predicts that the walls for $\cF(s)$ are given by 
\begin{equation*}
    \left\{\ \frac{1}{n} \, \mid \,  n=1,2,3,4, 6,8,10,12,16,25,27,28 ,31\ \right\}.
\end{equation*}
In the sequel paper \cite{psw2}, we will apply K-moduli theory to study the above  HKL model  and verify the above prediction, providing a new example of the HKL program worked out completely.  
Combined with \cite{LaO18}  and \cite{ADL21}, our work also completes the study of comparison for compactifications of the moduli space of hyperelliptic K3 surfaces.

We expect that our analysis on the arithmetic side will also help to find K-moduli  walls in other cases if the K-moduli space is birational to a ball quotient or a locally symmetric variety of orthogonal type. 

% \addtocontents{toc}{\SkipTocEntry}
\subsection*{Further directions}

One direction is to study the topology and intersection theory of moduli spaces of K3 surfaces  with an antisymplectic involution.  Combining wall-crossing results in \cite{psw2}, the birational map has an explicit resolution that transfers the arithmetic stratification described in the paper to the GIT strata. If a wall-crossing type formula for cohomology can be deduced from this explicit resolution of the birational period map, the cohomology of Baily-Borel compactification could be obtained. Note that a similar strategy to compute the cohomology of  moduli space of cubic $3$-folds and $4$-folds has been worked out in \cite{CGHL} \cite{Si}. It is also likely to compute the Chow ring of the moduli space $\cF$ based on our geometric description. A similar problem for moduli spaces of  quasi-polarized K3 surfaces of genus $2$  has been worked out in \cite{COP} and also the tautological property of Chow ring for these moduli spaces as in \cite{Si24}. 

Another direction is  to  investigate  possible arithmetic stratifications for certain toroidal compactifications of $\cF$ studied in \cite{AE23}, which may give a prediction of wall-crossings for KSBA compact moduli spaces of the surface  pairs $(X,cC)$ for $c>\frac{1}{2}$.  

\addtocontents{toc}{\SkipTocEntry}
\subsection*{Organization of the paper} In the paper, we work over $\bC$. In \S \ref{2-K3}, we set up the moduli space $P_d$ of del Pezzo pairs and discuss  the relation of $P_d$ to the moduli space of K3 surfaces with antisymplectic involution.  We then introduce the Baily-Borel compactification  $P_d^\ast$.  In \S  \ref{sec:GIT},  we construct GIT spaces $\overline{P}^{\mathrm{GIT}}_d$  and study the stability, which provides  the GIT (partial)compactification for $P_d$. 
In \S \ref{sect4} we give a classification of $O(\Lambda_d)$ orbit of $(-2)$-vectors for all $d$. 
As an application, we compute the class group $\pic(P_d^\ast)_\bQ$ of $P_d^\ast$. 
In \S \ref{sect:HKL model}, we set up HKL models $\cF(s)$ for $P_8^\ast$ and complete the arithmetic prediction of wall-crossing of $\cF(s)$.  In the two appendices \S \ref{Appa} and \S \ref{Appb},  we give the computational details on arithmetic stratification and  geometric degenerations corresponding to the arithmetic stratifications.  % we first study several towers of morphism between locally symmetric varieties. This gives an  arithmetic stratifications of the hyperelliptic divisor $H_h$ and the unigonal divisor $H_u$. Then based on the arithmetic principle, we compute the pullback of $\Delta(s)=\lambda+sH_h+25sH_u$ along  the towers and give the predictions for wall-crossings of $\cF(s)$.  In \S \ref{App}, we collect the computations of pullbacks of Heegner divisors in the arithmetic stratification.

\addtocontents{toc}{\SkipTocEntry}
\subsection*{Acknowledgment}
We thank  Zhiyuan Li and Yuchen Liu for  helpful conversations. We also thank Junyan Zhao for useful conversations. 
Some parts of this work were written when F.S visited SCMS and he would like to thank the hospitality.
F.S was partially supported by LMNS (the Laboratory of Mathematics for Nonlinear Science, Fudan University), the Fundamental Research Funds for the Central Universities,
NSFC (No. 12025106) and Shaanxi Fundamental Science Research Project for Mathematics and Physics (Grant No. 25JSZ007).
This project was also supported by NKRD Program of China (No.2020YFA0713200), NSFC Innovative Research Groups (Grant No.12121001) and General Program (No.12171090).

\section{Preliminaries }\label{2-K3}

\subsection{Basics on Locally Symmetric Varieties}

Let $\Lambda$ be an even lattice of signature $(2,n)$. The period domain associated to $\Lambda$ is defined as \[
    \cD_{\Lambda} \defeq \mybrace{z\in \bP(\Lambda \otimes \bC)}{z^2=0, z.\overline{z}>0}^+   
\]  where the superscript $\{\dots\}^+$ denotes a choice of a connected component.
%It is a Hermitian symmetric domain of Type IV and consists of two diffeomorphic components.  Denote by $ \cD_{\Lambda}^+$ one of them.  
Let $O^+(\Lambda) \le O(\Lambda)$ be the group of isometries preserving the component and 
\[ \widetilde{O}^+(\Lambda):=\{ \gamma\in O^+(\Lambda) |\  \gamma = \id \ \hbox{on}\ \Lambda^\ast/\Lambda\  \}. \]

\begin{thm}[Baily-Borel \cite{BB66}]
    For any finite index subgroup $\Gamma$ of $O^+(\Lambda)$, the orbit space $\Gamma \setminus \cD_{\Lambda} $  is a quasi-projective variety of dimension $n$.
\end{thm}
The quotient space $\Gamma \setminus \cD_{\Lambda} $  is called the locally symmetric space of orthogonal type, and it is the $\bC$-points of a connected Shimura variety of orthogonal type with some level structure. 
Indeed, Baily-Borel proved that  $\Gamma \setminus \cD_{\Lambda} $ admits a compactification $(\Gamma \setminus \cD_{\Lambda})^\ast $, known as the Baily-Borel compactification. 
The boundaries $(\Gamma \setminus \cD_{\Lambda})^\ast-\Gamma \setminus \cD_{\Lambda}$ are modular curves and points, which are determined by primitive isotropic sublattices of rank two and one in $ \Lambda$ respectively, see for example \cite[Section 3,4]{AE21} for more details.
Let $\breve{\cD}:=\{z\in \bP (\Lambda \otimes \bC) \mid z^2=0 \} $ be the compact dual of $\cD$, which is a smooth quadratic hypersurface. Let $\lambda$ be the descent of the line bundle $\cO_{\breve{\cD}}(1)|_{\cD}$ on $\Gamma \setminus \cD_{\Lambda}$. 
It is called the Hodge line bundle. 
Baily and Borel \cite{BB66} proved that 
\[ (\Gamma \setminus \cD_{\Lambda})^\ast \cong \proj \big ( \mathop{\bigoplus} \limits_{m \ge 0 } H^0(\Gamma \setminus \cD_{\Lambda}, \lambda^{m})\ \big)\]  
where the space $H^0(\Gamma \setminus \cD_{\Lambda}, \lambda^m)$ is exactly the space of modular forms of weight $m$ on $\Gamma \setminus \cD_{\Lambda}$.
% \subsubsection{Boundary of Baily-Borel compactification} Let \[ \breve{\cD}:=\{z\in \bP \Lambda \otimes \bC\  |\ z^2=0 \}  \] be the compact dual of $\cD$, which is a smooth quadratic hypersurface and let $\partial \cD$ be the boundary component of $\cD\subset \breve{\cD}$ in the Satake topology. A primitive isotropic sublattice $I \subset \Lambda$  will give a rational boundary component  \begin{equation*}  B_I= \begin{cases}  \{ z=x+\sqrt{-1}y \in \partial \cD | \{x,y\}_\bR \cong I\otimes \bR  \}    &  \rank (I)=2; \\ \{ \hbox{point} \}  &  \rank (I)=1. \end{cases} \end{equation*}Then there is natural identification  \[\cF_d^\ast= O(\Lambda)\setminus (\cD \sqcup \mathop{\cup} \limits_{I\ } B_I)\] where $I$ runs over $O(\Lambda)$-invariant primitive isotropic sublattice.  In particular, the boundaries \[ \cF_d^\ast-\cF_d=  O(\Lambda)\setminus \mathop{\cup}\limits_{I} B_I \cong  \mathop{\cup}\limits_{I} \mathrm{Stab}_{O(\Lambda)}(I) \setminus B_I\] consist of finitely many \begin{itemize}  \item points (also called $0$-cusps) for $\rank (I)=1$ and  \item  modular curves  (also called $1$-cusps) for $\rank (I)=2$. \end{itemize} From the perspective of Hodge theoretic degeneration of K3 surfaces, these points in $B_I$ correspond to Type II degeneration if $\rank (I)=2$ and to Type III degeneration if $\rank (I)=1$.See \cite[Section 3,4]{AE21} for more details. Let $\lambda$ be the line bundle obtained by the descent of the tautological line bundle $\cO_{\breve{\cD}}(-1)$. It is called the Hodge line bundle on the locally symmetric variety. 

\subsubsection{Heegner divisors}
$\Gamma \setminus \cD_{\Lambda} $  contains many natural locally symmetric subvarieties as effective cycles on  $\Gamma \setminus \cD_{\Lambda} $ which are very interesting objects in arithmetic geometry. In this paper, we are mainly interested in such codimension $1$ cycles. 
Let us give a precise definition. 

\begin{defn} Notation as above, denote by $\cD_v \defeq \{z\in  \cD_{\Lambda} \ |\   \langle z,v \rangle =0\}$ the hyperplane cut out by the vector $v\in \Lambda^\ast$ with $v^2<0$. 
Let $\Gamma\subset O^+(\Lambda)$ be a subgroup of finite index.
\begin{enumerate}

\item Let $v\in \Lambda$ be primitive with $v^2<0$. 
A Heegner divisor $H_{v,\Gamma}$ associated to $v$ on $\Gamma \setminus \cD_{\Lambda} $ is defined by \[
H_{v,\Gamma} := \Gamma  \setminus \sum_{g \in \Gamma } \cD_{g\cdot v} 
 .\]
 We abbreviate it as $H_v$ when $\Gamma $ is specified.

\item  Let $\beta \in G(\Lambda), m\in \bQ_{<0} $. 
If $\Gamma \subset \widetilde{O}^+(\Lambda)$ is a finite index subgroup, then \(
\sum_{\substack{v\in\Lambda^\ast,\, v^2=2m  \\ v\in \beta + \Lambda
  }} \cD_v .
\) is a $\Gamma $-invariant collection of hyperplanes with finitely many orbits. 
A Heegner divisor with discriminant $(\beta,m)$ on $\Gamma \setminus \cD_{\Lambda} $ is defined by
\begin{equation*}
  H_{\beta,m}:=\Gamma  \setminus \sum_{\substack{v\in\Lambda^\ast,\, v^2=2m  \\ v\in \beta + \Lambda
  }} \cD_v .
\end{equation*} 
The primitive Heegner divisor $P_{\beta,m}$ is defined similarly,
with $v$ required to be primitive in $\Lambda^\ast$.
\end{enumerate}
\end{defn} 
In general, a Heegner divisor $H_{\beta,m}$ can be non-reduced and reducible, while $H_v$ is always irreducible.  
Indeed, $H_{\beta,m}$ can be written as a sum of Heegner  divisors associated to vectors
\[H_{\beta,m}= \mathop{\sum}  \limits_{v} a_vH_v. \]
By \cite{BM19}, there is a relation (referred to as the triangular relation in the literature)  
   \begin{equation}\label{trirel}
       H_{\beta,m}=\mathop{\sum} \limits_{\substack{r\in \bZ_{>0} \\ r^2| m }}  \mathop{\sum} \limits_{\substack{\delta\in G(\Lambda)\\ r\delta=\beta}}  P_{\delta, \frac{m}{r^2}}
   \end{equation}
   between the Heegner divisors and the primitive Heegner divisors. 

\begin{lem}[Eichler's criterion] \label{Eichler} If $\Lambda$ is an even lattice containing $U^{\oplus 2}$, then two nonzero primitive elements $v,w\in \Lambda$ are in the same $\widetilde{O}^+(\Lambda)$-orbit if and only if $v^2=w^2,\ v^\ast=w^\ast $.
Here $v^* \defeq [v/ \divi(v)] \in G(\Lambda)$ where $\divi(v) \defeq \min\{\langle v,w\rangle \in \bZ_+ \mid  w\in \Lambda\}$.
\end{lem}
\begin{proof}
See  \cite[Proposition 3.3]{GHS09}.
\end{proof}
Define $\sh(\Lambda):=\widetilde{O}^+(\Lambda) \setminus \cD_\Lambda$.   Using Eichler's criterion,  Bruinier--M\"{o}ller \cite[Lemma 4.3]{BM19} proved that  the reduced part of the primitive Heegner divisor  $P_{\beta,m}$ on $\sh(\Lambda)$ is irreducible if $\Lambda$ contains two copies of the hyperbolic lattice. 
Note that any $\ell \in \Lambda^*$ is primitive if and only if $\ell = v/ \divi(v)$ for some unique primitive $v\in \Lambda$.
Thus, for $v\in \Lambda$ a primitive vector with $ v^\ast = \beta  , \ q(v^\ast)=m$, we have 
\begin{equation*}
  P_{\beta,m}=H_v  
\end{equation*}
as divisors on $\sh(\Lambda)$.   
%Note that there is a natural finite morphism\begin{equation*} \sh(\Lambda) \rightarrow \Gamma \setminus \cD_{\Lambda}.\end{equation*}for any $\Gamma$ containing $\widetilde{O}^+(\Lambda)$ as a finite index subgroup. 

%In particular, if $\Lambda$ is a $2$-elementary lattice, then by \cite[Section 4.3]{Ma}, the finite morphism $\sh(\Lambda) \rightarrow O(\Lambda) \setminus \cD_{\Lambda}$ has degree $|O(G(\Lambda))|$. 

\subsubsection{Modular forms on locally symmetric varieties }
We now review the  study  of the  Picard group $\pic(\sh(\Lambda))_\bQ$ via vector-valued modular forms by the work of Borcherds \cite{Bor}, Bruinier \cite{bruinier} and Bergeron-Li-Millson-Moeglin \cite{BLMM}. 
Let $\Lambda$ be an even lattice of signature $(2,n)$ containing $U^{\oplus 2}$ as a summand. 
Let $\bC[G(\Lambda)]$ be the complex vector space spanned by the discriminant group $G(\Lambda)$ and denote $\MP_2(\bZ)$ the metaplectic group. 
Then there is a so-called Weil representation  \[ \rho_\Lambda: \MP_2(\bZ) \rightarrow \GL(\bC[G(\Lambda)]) .\]
Denote the $\acup(\rho_\Lambda,k)$ space of almost cuspidal $\bC[G(\Lambda)]$-valued modular forms of weight $k$, that is, the space of vector-valued holomorphic functions
\begin{equation*}
    f: \BH \rightarrow \bC[G(\Lambda)]
\end{equation*}
satisfying certain modularity axioms with respect to $\rho_\Lambda$. 

\begin{thm}[\cite{Bor}, \cite{bruinier} and \cite{BLMM}] \label{modularform1}
Notation as above, then the Picard group $ \Pic_\bQ(\sh(\Lambda))$ is generated by Heegner divisors  and  there is an isomorphism of $\bQ$-vector spaces
\begin{equation*}
    \Pic_{\bQ}(\sh(\Lambda)) \rightarrow \acup(\rho_\Lambda,1+\frac{n}{2})^\vee_{\bQ}
\end{equation*}
by sending the Heegner divisor $H_{\xi,m}$ to the linear map 
\begin{equation*}
    f= \mathop{\sum} \limits_{m \in \bQ, \xi \in  G(\Lambda) } c_{m,\xi}\cdot  q^m e_\xi \ \mapsto \  c_{m,\xi}
\end{equation*}
where $\{e_\xi\}$ is the standard basis of $\bC[G(\Lambda)]$ given by the elements $\xi$ of  the discriminant group $G(\Lambda)$.
\end{thm}

In practice,  Raum (\cite{Rau16}) gives an algorithm for computation based on Jacobi forms for definite lattices, which can be translated into a Sage package by Williams \cite{Williamsa} in $\mathbf{Sage}$ \cite{SageDevelopers2022}. Theorem \ref{modularform1} will be used  in \S \ref{subsec: h arith strata}, \S \ref{subsec:u arithstra}, \S \ref{arithpred} and Appendix \ref{Appb} to conduct arithmetic computation on the restriction of Heegner divisors on the locally symmetric variety from moduli space of del Pezzo pairs.

\begin{rem}
     If $\widetilde{O}^+(\Lambda)$ is replaced by another arithmetic group $\Gamma$, then the above theorems may fail for $\Gamma \setminus \cD_\Lambda$. 
\end{rem}

\subsection{Moduli of del Pezzo pairs and K3 surfaces  with antisymplectic involution}
\begin{defn}\label{delpezzopair}
Let $X$ be a smooth projective surface with ample $-K_X$.  
The pair $ (X,C)$ is called a smooth del Pezzo pair of degree $d$ if $(-K_X)^2=d$ and $C \sim -2K_X$ is a smooth curve. 
Let $\pi: (X,C) \rightarrow B$ be a flat proper smooth morphism over a base scheme $B$ and $C\subset X$ a divisor. 
We say $\pi$ is a family of del Pezzo pairs of degree $d$ if
\begin{enumerate}
    \item $-K_{X/B}$ is $\pi$-ample,
    \item $C\sim-2K_{X/B}$ is a relative effective Cartier divisor, flat over $B$,
    \item and  $(X_t,C_t)$ is a smooth del Pezzo pair of degree $d$ for any closed point $t\in B$.
\end{enumerate}  
\end{defn}

Consider the moduli stack $\mathcal{P}_d:\ \csch(\bC) \rightarrow \cset $ defined by 
\begin{equation*}
\begin{split}
    \mathcal{P}_d(B)=\{ \pi: (X,C) \rightarrow B\ |\ \pi \  \hbox{is a family of smooth del Pezzo pair of degree } d\  \} /\sim. 
\end{split}
\end{equation*}
Then $\cP_d$ is an open substack of one irreducible component of the KSBA moduli stack of slc pair $(X,(\frac{1}{2}+\epsilon) C)$ for $0< \epsilon \ll 1$, which is a separated Deligne-Mumford stack.
We denote its coarse moduli space as $P_d$.

Given a smooth del Pezzo pair $(X,C)$ of  degree $d$, one can construct a K3 surface $Y$ as the double covering
\[\mu: Y \rightarrow X \]
branched along the curve $C\sim -2K_X$. 
Such a K3 surface has a natural involution $\tau: Y \rightarrow Y$ which is antisymplectic, whose fixed locus is just the smooth curve $C$ of genus $d+1$ by the adjunction formula. 

\begin{defn}
An involution $\tau: Y \rightarrow Y$ for a K3 surface $Y$ is called \emph{antisymplectic} if $\tau^\ast \sigma=-\sigma$, where $\sigma$ is the holomorphic volume form and  $H^{2,0}(Y)=\bC \cdot\sigma$.
\end{defn}
 
%In this subsection we recall the construction the moduli space of such K3 surfaces with certain type involution and show it will provide arithmetic compactification of $P_d$. The del Pezzo pairs are closely related to K3 surfaces with antisymplectic involution:

To construct the moduli space of K3 surfaces with antisymplectic involutions, we need to fix the type of the induced action on K3 lattice  $L=U^3\oplus E_8^2 $.
We recall the definition of $\rho$-markable K3 with non-symplectic action studied in \cite{AEH}.
 
\begin{defn}
   Fix $\rho \in O(L)$ of order $2$.  
   A $\rho$-markable K3 is a pair $(Y,\tau)$, where $Y$ is a smooth K3 surface and $\tau: Y \rightarrow Y$ is an antisymplectic involution such that there is an isometry $\phi: H^2(Y,\bZ) \rightarrow L$ satisfying  
    \[  \tau^\ast = \phi^{-1} \circ \rho  \circ \phi   . \]
\end{defn}
Denote by $\cF_\rho$ the moduli stack of $\rho$-markable K3 surfaces with antisymplectic involutions. 
%That is,  \begin{equation*}  \begin{split} \cF_\rho(B)=\big\{ \pi: (\fX,\tau_B) \rightarrow{} B \big\}/\sim. \end{split} \end{equation*}
Then $\cF_\rho$ is a Deligne-Mumford stack with coarse moduli space $F_\rho$. 
Let $N,T$ be the $\pm 1$-eigenspaces of $\rho$ respectively which satisfy $T = N^\perp$.
Set $ \cD_\rho:=\cD_T  \subset \cD_L$
the subdomain of the period space  $\cD_L$ of complex K3 surfaces. 
Then  $\cD_\rho$ is a Type IV Hermitian domain.
It is shown in \cite[Theorem 2.9]{AEH} and \cite[Section 2C]{AE23} that the group
\begin{equation*}
    \Gamma_\rho:=\{ g\in O^+(L)\mid g\circ \rho =\rho \circ g \}
\end{equation*} surjects onto $O^+(T)$ which acts on $\cD_\rho$, and the period map 
\begin{equation} \label{eq:fine moduli}
   \mathfrak{p}_\rho: F_\rho \rightarrow  \Gamma_\rho \setminus \cD_\rho  \cong O^+(T) \setminus \cD_T
\end{equation}
is an open immersion with the image $\Gamma_\rho \setminus (\cD_\rho - \bigcup_{\delta\in T,\delta^2=-2} \delta^\perp) $.

%\begin{rem} Note that for a $\rho$-markable K3 $(X,\varphi)$, the Neron-Severi group $\NS(X)$ of $X$ always contain $N$ the invariant sublattice of $\rho$. Generically,   $\NS(X)\cong N$.  For those $(X,\varphi)$ such that $\rank(\NS(X)/N) \ge 1$, they may form  Heegner divisors on   $ \Gamma_\rho \setminus \cD_\rho$. \end{rem}

The type $\rho$ is uniquely determined by the invariant (equivalently, anti-invariant) lattice $N \subset L$(resp, $T$) which was completely classified by \cite[Section 3.6]{Nikulin} (see also \cite{AN06}).

\begin{defn}
An even indefinite lattice $\invLattice$ is called a \emph{$2$-elementary lattice} if the discriminant group has the form \[G(\invLattice):=\invLattice^\ast/\invLattice \cong (\bZ/2\bZ)^a.\] 
\end{defn}

\begin{pro}[\cite{Nikulin}] A $2$-elementary lattice $\invLattice$ is uniquely determined by the main invariants $(r,a,\delta)$, where
\begin{enumerate}
    \item $r$ is  the rank of  $\invLattice$;
    \item  $a\in \bN $ the length of the discriminant group $G(N)$;
    \item $\delta: =  \begin{cases}
 0 & \text{if } q(G(\invLattice))\subset \bZ,\\
1 & \text{otherwise} .
\end{cases}$
\end{enumerate}
\end{pro}

Let $(Y,\tau)$ be a K3 surface with antisymplectic involution $\tau$. 
The quotient surface $X\defeq Y/\tau$ is smooth.
Then the invariant lattice $N_\tau$ of $\tau$ is 2-elementary and is uniquely determined by its invariants $(r, a,\delta)$. 
Let 
\begin{equation}\label{genus}
    g=\frac{22-r-a}{2},\ \ k=\frac{r-a}{2}.
\end{equation}
The triple $(g,k,\delta)$ is an alternative set of invariants of $N_\tau$. There are three cases:
\begin{enumerate}
    \item the fixed locus $Y^\tau=C_g+R_1+\cdots +R_k$ is a sum of an irreducible  nonsingular curve $C_g$ and $k$ irreducible nonsingular rational curves $R_1,\cdots R_k$. The surface $X$ is rational.
    \item $ Y^\tau= C_0 \sqcup C_1$ is the union of two elliptic curves. Then $(r,a,\delta) = (10,8,0)$ and $ (g, k, \delta) = (2, 1, 0)$. The surface $X$ is rational elliptic.
    \item $ Y^\tau= \emptyset$. Then $(r,a,\delta) = (10,10,0)$ and $(g,k,\delta) = (1,0,0)$. The surface $X$ is Enriques.
\end{enumerate}

From the above double covering construction, the moduli space $P_d$ of del Pezzo pairs of degree $d$ is closely related to the moduli space of $\rho$-markable K3 pairs for a certain involution $\rho$ . 
Here the involution $\rho= \rho(d)$  has invariant $(r=10-d,a=10-d,\delta =1 )$ and the corresponding anti-invariant lattices $T_d$ are classified in the Table \ref{table:anti-inv lattice} (see \cite{Yoshikawa20}).
\begin{table}[ht]
    \centering
    \begin{tabular}{ |c | c | c|c |c| c|c|}
    \hline
     $d$ &  $  a  $ & $ \delta$ &  $ T_d\ \ $ & Discriminant group $G(T_d)$ & \text{ complements in $I\!I_{2,26}$}  \\ \hline
     1 & $9$ & $1$  & $U^2 \oplus A_1^9$    & $(\bZ/2\bZ)^9$ & $E_7 \oplus A_1^8$ \\ \hline
     2 & $ 8$ & $1$  & $U^2 \oplus D_4 \oplus A_1^6$  &     $(\bZ/2\bZ)^8$ & $E_7 \oplus A_1^7$\\ \hline
    3 & $ 7$ & $1$  & $U^2\oplus D_6 \oplus A_1^5$ &   $(\bZ/2\bZ)^7$   & $E_7 \oplus A_1^6$\\ \hline
     4 & $ 6$ & $1$  & $U^2\oplus E_7 \oplus A_1^5$  &  $(\bZ/2\bZ)^6$ & $E_7 \oplus A_1^5$\\ \hline
      5 & $ 5$ & $1$  & $U^2\oplus E_8 \oplus A_1^5$ &    $(\bZ/2\bZ)^5$ & $E_7 \oplus A_1^4$ \\ \hline
      6 & $ 4$ & $1$  & $U^2\oplus D_4 \oplus E_8 \oplus A_1^2$ &     $(\bZ/2\bZ)^4$  & $E_7 \oplus A_1^3$\\ \hline
      7 & $ 3$ & $1$  & $U^2\oplus D_6 \oplus E_8 \oplus A_1$  & $(\bZ/2\bZ)^3$  & $E_7 \oplus A_1^2$ \\ \hline
      8 & $ 2$ & $1$  & $U^2\oplus E_7 \oplus E_8 \oplus A_1$   & $ (\bZ/2\bZ)^2$   & $E_7 \oplus A_1$ \\ \hline 
      9 & $ 1$ & $1$  & $U^2 \oplus E_8^2 \oplus A_1$   & $ \bZ/2\bZ$   & $E_7$ \\ \hline 
      % 8\,$'$ & $ 2$ & $0$  & $U^2\oplus D_{16}$   & $ (\bZ/2\bZ)^2$   & $D_8$ \\ \hline
  \end{tabular} 
%   \medskip
  \caption{List of anti-invariant lattices for $\delta=1$}
  \label{table:anti-inv lattice}
\end{table}

\begin{rem}
The only possibility for $\delta =0$ is the case $\rho(d)$ comes from the double cover of $\bP^1 \times \bP^1$, another del Pezzo surface of degree $8$, and the corresponding $T\cong U^2\oplus D_{16}$.

\end{rem}

To make the notation consistent, we denote the corresponding moduli space by $F_{\rho(d)}$ and the monodromy group $\Gamma_d:=\Gamma_{\rho(d)}$.
We have the following:
\begin{pro} \label{prop:moduli meaning}
  For $1\leq d \leq 9$, $\rho= \rho(d) \in O(L)$ as above, there is an open immersion
    \begin{equation*}
      \mathfrak{p}_d:   P_d  \rightarrow F_{\rho(d)}    \end{equation*}
    defined by sending the isomorphism classes of a del Pezzo pair $(X,C)$ of degree $d$ to the  $O(\Lambda)$-orbit of the period point  of its associated K3 surface with anti-involution $(Y,\tau)$. 
\end{pro}
\begin{proof}
This follows from \cite[Theorem 4.2]{AEH} directly.     
\end{proof}

\section{GIT compactification for 
\texorpdfstring{$P_d$}{P_d}}\label{sec:GIT}

In this section, we construct various GIT (or partial) compactifications $\overline{P}_d ^{\rm GIT}$  for the moduli space $P_d$ of smooth log del Pezzo pairs. 
% The proof of the claim is almost the local computation. If $Z=p_1+p_2$ consists of two distinct points, then locally  $\fm_{p_{i}}=\langle x,y \rangle \subset \bC[x,y]$ and thus \[\dim _\bC \cO_X/I_Z^2=\dim _\bC (\bC[x,y]/\langle x,y \rangle ^2)^{\oplus 2}=6 .\]
% If $Z$ is a closed subscheme of length $2$ supported at a point $p\in X$, then locally  $I_Z=\langle x,y \rangle \subset \bC[x,y]$. Then $\fm_p^2 \subset I_Z\subset \fm_p$. Thus we may assume $I_Z=\langle x^2, xy,y^2,x+ay \rangle, a\neq 0$ and then \begin{equation*}
%     \begin{split}
%         I_Z^2=\langle x^4, x^2y^2,y^4,(x+ay)^2, x^3y,xy^3,x^3+ax^2y,xy^2+ay^3 \rangle
%     \end{split}
% \end{equation*}
% and thus a direct calculation will show $\dim_\bC \bC[x,y]/( I_Z^2)=7$. This finish the proof of claim. 
They will be constructed separately because the relevant geometric structures differ with the degree.

We provide the stability analysis for general elements, and a fundamental tool is the Hilbert-Mumford numerical criterion. 
Recall that for a polarized projective variety $(Z,L)$ under the action of a reductive group $G$, there is a numerical function $\mu^L$ (abbreviated as $\mu$ when $L$ is fixed) called the Hilbert-Mumford index on the space of one-parameter subgroups ($1$-PS) of $G$
\[ \mu(z;\cdot):\   \Hom  (\bG_m,G) \rightarrow \bR \]
for any given point $z\in Z$.  
By \cite[Section 2.1]{MFK94}, $z\in Z$ is semistable if and only if $\mu(z;\lambda) \ge 0$ for any $1$-PS $\lambda$. 
$z\in Z$ is stable if and only if $\mu(z;\lambda) > 0$ for any nontrivial $1$-PS $\lambda$. 
We will apply this numerical criterion to study the GIT space constructed in the following subsections.

\subsection{GIT model for degree $8,7$ del Pezzo pairs}  \label{gitcons}
%Let $\Delta \subset \bP^2 \times \bP^2$ be the diagonal and $p_i:\bP^2 \times \bP^2 \rightarrow \bP^2$ be the projection. 
Let $X:=\bP^2$ and $\cZ \subset X^{[n]} \times X$ be the universal closed subscheme of length $n$. Denote $p_1: X^{[n]} \times X \rightarrow X^{[n]} $ and $p_2: X^{[n]} \times X \rightarrow X $ the two projection morphisms, then there is 
\[ 0 \rightarrow \cI_{\cZ}^2\otimes p_2^\ast \cO_X(6)  \rightarrow p_2^\ast \cO_X(6) \rightarrow (\cO_{X^{[n]}\times X}/ \cI_{\cZ}^2)\otimes p_2^\ast \cO_X(6) \rightarrow 0. \]
Let $\cE_n:=p_{1\ast }(\cI_{\cZ}^2\otimes p_2^\ast \cO_X(6))$. 
 If $Z\subset X$ is any closed subscheme  of length $n\le 2$, then  $I_Z^2$ defines a closed subscheme of the same length $\le 7$ for fixed $n$ by direct computation. (This fails when $n\geq3$.)
 Note that $ \cO_X(6)$ is $6$-very ample line bundle on $X$ since $\cO_X(1)$ is very ample. Thus, for any $0$-dimensional closed subscheme $Z\subset X$ of length $\le 7$, one has $ H^0(X,\cO_X(6))\twoheadrightarrow H^0(Z,\cO_X(6)|_Z)$. 
In particular, one then has \[
h^0(X,\cO_X(6)) = h^0(X,I_Z^2\otimes \cO_X(6))+ h^0(X,\cO_X/I_Z^2(6))
\]
and $h^0(X,\cO_X/I_Z^2)$ (hence also $h^0(X,I_Z^2\otimes \cO_X(6))$) remains the same for any $[Z]\in X^{[n]}$.
By Grauert's theorem, this implies that  $\cE_n$ is a vector bundle, and its fiber at the closed subscheme $Z\subset X$  is the vector space  $H^0(X,I_Z^2\otimes \cO_X(6))$, the plane sextics with vanishing order at least $2$ along $Z$. 
Set \[\pi: \ \bP \cE_n:=\proj_{\cO}(\sym^\ast(\cE_n^\vee) ) \rightarrow X^{[n]}  \]
be the projective bundle whose fiber $Z\in X^{[n]}$ is  the space of the one–dimensional subspaces of $\cE_n|_Z=H^0(X,I_Z^2\otimes \cO_X(6))$(cf.~ \cite[Appendix A]{Lazarsfeld1}). 
Let $L_t\coloneqq (\cO_{X^{[n]}}(1)-\epsilon E) +t\cdot \cO_{\bP \cE_n}(1), 0<\epsilon\ll 1$, where $E= 0$ for $n=1$, and $E$ is the exceptional divisor of the Hilbert-Chow morphism for $n=2$, and $\cO_{X^{[n]}}(1)$ is the pullback of the natural polarization from the symmetric product along the Hilbert-Chow morphism.
As the nef cone of $(\PP^2)^{[n]}$ is bounded by $\cO(1)$ and $(n-1)\cdot\mathcal{O}_{X^{[n]}} (1)-\frac{1}{2}E$ (see \cite[Corollary 3.2]{ABC+_minimal_2013}), one see $(\cO_{X^{[n]}}(1)-\epsilon E) $ is ample on $X^{[n]}$.
Since $\cE_n^\vee$ is a quotient of the trivial vector bundle $p_{1\ast }( p_2^\ast \cO_X(6))$, $\cO_{\bP \cE_n} (1)$  is nef by \cite[Proposition 6.1.2]{Lazarsfeld2}.  
We conclude that $L_t$ is ample for all $t>0$ by Nakai-Moishezon's ampleness criterion.  

Now, we consider the VGIT quotient
\[ V_n(t)\defeq \bP \cE_n \q_{L_t}  \SL(3), \] 
which will provide the GIT compactification of $F_{\rho(d)}$ for certain range of $t$.
Before proceeding, we make some preparation for the detailed stability analysis. Since the stability of elements in the same orbit remains the same, it suffices to check the stability of the representative $([Z],f) \in \PP \cE_n$ where \[ 
\Supp [Z]=
\begin{cases}
   p= [0,0,1], & n=1, \text{~or~} n=2 \text{~and~} [Z] \text{~non-reduced;}\\
   p=[0,0,1],\,q=[1,0,0]  & n=2 \text{~and~} [Z] \text{~reduced~.}
\end{cases}
\] 

\begin{lem}
    The stability of $([Z],f)$ is completely determined by the Hilbert-Mumford weight for 1-PS preserving $\Supp[Z]$.
\end{lem}

\begin{proof}
Let $([Z_0],f_0) = \lim_{s\to       0}\lambda(s)\cdot([Z],f)$ be the limit point.
For the $n=1$ case and the $n=2$, $[Z]$ nonreduced case, if $[Z_0]\neq [Z]$, there exists $g\in \SL(3)$ such that $g\cdot [Z_0] =[Z]$.
Note that \[\begin{split}
    \mu([Z],f;\lambda) &= \mu([Z_0],f_0;\lambda)  \\
    &=\mu([Z],g\cdot f_0; g\lambda g^{-1}).
\end{split}\]
The polynomial $g\cdot f_0$ is S-equivalent to $f$ and $\lim_{s\to 0} (g\lambda \inverse{g})\cdot [Z] = [Z]$.
So we may assume the limit of $\lambda\cdot [Z]$ is still $[Z]$.
Write $$\lambda(s) = P\begin{pmatrix}
    s^{a_1} & & \\
    & s^{a_2} & \\
    & & s^{a_3}
\end{pmatrix} \inverse{P} \text{~and~} P = (P_1,\,P_2, \, P_3)\in \GL(3).$$ where $a_1+a_2+a_3=0$.
Suppose $P\cdot v =\sum_{i=1}^3 v_iP_i= (0,0,1)^t$.
Then \[
\begin{bmatrix}
    0 \\ 0\\ 1 
\end{bmatrix}
= \lim_{s\to 0} \lambda(s)\cdot \begin{bmatrix}
    0 \\ 0\\ 1 
\end{bmatrix} = \lim_{s\to 0 }\sum_{i=1}^3 s^{a_i} v_i P_i\in \PP^2.
\] 
The limit is equal to $v_{i_0}P_{i_0}$ when the minimum $a_{i_0}$ of $\{a_i\}$ is unique, or equal to $v_{i_1}P_{i_1}+ v_{i_2}P_{i_2}$ when $a_{i_0} = a_{i_1}<0$.
By the linear independence of $P_i$, one knows that all $v_i=0$ except $v_{i_0}$ for the first case, and $v_{i_1}$ and $v_{i_2}$ for the second case.
In both cases, one can directly check that $\lambda(s)$ preserves $[0,0,1]$.

For the case  $\Supp [Z]= \{[1,0,0],[0,0,1]\}$, if $([Z],f)$ is destabilized by 1-PS with the limit lying in the non-reduced locus, it is S-equivalent to a non-reduced point, and we've already treated this case above.
So we only consider whether $([Z],f)$ is destabilized by 1-PS with the limit still being reduced, and the same argument applies.
\end{proof}

\subsubsection{The case $d=8$  stability analysis}
In degree $d=8$ case, let $Bl_p \bP^2 \rightarrow \bP^2$ be the blowup at $p=[0,0,1]\in \bP^2$, then $f$ takes the form
\begin{equation}\label{nodalequ}
    z^4f_2(x,y)+z^3f_3(x,y)+\cdots +f_6(x,y).
\end{equation}
For the $d=8$ case, one gets the following results.
\begin{pro}\label{pro:GIT8anly}
When $d=8$ (i.e. $n=1$), the following holds:
    \begin{enumerate}
        \item for $0<t<\frac{1}{3}$, $\sstable{\bP \cE}=\emptyset$;
        \item for $\frac{1}{3}<t< \frac{5}{12}$, $V_1(t)$ compactifies $P_8$ and has Picard rank one.
    \end{enumerate}
\end{pro}

\begin{proof}
    It suffices to compute the Hilbert-Mumford weight of $([0,0,1],f)$ for the diagonalized 1-PS $\lambda(s)=(s^{-1},s^{-1},s^2)$ with respect to the polarization $L_t$. 
    
    When $t<\frac{1}{3}$, we see that 
    \[
    \begin{split}
            \mu_t(\lambda) &= -2+ t\cdot \mu(f,\lambda)\leq 6t-2< 0 .
    \end{split}
    \]
    Hence, $\sstable{\bP \cE}$ is empty for $0<t<\frac{1}{3}$, $\sstable{\bP \cE}=\emptyset$.
    
    When $0<t<\frac{2}{3}$, if  $f$ has no $z^4f_2(x,y)$ term, then 
    $\mu_t(\lambda)  \leq -2+ 3t<0$.
    So we only consider those $f$ with nonzero $z^4 f_2(x,y)$ terms.
    Suppose $f_2$ has rank one, which can be assumed to be $x^2$ up to the coordinate change of $x$ and $y$. 
    Then we have \[
    \mu_t(\lambda'(s)) \leq -5 + 12t\  \text{~for~} \ \lambda'(s)= (s^{-4},s^{-1},s^5)
    \]
    where the equality holds if there exists a nonzero $z^3 y^3$ term.
    Thus, we know that such a pair $(p,f)\in \bP \cE$ is unstable until $t\geq \frac{5}{12}$.
    
    When $\frac{1}{3}\leq t< \frac{5}{12}$, we only need to consider $f$ with $\rk f_2 =2$, i.e. those $(p,f)$ with $f$ having at worst $A_1$-singularity at $p$.
    Suppose $f$ is destabilized by some 1-PS in $\stab_{[0,0,1]}(\SL(3))$, then up to the coordinate change of $x$ and $y$, we may assume the 1-PS takes the form \[
    \begin{pmatrix}
        s^a &     & \\
            & s^b & \\
          \text{const}\cdot (s^a-s^{-a-b}) & \text{const}\cdot (s^b-s^{-a-b}) & s^{-a-b}
    \end{pmatrix} , \; a\ge b.
    \]
    Then there are three possible cases:
    \begin{enumerate}
        \item $a\ge b \ge -a-b$ or $a\ge -a-b \ge b$: if $f$ has at worst $A_1$ singularity at $[1,0,0]$, then \[
        \mu_t = a+b +t \cdot \mu(f) \geq (a+b)+\frac{1}{3}\cdot 3a = 2a+b\geq 0
\]

        \item $-a-b\ge a \ge b$: since $f$ has rank $2$ at $[0,0,1]$, we have\[
        \mu_t = a+b +t \cdot \mu(f) \geq (a+b)+\frac{1}{3}\cdot [4(-a-b)+a+b] = 0
\]
    \end{enumerate}

 Combined together, one sees that when $\frac{1}{3}< t< \frac{5}{12}$, the $(p,f)$ is semistable for those $f$ that have at worst  $A_1$-singularities.  
 In this case, the corresponding surface pair $(\mathbb{F}_1,C)$ admits at worst $A_1$-singularities on the boundary curve.
 Thus, $V_1(t)$ compactifies $P_8$.
 Denote $W$ the locus of $\rk f_2 <2$, which is a divisor in $\bP\cE$, then $W\subset \bP \cE^{\rm us}_t$ and $\mathrm{CH}_{25}(\bP\cE^{\rm us}_t) = \bZ \cdot W$ by the above analysis.
(This corresponds to the hyperelliptic divisor $H_h$ in $F_{\rho(8)}$.)

We have the following exact sequence \begin{equation}\label{eq:pic8}
    \mathrm{CH}_{25}(\bP\cE^{\rm us}_t)\xrightarrow{\psi} \Pic  (\bP\cE)  \to \Pic (\bP\cE^{\rm ss}_t )\to 0
\end{equation}

Since $V_1(t) = \bP\cE^{\rm ss}\q_{L_t}\SL(3)$ is nontrivial and projective, by \cite[Proposition 5.1, \S 2.1]{KKV89}, we have \begin{equation} \label{eq:pic8'}
\Pic V_1(t) \hookrightarrow \Pic_{\SL(3)}\bP\cE^{\rm ss} \hookrightarrow \Pic \bP\cE^{\rm ss} 
\end{equation}
as $\SL(3)$ has only trivial characters.
Equation \eqref{eq:pic8} shows that $ \bP\cE^{\rm ss}$ has Picard rank $\leq 1$ as the image of $\psi $ is nontrivial and $\rk \Pic \bP\cE =2$.
One then concludes that $\rk \Pic V_1(t) =1 $.

\end{proof}

Therefore, we have a reductive GIT space which (partially) compactifies $P_8$
\begin{equation}  \label{gitd8}
    \overline{P}_8^{\rm GIT} \defeq V_1(t),\quad \frac{1}{3}<t<\frac{1}{3}+\epsilon.
\end{equation}

We remark that $\overline{P}_8^{\rm GIT} \cong \bP V\q G$ is defined in \cite[Section 5.1]{psw2} by the reduction procedure above.

\subsubsection{The case $d=7$  stability analysis}
For the $d=7$ case, consider the following commuting diagram:
\[
\begin{tikzcd}
Bl_\Delta (\bP^2)^2 \arrow[r, "h'"] \arrow[d, "\pi'"']   & (\bP^2)^2 \arrow[d, "\pi"] \\
{(\bP^2)^{[2]}} \arrow[r, "h"']                        & (\bP^2)^{(2)}             
\end{tikzcd}
\]
where $h$ is the Hilbert-Chow morphism (which is a crepant resolution) and $\pi$ is the finite quotient morphism.
Since $\pi'$ is also a finite quotient map, to check the stability of $([Z],f)$ with respect to $L_t$, one can instead check the stability of points in $(\inverse{(\pi')}([Z]),f)$ with respect to ${\pi'}^*L_t$.
Let $E_\Delta$ be the exceptional divisor over $\Delta$.
We have $E_\Delta\xrightarrow{\pi'} E $ is an isomorphism and 
\[
\widetilde{L}_t \defeq \pi'^*L_t = h'^*\cO(1,1) - 2\epsilon E_\Delta + t\cdot \pi'^*\cO_{\bP\cE}(1),\ \  0<\epsilon \ll 1
\] 
and 
\[\widetilde{L}_t|_{(\bP\cE|_{E_\Delta}) }= \cO_{\Delta}(1) +\epsilon \cO_{E_\Delta}(2) + t\cdot \pi'^*\cO_{\bP\cE}(1),\] 
where $E_\Delta \to \Delta$ is the projective bundle $\bP N_{\Delta/\bP^2\times \bP^2}$ and $\cO_{E_\Delta}(1)$ is the corresponding hyperplane class.

As discussed before, the stability analysis will be divided into two parts, $([Z],f)$
with $[Z] \in E\subset (\bP^2)^{[2]}$ or not.
For $[Z]\notin E$, it suffices to check the stability of $\left((p,q),f\right)$ with $p=[0,0,1]$ and $q=[1,0,0]$ with respect to $\widetilde{L}_t$.
In this case, the polynomial $f$ takes the form
\begin{equation}\label{deg7equa}
 \sum_{i=0}^4 y^{6-i}f_i(x,z) + xyz g_3(x,z) +x^2 z^2 h_2(x,z) =0
\end{equation}
where $f_i,g_i,h_i$ are  homogeneous polynomials of degree $i$ in $x,z$.

For $[Z]\in E$, after identifying $E$ and $E_\Delta$, we may assume $h'([Z]) = [0,0,1]\in \Delta \cong \bP^2 $ and the fiber of $h'$ is isomorphic to $\bP^1 \cong \bP(T_p\bP^2)$ with coordinates $[a,b]$, corresponding to tangent vectors $a\partial_x + b\partial_y$ at $[0,0,1]$.
The corresponding ideal $I_Z\subset k[x,y]$ is given by 
\[I_Z = (-ay+ bx, (x,y)^2).\]
Note that the fiber (hence the  entire exceptional divisor $E$) is contained in one single $\SL(3)$-orbit, hence we may assume $I_Z = (x,y^2)$ (corresponding to $[0,1]\in \bP^1$).
Then the $f\in H^0(I_Z^2\otimes\cO_{\bP^2}(6))$ takes the form \[
a_1z^4x^2 +z^3xf_2(x,y) + z^2 f_4(x,y) + z f_5(x,y) +f_6(x,y),\quad a_1\in k
\]
where $f_i$ are homogeneous polynomials of degree $i$.

\begin{lem}\label{lem:1-PS}
    For $n=2$ and reduced case, assume $\Supp[Z] = \{p,q\}$ as before, then \[\stab_{[Z]}\SL(3) \cong G \rtimes \langle g_0 \rangle\]
    where \[
    G= \left\{\begin{pmatrix}
        * & * & \\
          & * & \\
          & * & *
    \end{pmatrix} \in \SL(3)
    \right\} \text{\rm ~and~} 
    g_0 = \begin{pmatrix}
        0 & 0  & 1 \\
        0 & -1 & 0 \\
        1 & 0  & 0
    \end{pmatrix}
    \]
    Moreover, \[
    \hom{}{\bC^*}{ \stab_{[Z]}\SL(3)} = \left\{
\lambda(s) = \begin{pmatrix}
        s^m & a_1(s^m-s^i) & \\
          & s^i & \\
          & a_2(s^i-s^j) & s^j
    \end{pmatrix} \in G \,\middle|\, a_i \in k
    \right\}
    \]
\end{lem}

\begin{proof}
    The first statement follows from the fact that any $g\in \stab_{[Z]}\SL(3)$ satisfies either $g$ or $g_0g$ fixed $p,q$.
    This also implies that $g^2\in G$ for any $g\in \stab_{[Z]}\SL(3)$.
    For any 1-PS $\lambda$ in $\stab_{[Z]}\SL(3)$, $\lambda(s) = \lambda(\sqrt{s})^2 \in G$.
    The second statement then follows from direct computation.
\end{proof}

\begin{pro}\label{pro:GITanly7}
    When $d=7$ (i.e. $n=2$), the following holds:
    \begin{enumerate}
        \item For $0<t<\frac{1}{3}$, $\sstable{\bP \cE}=\emptyset$.
        % and for $t=\frac{1}{3}$, $V_2(t)$ is of dimension one.
        \item For $\frac{1}{2}-\epsilon < t < \frac{1}{2}$, $V_2(t)$ compactifies $P_7$ and has Picard rank two, where $0<\epsilon\ll 1 $.
    \end{enumerate}
\end{pro}

\begin{proof}
    First, we show that if $[Z]\in E$, then $([Z],f)$ is unstable for $0< t < \frac{1}{2}$.
    Take diagonal 1-PS $\lambda(s)=(s^{-a},1,s^a)$ with $a>0$, then 
    \[
    \begin{split}
            \mu_t([Z],f;\lambda(s)) &= (-1+2\epsilon)\cdot a + t\mu(f;\lambda)\\
            &\leq (2t+2\epsilon-1)\cdot a 
            <0 ,\quad \text{for~}0<t<\frac{1}{2} \text{~and~} 0<\epsilon \ll 1.
    \end{split}
    \]
    Now assume $0<t<\frac{1}{2}$, we may only consider $([Z],f)$ with $[Z]\notin E$.
    Then it suffices to assume $\Supp[Z] = \{p,q\}$ and test the stability of $((p,q),f)$ with respect to $\widetilde{L}_t$.
    Take diagonal 1-PS $\lambda'(s) = (s^a,s^{-2a},s^a), \,a>0$ and we have \[
    \begin{split}
    \mu_t((p,q),f;\lambda') &= -2a+t\mu(f;\lambda') \\
    &\leq 2a\cdot (3t-1)
    <0 ,\quad \text{for~}0<t<\frac{1}{3}.
    \end{split}
    \]
    This proves the statement in (1).
    
For the second statement, we first prove that the unstable locus of $\PP\cE$ with $[Z]\notin E$ is of codimension larger than two, and then apply similar techniques to compute the Picard rank of the GIT quotient, as in the degree $8$ case.
When $\frac{1}{3}<t<\frac{1}{2}$, by Lemma \ref{lem:1-PS}, one can directly compute the Hilbert-Mumford weight for $([Z],f)$ with $[Z]\notin E$.
This could be divided to the following cases:
\begin{enumerate}[label = (\roman*).]
    \item $i\geq \max\{m,j\}$: for equation $f$ in \eqref{deg7equa}, we have \[
    \mu_t([Z],f;\lambda) = i+t\cdot \mu(f;\lambda).
    \]
    We have $\mu(f;\lambda)\geq 0$ if any term in $y^6, y^5f_1(x,z),y^4f_2(x,z)$ exists.
    This means that in the locus destabilized by some 1-PS of this form, the corresponding $f$ has at least vanishing order three at $[0,1,0]$. 
    
    \item $i<\max\{m,j\}$: we may assume $\max\{m,j\} = m$. 
    Then there are two possibilities:
    \begin{enumerate}
        \item $m>i \geq j$:
        in this case,  if any term in $x^4z^2, x^4yz, x^4y^2$ exists, we have  
        \[
        \begin{split}
            \mu_t([Z],f;\lambda) &= i+t\cdot \mu(f;\lambda)  
            \geq i + t\cdot (4m + 2j) \\
            &\geq \frac{1}{3}(4m+2j+3i) = \frac{1}{3}(m-j)>0.
        \end{split}
        \]
        This means that in the locus destabilized by some 1-PS of this form, the corresponding $f$ has vanishing order at least three at $[1,0,0]$.

        \item $m\geq j>i$: in this case,   we have\[
            \mu_t([Z],f;\lambda) \geq 0,
        \] if terms $x^4z^2$ or $x^3z^3$ exist.
        If not, one sees further that the locus with no term $x^2z^4$ will be destabilized by 1-PS   $\lambda(s) = (s^2,s^{-3},s)$ as \[
        \mu_t([Z],f;\lambda)  \leq  -3+ 6\cdot t <0, \text{~for~} \frac{1}{3}<t<\frac{1}{2}.
        \]
        This locus parametrizes $f$ containing $L_{pq}$ as an irreducible component. 
        One may consider the locus with a nonzero $x^2z^4$ term but vanishing $x^4z^2$ and $x^3z^3$ terms.
        In this case, if the term $x^4yz$ is nonzero, then $([Z],f) $ will not be destabilized by 1-PS of this form for $\frac{1}{2}-\epsilon < t < \frac{1}{2}$ as
        \[
        \begin{split}
            \mu_t([Z],f;\lambda) &= i+t\cdot \mu(f;\lambda)  
            = i + t\cdot \max\{4m+i+j, 2m+4j\}  \\
            &\geq \frac{1}{2} \max\{3m+2i, 2m+4j+2i\}  -\epsilon \cdot A \\
            &= \frac{1}{2} \max\{m-2j, 2j\}-\epsilon \cdot A>0.
        \end{split}
        \]
        where $A= \max\{4m+i+j, 2m+4j\}  >0$.
        And if the term $x^4yz$ vanishes, take $\lambda(s) = (s,s^{-1},1)$, we have \[
        \mu_t([Z],f;\lambda) = -1+t\cdot \mu(f;\lambda)   < -1 + \frac{1}{2} \cdot \max\{2,2,2\} = 0.
        \]
    This means that in the locus destabilized by some 1-PS of this form, the corresponding curve on $Bl_{p,q} \bP^2$ is  either reducible containing the exceptional curve $E_p$ or tangent to the strict transform $\widetilde{L}_{pq}$ of the line passing through $p,q$ with $A_1$-singularities.
    And the locus the corresponding curve on $Bl_{p,q} \bP^2$ is smooth and tangent to the strict transform $\widetilde{L}_{pq}$ of the line passing through $p,q$ is semistable for $t\in (\frac{1}{2}-\epsilon, \frac{1}{2})$.        
    \end{enumerate}

    \end{enumerate}

Combined above, one sees that when $\frac{1}{2}-\epsilon < t < \frac{1}{2}$, the smooth curve $C\in |-2K_{Bl_{p,q}\bP^2}|$ or generically with $A_1$ singularities is GIT semistable.
Thus $V_2(t)$ gives a compactification of $P_7$ for $\frac{1}{2}-\epsilon < t < \frac{1}{2}$.
And the divisorial part of the unstable locus is just the locus  $([Z]=[p,q],f)$ with $ [Z] \in E$.
We have the following exact sequence \begin{equation}\label{eq:pic7}
    \mathrm{CH}_{24}(\bP\cE^{\rm us})\xrightarrow{\psi} \Pic  (\bP\cE)  \to \Pic (\bP\cE^{\rm ss} )\to 0.
\end{equation}
Hence, $\bP\cE^{\rm ss}$ has Picard rank two as $\bP\cE$ has Picard rank three.
As in the case of degree $8$, again by \cite[Proposition 5.1, \S 2.1]{KKV89}, we have \begin{equation} \label{eq:pic7'}
0 \to \Pic V_2(t) \rightarrow \Pic_{\SL(3)}\bP\cE^{\rm ss} \to  \prod_{x\in \cC}\chi(\SL(3)_x),
\end{equation}
where $\cC$ is a set of representatives of closed orbits and $\Pic_{\SL(3)}\bP\cE^{\rm ss} \hookrightarrow \Pic \bP\cE^{\rm ss}$.
Note that the two generators of $\pic \bP\cE^{\rm ss}$ are all $\SL(3)$-linearized, one sees that $\Pic_{\SL(3)}\bP\cE^{\rm ss}$ has rank two.
In particular, since the GIT quotient $V_2(t)$ remains the same for  $t \in (\frac{1}{2}-\epsilon,\frac{1}{2})$.
Choose $t_1, t_2$ in this range, and there is an integer $N$ such that $L_{t_1}^{\otimes N}$, $L_{t_2}^{\otimes N}$ descend to $V_2(t)$, which may serve as a basis of $\Pic_{\SL(3)}\bP\cE^{\rm ss}$ in $\bQ$-coefficients.
So $V_2(t)$ has Picard rank two for $\frac{1}{2}-\epsilon < t < \frac{1}{2}$.   

\end{proof}

Therefore, we have a reductive GIT space that compactifies $P_7$:
\begin{equation}  \label{gitd7}
    \overline{P}_7^{\rm GIT} \defeq V_2(t),\quad \frac{1}{2}-\epsilon <t<\frac{1}{2}.
\end{equation}
Combined with the stability analysis, we get the following results.
\begin{pro} \label{pro:gitd78}
For $d=7,8 $, let $\overline{P}^{\rm GIT}_d$ be the space defined in (\ref{gitd8}) or (\ref{gitd7}). 
Then $\overline{P}^{\rm GIT}_d$ is a partial compactification of $P_d$ and has Picard number $1$ for $d=8$ and $2$ for $d=7$.
\end{pro}

\subsection{GIT models for $d\leq 6$}
\subsubsection{GIT model for degree $5,6$ del Pezzo pairs}
Let $(X,C)$ be a  smooth del Pezzo pair of $d=5,6$, then $\aut(X)$ is reductive, and such a surface $X$ is rigid (cf.~\cite[Introduction and Theorem 2.1]{totaro_bott_2020}). 
Therefore, we set 
\begin{equation}\label{gitd56}
    \overline{P}_d ^{\rm GIT}:=|-2K_X| \q \aut(X).
\end{equation}

\begin{pro}\label{pro:git56}
If $d =5\ \hbox{or}\   6$,  $\overline{P}^{\rm GIT}_d$ is a partial compactification of $P_d$ and has Picard number $1$.
\end{pro}
\begin{proof}

For degree $5$ case, $\aut(X)\cong\mathfrak{S}_5 $ is finite by \cite[Theorem 8.5.8]{dolgachev_2012}. 
We have \[\abs{-2K_X} \q \aut(X)  \cong \abs{-2K_X} / \aut(X).\]
As the locus in $\abs{-2K_X}$ parametrizing curves with at least cuspidal singularities has codimension $\geq 2$, the statement follows.

For degree $6$ case, by \cite[Theorem 8.4.2]{dolgachev_2012} we have $\aut(X) \cong \big((\bG_m)^2 \rtimes \mathfrak{S}_3\big) \times \mathfrak{S}_2 $.
Without loss of generality, we may assume $X\cong Bl_{p_1,p_2,p_3} \bP^2$ where $p_i = [1,0,0], [0,1,0]$ and $[0,0,1]$ respectively.
Then the $(\bG_m)^2 \subset \SL(3)$ acts on $\bP^2$ diagonally while fixing $p_i\,(i=1,2,3)$, and thus acts on $X$.
Then $\abs{-2K_X} \cong \abs{\cO_{\bP^2}(6)\otimes \prod_{i=1}^3 \fm_{p_i}^2 }$, and a general $f$ takes the form \[
f=\sum_{\substack{0\leq i,j,k \leq 4 \\ i+j+k=6}} x^i y^j z^k
.\]
Direct stability analysis shows that $C \in \abs{-2K_X}$ is unstable if and only if $C = C_0 + E$ is reducible, where $E$ is an exceptional curve over some $p_i$. (This is also equivalent to $f$ having at least $D_4$ singularities at some $p_i$.)
Similar computation as in \eqref{eq:pic7'} shows that $\overline{P}_6^{\rm GIT}$ has Picard rank one.
\end{proof}

\subsubsection{GIT model for degree $3,4$ del Pezzo pairs}
If $d=3$, let   $\fX \subset \bP^3 \times \bP H^0(\bP^3,\cO_{\bP^3}(3))$ be the universal cubic surfaces with two natural projections 
\[ 
\pi_1\colon \fX \rightarrow \bP^3,
\quad \pi_2 \colon  \fX \rightarrow \bP^N \defeq \bP H^0(\bP^3,\cO(3))^\vee .
\]
Then $\cE:=\pi_{2\ast}(\pi_1^\ast \cO_{\PP^3}(2))$ is a locally free sheaf of rank $10$ on $\bP^N$ as all the higher direct images of $\pi_2$ vanish.
Indeed, $\cE\cong\cO^{\oplus 10}$. Thus we get a projective bundle $p:\bP \cE \cong \bP^9 \times \bP^N  \rightarrow \bP^N$ over $\bP^N$.  
Let $G =\aut(\PP^3)$ act on $\bP^N$ by the action induced from $\PP^3$ and  %By \cite[Proposition 2.7 ]{Benoist14}, the $\bQ$-line bundle 
$L_t:=\cO_{ \bP^9 \times \bP^N}(t,1),\ \ t\in \bQ$
be an ample $\bQ$-line bundle  for $ 1 \gg t >0$. Then we construct the GIT quotient space
\begin{equation}\label{gitd34}
   \overline{P}_3^{\rm GIT} \coloneqq \bP \cE \q_{L_t} G,\; 0<t<\epsilon\ll 1
\end{equation}
The same construction works for $d=4$, where one only needs to replace $\bP H^0(\bP^3,\cO_{\bP^3}(3))$ with the Grassmannian $\mathrm{Gr}(2,H^0(\bP^4,\cO_{\bP^4}(2)))$, and $\cE = \pi_{2*}(\pi_1^*(\cO_{\PP^4}(2))$.
Then 
\[\overline{P}^{\rm GIT}_4 \coloneqq \bP \cE \q_{L_t} G, \; 0<t<\epsilon\ll 1\] where $L_t = \cO_{\rm Gr}(1) + t\cdot \cO_{\bP \cE}(1)$. 

\begin{pro}\label{pro:git34}
If $d = 3,4$, let $\overline{P}^{\rm GIT}_d$ be the space defined in (\ref{gitd34}), then  $\overline{P}^{\rm GIT}_d$ is a   compactification of $P_d$ and the Picard number of $\overline{P}^{\rm GIT}_d$ is $2$. 
\end{pro}

\begin{proof}
For $d=3$, $\bP \cE \cong \cO_{\bP^N}^{\oplus  10}$ and thus $\bP \cE=|\cO_{\bP^3}(3)| \times |\cO_{\bP^3}(2) |\cong \bP^{19} \times \bP^9$. 
We claim that if $X=\{ c=0\}$ is a smooth del Pezzo surface of degree $4$ and $C=\{c=q=0\}\in |-2K_X|$ is a curve with only an $A_1$ singularity, or 
$\{c=0\}$ has $A_1$ singularities and $\{c=q=0\}$ is smooth, then $(c,q)$ is GIT stable. 
Note that if the rank of $q$ is less than or equal to $2$, then $C$ is reducible. 
Thus one can assume $\rank(q)\ge 3$. If $\rank(q)=4$, it is known that the smooth quadrics and nodal cubics in $\bP^3$ are GIT stable (cf.~\cite[Theorem 7.14]{mukai_invariants_2003}), then 
\begin{equation*}
\mu_t(c,q\,;\lambda)=\mu(c\,;\lambda)+t\cdot \mu(q\,;\lambda)>0
\end{equation*}
for any 1-PS in $G=\SL(4)$ and $t>0$. Now we only consider  $\rank(q)=3$. We argue by contradiction.  If the pair $(c,q)$ is GIT-unstable, there is a 1-PS $\lambda$ in the Weyl chamber of $\mathfrak{g}^\ast$ such that
\begin{equation}\label{constr}
    \mu(q\,;\lambda)<0, \ \mu(c\,;\lambda)> 0.
\end{equation}
Note that $c$ has at least $3$ monomials of degree $3$ and $q$ has at least $2$ monomials of degree $2$. By running a computer program\footnote{The code can be found in \url{https://changfeng1992.github.io/SiFei/GIT (3,2) pair in P^3.py}. }, there is no solution for the $\lambda$ satisfying the constraints in (\ref{constr}).

For $d=4$, recall that a degree 4 del Pezzo surface $X=(Q_1 \cap Q_2)$ is semistable if and only if it has at worst $A_1$ singularities in $Gr$, and stable if $X$ is smooth, under the natural action of $\SL(5,\bbC)$ by \cite[Theorem 6.1]{Mabuchi1993}.
A curve $C\in |-2K_X|$ which is a complete intersection $Q_1\cap Q_2\cap Q_3$ with at worst $A_1$ singularities is also stable in $Gr(3,H^0(\cO_{\PP^4}(2))$ by the results of Fedorchuk-Smyth \cite[Theorem 3.1]{Fedorchuk2013}.
So for any 1-PS $\lambda$, we have \[
\begin{split}
    \mu(X;\lambda) &=\mu(Q_1 \wedge Q_2 ;\lambda) \geq 0  \\
    \mu(X;\lambda) + \mu (C;\lambda) &= \mu (Q_1 \wedge Q_2 ;\lambda)+\mu ([Q_3];\lambda) \\
    &= \mu (Q_1 \wedge Q_2 ;\lambda)+\min \left\{\mu (Q_3';\lambda) \,\middle|\, Q_3' \equiv  Q_3 \bmod \langle Q_1,Q_2\rangle \right\}\\
    &\geq \mu (Q_1 \wedge Q_2 \wedge Q_3;\lambda) > 0
\end{split}
\] by definition.
We only need to show that $(X,C)$ is stable for $0<t \ll 1 $.
One may assume $\mu(C;\lambda) \leq 0$ otherwise the statement holds trivially.
In this case, we have
\[
\mu_t(X,C;\lambda) \defeq  \mu(X;\lambda) +  t\cdot \mu(C;\lambda) \geq  \mu(X;\lambda) + \mu(C;\lambda) >0.
\]

\end{proof}

\subsubsection{GIT model for degree $2$ del Pezzo pairs}
 %The geometry is a little different in this case.
 By \cite[Section 6.3.3]{dolgachev_2012}, a smooth del Pezzo surface of degree $2$ is a double cover of $\bP^2$ branched along a quartic curve. 
So let   $\cB \in \lvert\cO_{\bP^2 \times |\cO_{\bP^2}(4)|}(4,2)\rvert$  be a universal plane quartic curve and  $ \fX \xrightarrow{\phi}  \bP^2 \times |\cO_{\bP^2}(4)|$  the double cover branched along $\cB$ with  two natural projections
 \begin{equation*}
    \pi_1: \fX \rightarrow |\cO_{\bP^2}(4)|=\bP^{14},\ \  \ \ \pi_2: \fX \rightarrow \bP^{2}.
 \end{equation*}
Then  $\cE:={\pi_1}_\ast \left(\omega_{\fX/\bP^{14}}^{\otimes -2}\right)$ is a locally free sheaf of rank $7$ and we define the GIT space
 \begin{equation}\label{gitd2}
   \overline{P}^{\rm GIT}_2:= \bP \cE \q_{L_t} \PGL(3):=\proj(R(\bP \cE,L_t)^{\PGL(3)}), \ \ \  0<t\ll 1 
\end{equation}
where $L_t = \cO_{\rm \bP^{14}}(1) + t\cdot \cO_{\bP \cE}(1)$. 
\begin{pro}\label{pro:git2}
%If $d =2$, let $\overline{P}^{\rm GIT}_d$ be the space defined in (\ref{gitd2}), then  
$\overline{P}^{\rm GIT}_2$ is a   compactification of $P_2$ and the Picard number of $\overline{P}^{\rm GIT}_2$ is $2$. 
\end{pro}

\begin{proof}
    For $d=2$, we show that $(X, C)$ is stable if $X$ has at worst $A_1$ singularities and $C$ is smooth, or  $X$ is smooth and $C$ is general with $A_1$-singularities.
By direct computation, one gets \[
\cE = H^0(\bP^2,\cO_{\bP^2}(2)) \otimes \cO_{\bP^{14}}(-2) \oplus \cO_{\bP^{14}}(-3).
\]
One can identify the fiber of $p\mcolon \bP\cE \to \lvert \cO_{\bP^2}(4) \rvert$ over $B\in \lvert \cO_{\bP^2}(4) \rvert$ with $\bP V$, where 
$$V= H^0(X_B, -2K_{X_B}) \cong H^0(\bP^2, {\phi_B}_\ast (-2K_{X_B}) )\cong H^0(\bP^2, \cO_{\bP^2}(2) \oplus \cO_{\bP^2}).
$$
View $X_B$ as a surface in $\bP(1,1,1,2)$ defined by $w^2 = B(x,y,z)$.
Given $[(q,a)] \in \bP H^0(\bP^2, \cO_{\bP^2}(2) \oplus \cO_{\bP^2}) $, we get a $C\in \lvert -2K_{X_B} \rvert $ which is defined by $aw+q$ on $X_B$.
When $a=0$, $C$ is invariant under the  involution $w\mapsto -w$ and the double cover of the plane conic $Q \defeq [q]\in \bP H^0(\bP^2,\cO_{\bP^2}(2))$ branched over $Q\cap B$.
For $(B,q,a) \in \bP\cE$, the Hilbert-Mumford index is given by
\[
\begin{split}
    \mu_t(B,q,a\,;\lambda) &\defeq  \mu(B;\lambda) +t \cdot \mu(q,a\,;\lambda) {\phantom{:}=~} \mu(B;\lambda) +t \cdot \mu(q\,;\lambda).
\end{split}
\]
% Note that for $(X, C)$ attached to $(B,q,a)$, $X$ is smooth if and only if $B$ is smooth, and 
% \[ C \text{~is smooth} \iff \begin{cases}
%     Q \text{~and~} Q\cap B \text{~are smooth},  & a=0; \\
%     (q^2-a^2B)  \text{~is a smooth plane quartic} & a\neq 0.
% \end{cases}\]
% If $X=X_B$ is smooth and $C$ has $A_1$-singularities, there are two cases when $a=0$:
% \begin{itemize}
%     \item The plane conic $q$ is smooth and tangent to $B$ at exactly one point with intersection multiplicity $2$.
%     \item The plane conic $q$ has  an $A_1$-singularity where $B$ does not pass through.
% \end{itemize}
% If $a\neq 0$, this means $(q^2-a^2 B)$ has $A_1$-singularities away from $Q\cap B$.
% If $C$ is smooth and $X$ has $A_1$ singularities, then $B$ has $A_1$ singularities and $q$ is smooth when $a=0$.
Then in all these cases we consider, according to the general wall-crossing principle for VGIT, $(X,C)$ is stable for $0<t\ll 1$ as plane quartic curves with $A_1$-singularities are stable (cf.~\cite[Example 7.13]{mukai_invariants_2003}). 
Thus $\overline{P}_2^{\rm GIT}$ is a compactification of $P_2$.
One also sees that the unstable locus of $\bP\cE$ is of codimension at least two.
The statement of Picard number follows from the same argument as in the degree $7$ case \eqref{eq:pic7'}.
\end{proof}

 \subsubsection{GIT model for degree $1$ del Pezzo pairs}
By \cite[Theorem 8.3.2 (v)]{dolgachev_2022} a  Du Val del Pezzo surface $X$ of degree $1$ is a  hypersurface of degree $6$ in $\bP(1,1,2,3)$ and thus not passing through the singular points $[0,0,0,1]$ and $[0,0,1,0]$. (cf.~\cite[Theorem 3.1]{OSS16}).
Let $(x,y,z,w)$ be the weighted homogeneous coordinates.
After a coordinate change, we may assume the equation of $X$ is given by
 \begin{equation} \label{normalformd1}
       (w^2+z^3)=zf_4(x,y)+f_6(x,y).
 \end{equation}
as in \cite[Section 5.3.1]{OSS16}.
  Let $V\subset H^0(\bP(1,1,2,3),\cO_{\bP(1,1,2,3)}(6))$ be the subvector space spanned by the monomials in the right hand side of the equation (\ref{normalformd1}).
 There is a natural  action $\GL(2)$ induced from the action $\aut(\bP(1,1,2,3))$ on $ V$ ($\GL(2)$ acts on the first two weight one coordinates $x,y$).
%Let $\cS \subset \bP V \times \bP(1,1,2,3)$ the universal  bidegree $(1,6)$ hypersurface over $ V $ and the the projection $\pi_1:\ \cS \rightarrow  V$ onto  $\bP V$ and $\pi_2:\cS \rightarrow \bP(1,1,2,3)$ the projection on $\bP(1,1,2,3)$.  Denote  $\cE:=\pi_{1\ast} \pi_2^\ast \cO_{\bP(1,1,2,3)}(2)$ and then it is not hard to see $\cE$ is the free sheaf $\cO_{ V}\otimes H^0(\cO_{\PP(1,1,2,3)}(2))$ of rank $4$.  
 Let $W$ be the subvector space of $H^0(\cO_{\PP(1,1,2,3)}(2))$ spanned by $x^2,xy,y^2$ to parametrize curves of the form $(z+g_2(x,y) =0)\cap X$.
 (If the coefficient of $z$ vanishes, the curves are union of two curves in $\abs{-K_X}$ which are not smooth and irreducible).
 Now we define the GIT space for degree one del Pezzo pairs \begin{equation} \label{gitd1}
   \overline{P}_1 ^{\rm GIT}:= ((V\times W)/\bC^*)\q \PGL(2)\cong \bP(2,2,2,2,2,3,3,3,3,3,3,3,1,1,1)\q \PGL(2)
\end{equation}

\begin{pro}\label{pro:GIT1}
$\overline{P}^{\rm GIT}_1$ is a compactification of $P_1$ and the Picard number of $\overline{P}^{\rm GIT}_1$ is $1$.    
\end{pro}

\begin{proof}
    Let $0\neq (f,g)\in V\times W$ and we know that \[
    \mu(f,g;\lambda) = \max\{\mu(f;\lambda), \mu(g;\lambda)\}.
    \]
    By \cite[Lemma 5.17]{OSS16}, we know that $\{f=0\}$ is stable if and only if it contains at worst $A_k$ singularities.
    In this case, $(f,g)$ is stable for arbitrary $g$.
    So the unstable locus is of codimension at least 2.
    Then we know that $P_1\hookrightarrow \overline{P}_1^{\rm GIT}$ and the Picard number of $\overline{P}^{\rm GIT}_1$ is $1$.    
\end{proof}

\section{Class group of Baily-Borel compactification for \texorpdfstring{$P_d$}{P\_d}} \label{sect4}
In this section, we classify the discriminant locus for all $d$ and compute the divisor class groups of the Baily-Borel compactification of $P_d$.

\subsection{Arithmetic results for \texorpdfstring{$P_d$}{Pd}.} \label{subsec:arithmetic cptfy}
 We denote the Baily-Borel compactification of  $P_d$ by $P_d^\ast:=\cF_d^\ast\coloneqq (\Gamma_d\backslash \cD_{\rho(d)})^*$. Recall  $\Lambda =T_d=(\invLattice_d)^\perp\subset L$ where 
\begin{equation*}
\invLattice_d = \bZ l \oplus \bZ e_1 \oplus \cdots \bZ e_{9-d} \cong A_1(-1)\oplus A_1^{9-d}
\end{equation*}
is the invariant lattice of $\rho: L \rightarrow L$ with $l^2=2,\ e_i^2=-2, \ \langle e_i,e_j\rangle= \langle e_i,l \rangle =0 $ for $i \neq j$.
Then the discriminant form %\footnote{We warn that $\eta, \xi_i$ are only temporary notations here in the proof and one shall not confuse with the notations in Section \ref{arithmeticstra}.}
on $G(\invLattice_d)$  is  given by
 \begin{equation}
     q(\eta)=\frac{1}{4}\  \mathrm{mod}\  \bZ,\quad q(\xi_i)=-\frac{1}{4} \  \mathrm{mod}\ \bZ
 \end{equation}
where $\eta:=\frac{l}{2}$ and $\xi_i=\frac{e_i}{2}$. By Proposition \ref{prop:moduli meaning} and \eqref{eq:fine moduli}, \[
    P_d \subset \cF_d = \Gamma_d \backslash  \big( \cD_{T_d} - \Delta_d \big) \subset \Gamma_d \setminus \cD_d
\] 
where the discriminant locus \[
    \Delta_d \defeq \Gamma_d \backslash \bigcup_{v^2=-2, v\in T_d} v^\perp .
\]

To prepare computing $\Cl_\bQ(P_d^\ast)$, we are led to study the irreducible components of discriminant locus $\Delta_d$, which is reduced to the computation of $\Gamma_d$-orbits of various vectors in $T_d$.
\begin{lem}
    The irreducible component of $\Delta_d$ parametrizes 2-elementary K3 surfaces $(Y,\tau)$ with $\NS(Y)^\tau = \NS(Y)$ of rank $11-d$.
\end{lem}

\begin{thm}\label{Torelli}  Through canonical isomorphism \[
\epsilon \colon (G(T_d), q_{T_d}) \xlongrightarrow{\cong} (G(N_d), -q_{N_d}),
\] 
the $O(T_d)$-orbit of $(-2)$ vectors $v$ in $T_d$ are described as follows
    \begin{enumerate}
        \item  $d=8$ and $7$, there are two $O(T_d)$-orbit of $(-2)$ vectors $v$ where $v^*$ corresponds to  \[\{ 0\},\  \hbox{or}\  \{ \eta \}\ .\]

        \item  $d=6$,  then $O(G(T_6))\cong S_3 \times \bZ/2\bZ$ and  there are two $O(\Lambda_6)$-orbit of $(-2)$ vectors $v$ where $v^\ast $   corresponds to
        \[ \{ 0\},\ \  \hbox{or}\ \ \{ \xi_1+ \xi_2+\xi_3,\, 
   \eta \}. \]
        
        \item  $d=5$, $O(G(\Lambda_5)) \cong S_5$.
        There are three $O(\Lambda_5)$-orbits of $(-2)$-vectors $v$ with $\epsilon_5(v^*)$ belonging to\[
        \{ 0\},\{\eta,\,\xi_i+\xi_j+\xi_k\}_{1\leq i< j< k\leq 4}  \text{~or~} \{\eta+\sum^4_{i=1}\xi_i\}
        \]
        
        \item  $d=4$, there are two $O(\Lambda_4)$-orbits of $(-2)$-vectors $v$ with $\epsilon_4(v^*)$ belonging  to\[
        \{ 0\}, 
        \text{~or~} \{\eta,\,\xi_i+\xi_j+\xi_k,\,\eta + \xi_i+\xi_j+\xi_k+\xi_l\}_{1\leq i<j<k<l\leq 5}.
        \]
        \item  $d=3$, there are two $O(\Lambda_3)$-orbits of $(-2)$-vectors $v$ with $\epsilon_3(v^*)$ belonging  to
        \[
         \{ 0\}, 
        \text{~or~} \{\eta,\,\xi_i+\xi_j+\xi_k,\, \eta + \xi_i+\xi_j+\xi_k+\xi_l\}_{1\leq i<j<k<l\leq 6}.
        \]
        \item  $d=2$, there are two $O(\Lambda_2)$-orbits of $(-2)$-vectors $v$ with $\epsilon_2(v^*)$ belonging  to\[
        \{ 0\}, 
        \text{~or~} \{\eta,\, \sum_{i=1}^7 \xi_i,\,\xi_i+\xi_j+\xi_k,\, \eta + \xi_i+\xi_j+\xi_k+\xi_l\}_{1\leq i<j<k<l\leq 7}
        \]
        \item  $d=1$, there are three $O(\Lambda_1)$-orbits of $(-2)$-vectors $v$  with $\epsilon_1(v^*)$ belonging  to\[
        \begin{split}
          &\big\{\eta,\,\xi_i+\xi_j+\xi_k,\,\eta + \xi_i+\xi_j+\xi_k+\xi_l \big\}_{1\leq i< j< k\leq 8} \sim  \big \{ \ \sum_{I} \xi_I\  \big \}_{|I|=7, I\subset \{1,2,3,4,5,6,7,8\} },
         \\
        \text{~or~} &\{0\} \text{~or~} \   \{ \eta +\sum_{i=1}^8 \xi_i\}.
        \end{split}
        \]
    \end{enumerate}
\end{thm}
\begin{proof}

Let $d_v:=\divib_{T_d}(v)$ for any vector $v\in T_d$, then $d_v\in \{ 1,2\}$ for $(-2)$-vector.
Note that for any  $(-2)$-vectors with $d_v=1$, $v^\ast=0\in G(T_d)$ and Eichler's criterion in Lemma \ref{Eichler} shows these vectors are  in the same $\widetilde{O}(T_d)$-orbit and thus $O(T_d)$-orbit. 
So it remains to classify the $O(T_d)$-orbit of $(-2)$-vector with $d_v=2$.  
Note that by Nikulin's results \cite[Theorem 3.6.3]{Nikulin}, there is a surjection $O(T_d) \twoheadrightarrow O(G(T_d)) $.
Hence it's enough to classify the orbits of discriminant vector $u\in G(\Lambda)$ with $q(u)=-\frac{1}{4}$ under the action of $O(G(T_d))$.
By the canonical isomorphism $\epsilon_d$, it's equivalent to classify the orbits of $\frac{1}{4}$-norm vectors in $G(\invLattice_d)$.

For $d=7,8$, the only $\frac{1}{4}$-norm element in $G(\invLattice_7)$ is $\eta$, thus $\phi(\eta)=\eta$ for any $\phi\in O(G(\invLattice_7))$.

For $d=6$, note that the map 
\begin{equation*}
    \begin{split}
        \eta \mapsto \xi_1+\xi_2+\xi_3,\ \ \xi_1 \mapsto \eta+\xi_2+\xi_3; \\
       \xi_2 \mapsto \eta+\xi_1+\xi_3,\ \ \ \xi_3 \mapsto \eta+\xi_1+\xi_2
    \end{split} 
\end{equation*}
defines an element in $O(G(\invLattice_6))$.

For $d=5$, Note that the $\invLattice_5 \cong A_1(-1) \oplus A_1^{\oplus 4}$ and $\Lambda_5 \cong U^{\oplus 2} \oplus E_8 \oplus A_1^5$.
This comes from \cite[Lemma 2.3, 2.4, Remark 2.6]{Mkondo11}.

Based on the results of $d=5$ case, the proof for remaining cases are similar but rather lengthy.
We give the detailed proof for $d=3$  as a guiding example.
For $d=3$, one only need to show that $\eta$ and ``$\eta+4\xi$'' type vectors are in the same orbit.
It suffices to show there exits $\phi\in O(G(\invLattice_d))$ such that $\phi(\eta) = \eta+\sum_{i=1}^4 \xi_i $.
One can check directly that \[
\begin{split}
\eta \longmapsto \eta + \sum_{i=1}^4 \xi_i,  \quad  \xi_i \longmapsto \eta + \xi_i +\xi_5 ,\ 1\leq i \leq 4 ,  \quad 
\xi_5 \longmapsto \sum_{j=1}^{5} \xi_j,\quad\quad \;\,  
\xi_6 \longmapsto \xi_6
\end{split}
\]
gives the desired $\phi$.
For $d=2$, the following \[
\xi_7\longmapsto \eta+ \sum_{j=1}^{6}  \xi_j,  \quad \eta \longmapsto \sum_{i=1}^7 \xi_i,  \quad  \xi_i \longmapsto \eta + \xi_i +\xi_7 ,\ 1\leq i \leq 6 
\]
provides $\phi \in O(G(N_2))$ that connects $\eta $ and $\sum_{i=1}^7 \xi_i$.
For $d=1$, it is easy to see that $\eta$ and $\eta+\sum_{i=1}^8 \xi_i$ are not in the same orbit, since any element with the same norm as $\xi_i$ is not orthogonal to $\eta+\sum_{i=1}^8 \xi_i$.
\end{proof}
\begin{pro} Let $d=8$, then
the boundaries $\cF_d^\ast-\cF_d$ consist of four modular curves ( Type II  boundaries) and two points (Type III boundary).  
The root lattices of Type II  boundaries are given by \[ E_8\oplus E_7\oplus A_1,\  \ E_7\oplus D_9,\ \ D_{15}\oplus A_1,\ \ A_{16}. \]
\end{pro}
\begin{proof}
The method is standard. By \cite[\S 2.1]{Scattone1987}, the Type II (resp. Type III) boundaries of $\cF^\ast$ correspond to the classes of isotropic planes $J$ (rep. isotropic lines) in $\Lambda$.
 By \cite[Lemma 4.1.2]{Scattone1987},  Type III boundary consists of two points.
 Let $J\subset \Lambda$ be an isotropic plane, then $(J)^\perp_\Lambda/J$ is a negative definite  of rank 16 with the same genus as $ E_7\oplus A_1$.
 Using the embedding of  $ E_7\oplus A_1$ into a Niemeier lattices  \footnote{A unimodular definite lattice of rank $24$.}  which are classified into $24$ types (cf.~ \cite[Theorem 3.5.1]{Scattone1987}), we get the following possibilities
 \begin{equation*}
     E_8\oplus E_7\oplus A_1,\  E_7\oplus D_9,\ D_{15}\oplus A_1,\ A_{16}.
 \end{equation*}
 This classification is a prior up to $O(\Lambda)$. But as $O^+(\Lambda)/\widetilde{O}^+(\Lambda)=\{ \pm \id \}$, we get the classification for the Type II boundaries.
\end{proof}

\subsection{Nodal and weak del Pezzo divisor}\label{subsec:Nodes}
If the pair $(X,C)$ has nodal curve $C$, then   the double cover $S_0$ of $X$ branched along $C$ has a nodal point, then the minimal resolution $S$ of $S_0$ has Neron-Severi group 
 \begin{equation*}  \left(\begin{array}{c|cc} 
    & \invLattice_d  & E \\ \hline
    \invLattice_d & \invLattice_d  & 0 \\
    E & 0  &-2  \\
      \end{array}\right)  \end{equation*}
where we also use $\invLattice_d$ to indicate the Gram matrix and $E$ is the exceptional divisor of the resolution of the nodal point.  The period point of such pairs is just the irreducible divisor $H_v$ with $v\in T_d$ and $\divi(v)=1$. 
We call it the nodal divisor and denote it by  $H_{\rm n}$.
It lies in the complement of $F_{\rho(d)}$ in $\Gamma_\rho \backslash \cD_\rho$.

Assume $X$ is a normal surface with a $A_1$ singularities and $-K_X$ is ample with $d=(-K_X)^2$. 
It is called a nodal del Pezzo surface of degree $d$ and it's the ample model of weak del Pezzo surfaces (whose $-K_X$ is only big and nef). 
Another type of Heegner divisors we consider is associated to the K3 surface pairs $[(Y, \tau)] \in \cF_{\rho(d)}$ with $X\defeq Y/\tau$ that are weak del Pezzo surfaces of degree $d$.
Geometrically, there are two types of constructions:  
\begin{enumerate}
    \item take a blowup of $\bP^2$ along  $8-d$ points $(x_1,\cdots , x_{8-d})$ and then blow up a point $x_{9-d}$ at the exceptional divisor $E_{8-d}$ at $x_{8-d}$, we get a weak del Pezzo surface $\tilde{X}$ of degree $d$, i.e., $-K_{\tilde{X}}$ is big and nef. 
The strict transform $\bar{E}_{8-d}$ of $E_{8-d}$ is a $(-2)$-curve.
By contracting $\bar{E}_{8-d}$ one gets the ample model $X$ of $\tilde{X}$ which is a del Pezzo surface with an $A_1$ singularity;

    \item take a blowup of $\bP^2$ along  $9-d$ points $(x_1,\cdots , x_{9-d})$ with three points on a line $l$.
    This case only occurs when $d\leq 6$.
    When $d<6$, the weak del Pezzo surfaces of this kind coincide with the first case, and they also coincide with the weak del Pezzo surfaces coming from a blow up of $\bP^2$ along $9-d$ points with six of them lying on a conic. 
\end{enumerate}

\begin{defn}
    The irreducible Heegner divisor $(\Gamma_d \setminus \mathop{\cup} \limits_{\gamma \in \Gamma_d} \cD_{d,\gamma v})^{\rm red}$ is called the weak del Pezzo divisor and denote it by $H_w$.
\end{defn}

%\begin{rem} Note that in the case $d=8, a=1$,  $H_m$ is empty since a degree $8$ del Pezzo is obtained by blowing up $\bP^2$ at one point. \end{rem}

\begin{pro}\label{complement}
    The divisors in $\cF_{\rho(d)}-P_d$ parametrize K3 surfaces that are double covers of weak del Pezzo surfaces and \[
    D_{\rho(d)} / \Gamma_{\rho} = \cF_{\rho(d)} \cup H_n \cup \bigcup_{v^2=-2, div(v) =2} H_v
    \]
\end{pro}

\begin{proof}
    This follows from \cite[Section 4.A]{AE23}.
\end{proof}

\subsection{The class number of Baily-Borel compactifications}

By the GIT stability analysis in previous subsections, we can take  $U_d\subset \overline{P}_d^{\rm GIT}$ the locus with nodes at worst. 
(There are cases where only curves are allowed to be nodal as the ambient surfaces are fixed, as $d=8,7,6,5$. For $d=4,3,2,1$, both surfaces and curves are allowed to be nodal.) 
Let $\widetilde{U}_d:= F_{\rho(d)}\cup H_{\rm n }\subset \Gamma_{\rho}\backslash \cD_{\rho}$.
Since the K3 surface obtained from the del Pezzo surface pair given by a point in $U$ has ADE singularities at worst, there is a birational "period" map from $U_d$ to $\overline{F}_{\rho(d)}^{\rm slc}$ as they both contain $P_d$.
The image is contained in $\widetilde{U}_d$ and the period map is actually an open immersion for $U_d$.
By \cite[Theorem 4.2, Corollary 4.3]{AEH}, there is a 1-1 correspondence between the K3 surfaces with the involution fixed curve and the del Pezzo surface pairs.
Thus, two GIT orbits $[(X_1,C_1)]$ and $[(X_2,C_2)]$ with the same period will differ by an isomorphism preserving $-K_X$.
This implies the two orbits coincide in the GIT quotient and the period map is injective on closed points on $U_d$.
Then by Zariski’s main theorem and normality of $\Gamma\backslash \cD$ we conclude that the period map is an open immersion on $U_d$.

\begin{pro} \label{GIT-summary}
    %There is an big open subscheme $U_d$ of $\overline{P}_d^{\rm GIT}$ such that  
    The complement $\widetilde{U}_d-U_d$ is a union of the irreducible Heegner divisors given by  
    the following (up to a locus with codimension $\geq 2$):
    \begin{itemize}
        \item $d=8$: $\widetilde{U}_d-U_d =  H_h$;

        \item $d=7,5$: $\widetilde{U}_d-  U_d =  H_{w}$;
        
        \item $d=6$: $\widetilde{U}_d- U_d =  H_{w} \cup H_{w'}$;
        
        \item $d=4,3,2,1$: $\widetilde{U}_d-U_d = \emptyset$;
        
    \end{itemize}

\end{pro}

\begin{proof}
     By \cite[Section 4.A]{AE23},  $F_{\rho(d)}\setminus P_d$ is the union of weak del Pezzo divisors, which parametrizes nodal K3 surfaces that are double covers of nodal del Pezzo surfaces branched along smooth boundary curves as described in \S \ref{subsec:Nodes}.
     And $H_n$ parametrizes nodal K3 surfaces obtained from a double cover of smooth del Pezzo surfaces branched along irreducible nodal curves.
     Then for $1\leq d\leq 4$, $U_d$ is a big open subset of $\widetilde{U}_d$ by Proposition \ref{pro:git34}, \ref{pro:git2} and \ref{pro:GIT1}.
     
     For $d=5,6,7$, $U_d$ is a big open subset of $P_d \cup H_n$ by (the proof of) Proposition \ref{pro:GITanly7} and \ref{pro:git56}.
     For $d=8$, $U_d$ is a big open subset of  $P_d\setminus H_h$, by (the proof of) Proposition \ref{pro:GIT8anly}.
\end{proof} 

Now we provide the geometric strategy to compute the class group of $P_d^\ast$.  
Similar methods have been used in \cite{GLT} to compute the class group of moduli spaces of the quasi-polarized K3 surfaces with Mukai models. 
From the localization of Chow groups and the birational contraction map $P_d^\ast  \dashrightarrow \overline{P}^{\rm GIT}_d$, we have 
\begin{equation} \label{Pic}
    \dim _\bQ \Cl(P_d ^\ast)_\bQ=\dim \Cl(\widetilde{U})+ \# (P^*_d - \widetilde{U})=\dim \Cl(U)+\#(\widetilde{U}-U)+ \# (P^*_d - \widetilde{U})
\end{equation}
where  $\# (-)$ means the number of irreducible divisorial components of the corresponding locus.

\begin{thm}
    The class group of $P_d^\ast$ is generated by the Hodge line bundle $\lambda$ and irreducible Heegner divisors described in Theorem \ref{Torelli}.  
    Moreover, the class numbers are given by
\begin{table}[ht]
    \centering
    \begin{tabular}{ |c|*{11}{c|} }
    \hline
     $d$ & $1$  & $2$   & $3$  & $4$  & $5$  & $6$  & $7$ & $8$ & $8'$  & $9$     \\ \hline
     $\dim_\bQ \Cl(P_d ^\ast)_\bQ $  & $3$ &  $3$& $3$ & $3$  & $4$ &   $ 4$   & $4$    & $3$  & $2$  & $2$ \\  \hline 
     \end{tabular}
    \caption{Class numbers of Baily-Borel compactifications $P_d ^\ast$}
    \label{Picard of BBmoduli}
\end{table}
\end{thm}
\begin{proof}
This follows from the formula (\ref{Pic}) and Proposition \ref{GIT-summary}.
\end{proof}

\begin{rem}\label{rem:generator}
    Indeed we can match the generators in the above localization argument with the following Table \ref{tab:generators}.
 Here, $H_u$ is an irreducible Heegner divisor that parametrizes K3 surfaces  from  blowups of $\bP(1,1,4)$. $H_{u'}$ parameterizes the K3 pairs from another index $2$ of the pezzo for   $d=5$ and $d=1$.  The generators are not essential to the computation of class numbers but will be crucial for studying the Hassett-Keel-Looijenga program for moduli spaces of del Pezzo surface pairs of all remaining degrees, which will be the forthcoming work.  %of the second and third authors with Junyan Zhao.
\end{rem}

 \begin{table}[ht]
    \centering
 \[\bgroup
% `\def\arraystretch{1.2}
 \begin{array}{|c|*{4}{c|}|}
	 \hline & & & \\[-1em]
	d & P_d^*-\widetilde{U} &  \widetilde{U}-U &  \Pic(U) \\ \hline
	 d=8 & H_u    &  H_h  & H_{\rm n}   \\ \hline
 d= 7 &  H_u   & H_w   & H_{\rm n}, H_h \\  \hline 
	 d=6 &  H_u   &  H_w,H_{w'}  &   H_{\rm n}  \\ \hline
 d= 5 &  H_u,H_{u'}   &  H_w & H_{\rm n}  \\  \hline 
 d=4 &   H_u  &  \emptyset   &   H_{\rm n}, H_w  \\ \hline
 d= 3 &   H_u  & \emptyset   & H_{\rm n},H_w \\  \hline 
	 d=2&  H_u   &  \emptyset   &  H_{\rm n}, H_w \\ \hline
 d= 1 &  H_u ,H_{u'}  &  \emptyset   &    H_{\rm n} \parallelsum  H_w \\  \hline 
 \end{array}\egroup
 \]
    \caption{Generators of  $\Cl(P_d^*)$ }
    \label{tab:generators}
\end{table}

\noindent \textbf{\large A motivating question}.
Comparing different compactifications for a moduli space is an interesting topic. It dates backs to the work of Shah (\cite{Shah80}) and Looijenga (\cite{Loo86}) on the moduli space of  genus $2$ quasi-polarized K3 surface, Laza's (\cite{Laza10})  work   on moduli space  cubic $4$-folds. And later it is the work of  Laza-O'Grady  and  Ascher-DeVleming- Liu (\cite{ADL21}) on  moduli space of  genus $3$ quasi-polarized K3 surface, also the forthcoming work \cite{FLLST} on moduli of K3 surface in lower genus.   In our situation, a parallel question is 
\begin{question}
Can we resolve the birational period  map $\overline{P}^{\rm GIT}_d \dashrightarrow  \overline{P}^{\ast}_d$ explicitly, even with modular interpretation? 
\end{question}
Suggested by the known works mentioned above, we may approach the question by combining arithmetic methods and K-moduli theory. The K-moduli spaces of del Pezzo pairs are expected to be identified with a suitable Hassett-Keel-Looijenga model for $P_d^\ast$, which has a natural connection between the (partial) GIT compactification of $P_d$ and the Baily-Borel compactification of $P_d$. In other words, wall-crossing of K-moduli spaces should provide a modular  resolution of the birational period map $\overline{P}^{\rm GIT}_d \dashrightarrow  \overline{P}^{\ast}_d$. In the following sections we will give an arithmetic stratification for $P_8^\ast$, which will provide a typical example to understand the general cases. 

\section{Hassett-Keel-Looijenga program for \texorpdfstring{$P_8$}{P8}} \label{sect:HKL model}

The moduli space of hyperelliptic K3 surfaces has two components, which generically parametrize double covers of $\PP^1\times \PP^1$ and  $Bl_p\PP^2$ respectively.
For the former, Laza and O'Grady established the Hassett--Keel--Looijenga program in \cite{LaO18} and verified the program using the theory of variation of GIT.
Later, Ascher, DeVleming, and Liu realized the program using K-moduli of $(\bP^1 \times \bP^1,cC)$ in \cite{ADL21}.
Motivated by their work, we aim to establish the HKL program for the remaining component $P_8^*$ and give the arithmetic predictions of wall-crossing here. 

\subsection{HKL models for $P_8^\ast$ } \label{degree8hkl}
In \S  \ref{sec:GIT} and \S\ref{sect4} we established that the GIT model $\overline{P}_8^{\rm GIT}$ and the Baily-Borel compactification $P_8^*$ differ primarily by two divisors, the hyperelliptic locus $H_h$ and the unigonal locus $H_u$, see table \ref{tab:generators}. 
To resolve the birational map $\overline{P}_8^\mathrm{GIT} \dashrightarrow P_8^*$, we adopt the Hassett-Keel-Looijenga (HKL) program strategy and construct an intermediate moduli space $\mathcal{F}(s)$ by varying the coefficients of these specific boundary divisors. First, we give  the definition of  two divisors as Heegner divisors, which is a preliminary preparation to construct stratification on $H_h$ and $H_u$. Recall that the discriminant group  of the transcendental lattice  $T_8= U^2 \oplus E_8 \oplus E_7  \oplus A_1$ of the K3 surface from  a smooth del Pezzo pair $(Bl_p\bP^2,C)$ is given by $G(T_8)=\bZ \eta \oplus \bZ \xi,\ \ 2\eta=0,\ 2\xi=0, \eta\in G(E_7)$ and $\xi\in G(A_1)$.
\begin{defn} \label{Heegner}
  We call a 
 primitive vector $v\in T_8$ is
 \begin{enumerate}
     \item nodal if $v^2=-2$ and $\divib_{\Lambda}(v)=1$ ;
    \item unigonal if $v^2=-2$ and $\divib_{\Lambda}(v)=2$, $v^\ast = \xi$ ;
    \item   hyperelliptic if 
    $ v^2 =-6,\ \hbox{and}\   v^\ast=\xi
    $.
Then the complement $v^\perp$ of $v$ in $T_8$ is isomorphic to the complement  $(E_7\oplus A_2)^\perp$
   in $\II$. 
\end{enumerate}
\end{defn}
Denote by $H_h$  (resp. $H_u$) the primitive Heegner divisor associated with the hyperelliptic (resp. unigonal) vector $v$.
%Then a general point in $H_h$ recovers a K3 surface from the double cover of a general smooth pair $(Bl_p \bP^2,C)$ with $C$ in some special position, while a general point in $H_u$  recovers a K3 surface from the double cover of a general smooth pair $(Bl_p \bbF_4,C)$.
Consider the projective scheme $\cF(s)$ defined by the $\proj$ of the section rings 
\begin{equation*}
   R(s):= R(\cF, \lambda+s(H_h+25H_u)),\ \ s\in [0,1] \cap \bQ.
\end{equation*}
The locally symmetric variety $\Gamma \setminus \cD$ is strictly log canonical (see \cite[Theorem 3.2]{Alexeev96}) and   the finite generation of the section rings $R(s)$ is not known a priori except for $R(0)$, which gives back the  Baily-Borel compactification. Let $\overline{P}^K(c)$ be the component of K-moduli spaces generically parametrizing pairs $(Bl_p\PP^2cC)$ (see \cite{Xu21} for a survey of the construction of K-moduli space). Now we make the following conjecture which provides an explicit resolution of the birational period map $\overline{P}_8^\mathrm{GIT} \dashrightarrow P_8^*$.
%But we may assume this and give the prediction for wall-crossings of log canonical models.

%Note that the pair $(X,cC)$ be viewed as a log Fano surface for  $c\in (0,\frac{1}{2})\cap \bQ$. Recent breakthroughs in (log) Fano geometry ensure that the K-semistable log del Pezzo surfaces $(X,cC)$ form a good moduli space  $ \overline{P}_{d}^K(c)$ known as the K-moduli space(see \cite{Xu21} for a survey).  By the interpolation property of K-stability, %\cite{Tian90} and \cite{OSS16}  $ \overline{P}_{d}^K(c)$  is birational to $P_d$ whenever it is non-empty. According to the HKL principle explained in \S \ref{arith}, we make the following conjectural picture to relate  all the compactifications for $P_8$ appeared so far. 
\begin{conj}[HKL for $P_8^\ast$] \label{conj:HKL8}
Notation as above,
\begin{enumerate}
    \item The section rings $ R(s)$ are finitely generated for all $ s\in [0,1] \cap \bQ$. In particular, $\cF(s)$ is a projective variety of dimension $18$. 
    \item  $\cF(s)$ will interpolate between $P_8^\ast $  and $\overline{P}^{\rm GIT}_8$  when $s$ varies, i.e., 
    \[\cF(0)=P_8^\ast,\ \  \ \ \cF(1)= \overline{P}_8^{\rm GIT} .\] Moreover, there are finitely many rational numbers $0<w_1 <w_2<\cdots<w_n <1$ (called walls ) such that if $s, s'\in (w_i,w_{i+1}) \cap \bQ$ ,  there is a natural isomorphism $\cF(s)\cong \cF(s')$. Denote $ \cF(w_i,w_{i+1}):=\cF(s)$ for any $s\in (w_i,w_{i+1}) \cap \bQ$, then there is a wall-crossing behavior  (called HKL wall-crossing):
    \begin{equation*}
    \xymatrix@R-1em@C-1em{
   \cF(0,w_1)\ar[dr] &    &  \cF(w_1,w_2)   \ar[dl] \ar[dr] &   & \ar[dl]\cF(w_{n-1},w_n) \ar[dr]  & &  \ar[dl]  \cF(w_{n},1) \\
   & \cF(w_1) &      &  \cdots  &  & \cF(w_n) }
\end{equation*}
where  the birational maps are  flips or divisorial contractions at each wall $w_i$. 
    \item There is an isomorphism $\overline{P}^K(c) \cong \cF(s)$ under the transformation 
    \begin{equation*}
        s=s(c)=\frac{1-2c}{56c-4}.
    \end{equation*}
    In particular, the K-moduli wall crossing will  coincide with HKL wall-crossing .
\end{enumerate} 
\end{conj}

%\begin{rem}For del Pezzo pair $(\PP^1\times \PP^1,C)$, Laza-O'Grady \cite{LaO18}described  the decomposition of the birational map between the GIT compactification and the Baily-Borel compactification using VGIT  and the K-moduli realization was settled by \cite{ADL21}.In \cite{psw2}, we will verify the Conjecture \ref{conj:HKL8} via K-moduli method. The remaining degree cases will be dealt with in the future. \end{rem}

\subsection{Arithmetic predictions for wall-crossing of HKL }\label{arith} To find  the walls in the Conjecture \ref{conj:HKL8} from arithmetic side, we recall the method via arithmetic stratification in \cite{LO19} and further extended by \cite{FLLST}. 
For simplicity, we explain the principle in the case where the boundary consists of  only one prime divisor.  
Let $H_v=\Gamma  \setminus \sum_{g \in  \Gamma} \cD_{g\cdot v}$ be the Heegner divisor associated to $v$ on a locally symmetric variety $\cF$  of orthogonal type and $\Delta(s)=\lambda+a(s) H_v$, up to scaling, we may assume $s\in [0,1]\cap \bQ, a(0)=0$. Suppose the section ring \[ R(s):=\mathop{\bigoplus} \limits_{m} H^0(\cF, m l_s\cdot \Delta(s)) \]
is a finitely generated $\bC$-algebra where $l_s\in \bZ_{>0}$ is sufficiently divisible such that $l_s \Delta(s)$ is Cartier and thus $\cF(s):= \proj( R(s))$ defines a projective variety .
\begin{defn}
A self intersection $H_v^{(m)}$ for  $H_v$ with depth $m\in \bZ_{>0}$
in the sense of Looijenga is defined to be
\begin{equation*}
    \Gamma  \setminus   \bigcup_{(v_1,\cdots,v_m) } \cD_{v_1} \cap \cdots \cap \cD_{v_m}
\end{equation*}
where $(v_1,\cdots,v_m)$ runs over all collection of $m$ linearly independent vectors in the $ \Gamma $-orbit of $v$.
    
\end{defn}
 Let $m_v$ be the biggest depth of self-intersections for $H_v$ and there is a natural stratification (see \cite[\S 6,\S 7]{Lo03} for more details)
\begin{equation}\label{Looijenga's stratification}
    H_v^{(m_v)} \subset H_v^{(m_v-1)}\subset \cdots \subset H_v^{(2)} \subset H_v^{(1)}=H_v.
\end{equation}
A philosophy of HKL is that the birational  maps of $\cF(s)$ when varying $s$ from $0$ to $1$ should 
transform the stratification (\ref{Looijenga's stratification}) from the deepest stratum  to the less deep one successively. 
If we expect that the birational  maps of $\cF(s)$  can be obtained by running log MMP, consider the restriction of  $\Delta(s)$  on the expected  flipped locus $ H_v^{(n)}$  
\begin{equation*}
    \Delta(s)|_{H_v^{(n)}}=(1-a_n(s)) \lambda+Z_n
\end{equation*}
If $Z_n$ is birationally contractible on $H_v^{(n)}$ (in particular, $Z_n$ is an extrmal  in the effective cone of $H_v^{(n)}$), then $\Delta(s)|_{H_v^{(n)}}$ should lose positivity.
Hence $1-a_n(s)=0$ and  $s_n\in [0,1]$ such that $1-a_n(s_n)=0$ should be a wall.  If $Z_n$ fails to be  birationally contractible, then one need to modify the stratification (\ref{Looijenga's stratification}) from  depth $n$ stratum. Based on this philosophy, \cite[Section 10]{FLLST} outlines an algorithm to find all the potential walls.

In \S \ref{subsec: h arith strata} and \S \ref{subsec:u arithstra}, combined the geometric degeneration of K3 surfaces in Appendix \ref{Appa}, we give an explicit stratification on $\cF$ with respect to $H_h$ and $H_u$.
We apply the algorithm in \S \ref{arithpred} and Appendix \ref{Appb} to give prediction of walls and conjecture that it gives a complete description of HKL for moduli space of del Pezzo pairs of degree $8$.

%Note that the  points lying in $P_8\cup \cF_8^\circ$ correspond to K3 surfaces coming from the double cover of $Bl_p \bP^2$ or its nodal degeneration. Then the locus consisting of the double cover of the admissible degenerations of $Bl_p \bP^2$ will give rise to some arithmetic divisors in $\cF_8$ in the following way: let $X$ be some possible singular degeneration of $Bl_p \bP^2$ and $\widetilde{X}$ be its minimal resolution.Take the double cover of $\widetilde{X}$ branched along a curve $B\in \lvert -2K_{\widetilde{X}}\rvert $, and we get a possible singular K3 surface $S_0$. Then the crepant resolution $S$ of $S_0$ lies in $\cF_8$. We first characterize the general period point in $H_h$ and $H_u$ as K3 surfaces $S$ obtained from $X=Bl_p \bP^2$ or $Bl_p\bP(1,1,4)$ by computing their Neron-Severi groups. 

\subsection{Restriction of Heegner divisors} In this subssection, we make some preparations for explicit  construction of arithmetic stratification on  $\cF$. 
Let $\Lambda:=U^{\oplus 2} \oplus \Lambda_0$ be an even lattice of signature $(2,n)$ where $ \Lambda_0$ is negative definite. 
Let $v\in \Lambda$ be a primitive vector with $v^2<0$. Then 
\[ \bZ v \oplus v^\perp \subset \Lambda \]
is an overlattice and the quotient group $\Lambda/(\bZ v \oplus v^\perp) $ is a finite group. Denote by $N$ the  order of  $\Lambda/(\bZ v \oplus v^\perp) $, then $ N^2\cdot |G(\Lambda)|=|v^2|\cdot |\det(v^\perp)|$. Thus for any $w \in \Lambda$, there is a unique decomposition
\begin{equation} \label{vectorprojection}
    m\cdot w=a_w\cdot v+ x_w \in \bZ v \oplus v^\perp,\ a_w\in \bZ 
\end{equation}
where $m\in \bZ_{>0}$ is the order of the image of $w$ in the quotient group $\Lambda/(\bZ v \oplus v^\perp) $.
\begin{defn}\label{proj_vector}
    We denote by $\pi_{v^\perp}(w) $ the projection of vector $w\in \Lambda$ onto $v^\perp$ where \[\bZ_{>0} \cdot \pi_{v^\perp}(w) = (\bQ_{>0}\cdot x_w) \cap v^\perp.\] 
\end{defn}
The natural lattice embedding $v^\perp \hookrightarrow \Lambda$ induces a natural morphism between locally symmetric varieties 
\begin{equation*}
  \iota_v:  \sh(v^\perp) \longrightarrow \Gamma_v \setminus D_v \longrightarrow \Gamma \setminus D.
\end{equation*}
If the irreducible Heegner divisor $H_w$  on $ \Gamma \setminus D$ is different from $H_v$, then we have the pullback formula (see also \cite[Proposition 5.2.4]{LO19})
\begin{equation*}
  \iota_v^\ast H_w=\mathop{\sum} \limits_{w' } H_{\pi_{v^\perp}(w')}
\end{equation*}
where $w'\in \Lambda$ runs over the $\Gamma$-orbit such that the rank $2$ sublattice $\langle v,w'\rangle $ of $\Lambda$ is negatively definite.

\begin{pro}
    There is a commutative diagram 
    \begin{equation*}
        \begin{tikzcd}
          \sh(v^\perp)  \arrow[r]  \arrow[d] & \arrow[d]   \sh(\Lambda)  \\
            H_v \arrow[r,hook]   &   \Gamma  \setminus \cD
        \end{tikzcd}
    \end{equation*}
\end{pro}
\begin{proof}
  Let $\Gamma_v:=\{ \gamma \in \Gamma \mid \gamma(v)=v \}$ be the subgroup preserving the integral vector $v$. 
  Then $\Gamma_v \le  \stab_\Gamma(\cD_v)$ and thus there is a surjective morphism \[ \Gamma_v  \setminus \cD_v \rightarrow H_v=\stab_\Gamma(\cD_v) \setminus \cD_v  \] since $H_v= \Gamma \setminus \sum_{\gamma \in \Gamma } \cD_{\gamma  \cdot v}$. Then the proposition follows from the claim: there is an embedding of arithmetic groups  $\widetilde{O}(v^\perp) \le \Gamma_v \le O(v^\perp)$.

Proof of the claim:
Let $\Lambda_v:=v^\perp$ be the orthogonal complement of $v\in \Lambda$, then $\Lambda_\bC=\bC v \oplus (\Lambda_v \otimes \bC)$. 
We can write $\gamma(v)=av+b w$ where $w\in \Lambda_v \otimes \bC$  for any $\gamma\in \Gamma_v$. 
For any $x\in v^\perp$, $$0=\langle\gamma^{-1}x,v\rangle=\langle x,\gamma(v)\rangle=b\langle x,w\rangle,$$ which implies $b=0$. 
In particular, $\gamma(v)=a v$ and $a\in \bZ$ since  $\gamma(v)\in \Lambda$ and $v$ is primitive.
Observe that $\langle\gamma(v),\gamma(v)\rangle=|a|^2\langle v,v\rangle$ shows $|a|=1$. 
Since  $\gamma \in  \widetilde{O}^+(\Lambda)$, it  preserves $\cD^+_{\Lambda}$.  
On the contrary, any $\gamma \in \stabop(\Lambda)$ such that $\gamma(v) = \pm v$ preserves $\cD^+_{v}$. 
Note that $\gamma(v) =-v $ is possible if $v^* = -v^* \in A_\Lambda$ and $v^2 <0$.
In this way, we prove
$$\Gamma_v = \{ \gamma \in \widetilde{O}^+(\Lambda)\ | \  \gamma (v)= \pm v \} \leq O^+(\Lambda_v).$$  
 Let $\gamma \in \widetilde{O}^+(\Lambda_v)$. Since  \[\Lambda_v  \oplus \bZ v \subset  \Lambda \subset \Lambda^* \subset \Lambda_v^* \oplus \bZ v^* \, , \]
we may first extend $\gamma$ to an automorphism of $\Lambda_v^* \oplus \bZ v^*$ by setting $\gamma (v)= v$. By abuse of notation, we still denote it by $\gamma$.  Any element in $\Lambda \subset \Lambda_v^* \oplus \bZ v^*$ can be  written  as $xv+w$ where $w\in \Lambda_v^\ast, x\in \frac{1}{\divi_\Lambda(v)}\bZ$. 
Since $\gamma |_{\Lambda_v^\ast/\Lambda_v}=id$, there is an $u\in \Lambda_v$ such that $\gamma(w)=w+u$. Thus,
\begin{equation*}
    \gamma (xv+w)=x\gamma(v)+ \gamma(w)= (xv+w)+u \in \Lambda.
\end{equation*}
This shows $\gamma(\Lambda) \subset \Lambda$ and it preserves orientation. Since $\gamma$ acts trivially on $(\bZ v^*/\bZ v) \oplus  (\Lambda_v^\ast/\Lambda_v) $ and $A_\Lambda \le (\bZ v^*/\bZ v) \oplus  (\Lambda_v^\ast/\Lambda_v)$, it also acts trivially on $A_\Lambda$. It is not hard to see such extension is unique. 
Thus, we get an embedding $
   \widetilde{O}^+(\Lambda_v)  \hookrightarrow \Gamma_v \hookrightarrow O^+(\Lambda_v)
 $.
\end{proof}
\subsection{Arithmetic stratification on hyperelliptic divisor } \label{subsec: h arith strata}

Let $\cF= \Gamma \backslash \cD_\Lambda$ be the locally symmetric space associated with the lattice $\Lambda$ in \S \ref{2-K3}. 
We start by taking $\Lambda  = U^2 \oplus E_8 \oplus E_7  \oplus A_1$.  
By Nikulin's results  \cite[Theorem 1.14.4]{Nikulin}), there is a unique primitive embedding  \[\Lambda \hookrightarrow I\!I_{2,26}\cong U^{\oplus 2} \oplus E_8^{\oplus 3},\ \  \hbox{and}\ \  U_1\defeq (\Lambda)^{\perp}_{\II_{2,26}} \cong  E_7 \oplus A_1\subset I\!I_{2,26} .\] 
Inspired by the prediction in \cite{LaO18ii}, we may add a root vector to the $E_7\oplus A_1$ to get $E_7\oplus A_2$ and continue the procedure. In this way, we have

\begin{enumerate}
    \item the  $A_n$  lattice towers
\begin{equation*} 
\begin{split}
    E_7\oplus A_1 \subset  E_7\oplus A_2 \subset  E_7\oplus A_3 \subset  \cdots  \subset E_7\oplus A_5;
\end{split}
\end{equation*}
\item the  $D_n$  lattice towers
\begin{equation*} 
\begin{split}
    E_7\oplus D_3=E_7\oplus A_3 \subset  E_7\oplus D_4 \subset  \cdots \subset  E_7\oplus D_{10};
\end{split}
\end{equation*}

\item the  $E_n$  lattice towers
\begin{equation*} 
\begin{split}
    E_7\oplus E_6 \subset  E_7\oplus E_7 \subset  E_7\oplus E_8.
\end{split}
\end{equation*}
\end{enumerate}
Suggested by the geometric degeneration computation in \S \ref{hk3} and \S \ref{uk3}, we should also have modified towers of lattices. Now we describle them by Dynkin diagram: let $\diamond$ and $\circ$ be the new roots added to the ADE Dynkin diagram and   $\diamond$ has intersection number $1$ with the fixed polarization class (omitted in the Dynkin diagram). The modified towers are the complements of some hyperbolic lattices in the K3 lattice given by
\begin{enumerate}
    \item $A_n'$-type:  %lattice towers are the complements of some hyperbolic  lattices in K3 lattice, where the hyperbolic lattices are given by the following diagram
\[  \xymatrix{  & &  & &   \diamond \ar@{-}[d]  \\    \bullet_{{n}} \ar@{-}[r] &  \bullet_{{n-1}} \ar@{-}[r]& \cdots  \ar@{-}[r] & \bullet_{4} \ar@{-}[r] &  \ar@{-}[r]  \bullet_3  &   \ar@{-}[r]  \bullet_2   & \bullet_1   } \]
 $A_n''$-type:  %lattice towers are determined by the   hyperbolic lattices whose diagrams are given by the following 
\[  \xymatrix{  & &  & & \diamond \ar@{-}[d]  \\    \bullet_{{n}} \ar@{-}[r] &  \bullet_{{n-1}} \ar@{-}[r]& \cdots  \ar@{-}[r] & \bullet_{4} \ar@{-}[r] &  \ar@{-}[r]  \bullet_3  &   \ar@{-}[r]  \bullet_2   & \bullet_1   \ar@{-}[r] & \circ  } \]
$A_n'''$-type:  %lattice towers
\[\label{fig:A_n'''} \xymatrix{  & \circ \ar@{-}[d]  &  & & \diamond \ar@{-}[d]  \\    \bullet_{{n}} \ar@{-}[r] &  \bullet_{{n-1}} \ar@{-}[r]& \cdots  \ar@{-}[r] & \bullet_{4} \ar@{-}[r] &  \ar@{-}[r]  \bullet_3  &   \ar@{-}[r]  \bullet_2   & \bullet_1   \ar@{-}[r] & \circ  } \]

    \item $D_n'$-type: %lattice towers are given by the complements of some hyperbolic  lattices in K3 lattice  where the hyperbolic lattices are given by the following diagram
\[  \xymatrix{ &  & \bullet_{{n-1}} \ar@{-}[d]  &  & &  &   \\  \diamond \ar@{-}[r] &  \bullet_{{n}} \ar@{-}[r] &  \bullet_{{n-2}} \ar@{-}[r]&   \cdots  \ar@{-}[r] & \bullet_{4} \ar@{-}[r] &  \ar@{-}[r]  \bullet_3  &   \ar@{-}[r]  \bullet_2   & \bullet_1      } \]

\end{enumerate}
These towers of lattices motivate us to define the following towers of locally symmetric  varieties.
\begin{defn} \label{defA}
\phantom{a}
\begin{enumerate}
    \item Let $\Lambda_{A_n}$ be the orthogonal complement of $E_7\oplus A_{n}$ in the Borcherds lattice $\II_{2,26}$. 
    The discriminant group $G(\Lambda_{A_n})=\bZ/2\bZ \times \bZ /(n+1)\bZ$ is generated by $\xi,\zeta$ with 
\[
\xi^2\equiv -\frac{1}{2} \,(\bmod\,2),\quad
\zeta.\xi=0,\quad
\zeta^2\equiv \frac{n}{n+1} \,(\bmod\,2).
\]
Denote by $\sh(A_n):= \widetilde{O}(\Lambda_{A_n})\setminus \cD_{\Lambda_{A_n}}$ the associated  $19-n$ dimensional locally symmetric spaces. If $n=5,6$,  denote by $\Lambda_{A_n}'$ the orthogonal complement of the Neron-Severi lattice (\ref{NSA'}) in the K3  lattice.  We call it the $A_n'$ tower. If $n\ge 7$, similarly we define $\Lambda_{A_n}'', \Lambda_{A_n}'''$ and  $\sh(A_n'')$, $\sh(A_n''')$ via the Neron-Severi lattices (\ref{NSA''}) and (\ref{NSA'''}) .  We call them the $A_n''$-tower and the $A_n'''$-tower.

\item  Let $\Lambda_{D_n}$ be the orthogonal complement of $E_7\oplus D_{n}$ in the Borcherds lattice $\II_{2,26}$. 
The discriminant group 
\begin{equation*}
G(\Lambda_{D_n})=  \left\{ 
    \begin{array}{lc}
         \bZ/2\bZ \times \bZ /4\bZ , &\  n \ \hbox{odd},  \\
        \bZ/2\bZ \times (\bZ/2\bZ \times \bZ/2\bZ), &\  n \ \hbox{even}. \\
    \end{array}
\right.
\end{equation*}
 is generated by $\xi,\eta$ with \[
 \xi^2 = -\frac{1}{2} \,(\bmod\, 2), ~ \eta^2=\frac{n}{4} \,(\bmod\,2),~ \xi.\eta=0
\] for odd $n$
 and generated by $\xi,\eta_1, \eta_2$ with 
\[\xi^2 = -\frac{1}{2} \,\bmod\, 2, \ \ ~ \eta_1^2 = \eta_2^2=\frac{n}{4} \,(\bmod\,2),~ \ \  
(\eta_1 + \eta_2)^2 =1,~\xi.\eta_i=0
\]
when $n$ even. Denote by $\sh(D_n):= \widetilde{O}(\Lambda_{D_n})\setminus \cD_{\Lambda_{D_n}}$
the associated $19-n$ dimensional locally symmetric variety.

For $n \ge 6$, we denote $\Lambda_{D_n}'$ the orthogonal complement of the Neron-Severi lattice (\ref{NSD'}) in the K3  lattice.  We call it the $D_n'$-tower and $\sh(D_n')$ the associated  locally symmetric variety.

\item  Let $\Lambda_{E_n}$ be the orthogonal complement of $E_7\oplus E_{n}$ in the Borcherds lattice $\II_{2,26}$. 
The discriminant group $G(\Lambda_{E_n})=\bZ/2\bZ \times \bZ /(9-n)\bZ$ is generated by $\xi,\eta_E$ with 
\[
\xi^2 = -\frac{1}{2} \,(\bmod\, 2), \ \ ~ \xi.\eta_E=0,\ \ \eta_E^2=
\frac{n+2}{6} \,(\bmod\, 2), \  n=6,7.
\]
Denote by $\sh(E_n):= \widetilde{O}(\Lambda_{E_n})\setminus \cD_{\Lambda_{E_n}}$ the associated  dimensional locally symmetric variety.
\end{enumerate}

\end{defn}
By the primitive embeddings $\Lambda_{A_n} \subset \Lambda_{A_{n-1}},~ \Lambda_{D_n} \subset \Lambda_{D_{n-1}}$ and $\Lambda_{E_n} \subset \Lambda_{E_{n-1}}$, we have 
natural morphisms between the locally symmetric spaces
\begin{equation*}
   \sh(A_n) \rightarrow \sh(A_{n-1}),\ ~
   \sh(D_n) \rightarrow \sh(D_{n-1}),\ ~
   \sh(E_n) \rightarrow \sh(E_{n-1})
\end{equation*}
which are finite onto their images in general.  The morphism $\sh(A_n) \rightarrow \cF$ can be factored into a series of compositions of morphisms $\sh(A_n) \rightarrow \sh(A_{n-1})$ from the above lattice towers.  Similar results hold for other lattice towers.  This induces the diagrams (\ref{diah}) and (\ref{diau}).

\subsubsection{A-type tower and their modification}
\begin{defn} \label{Avector}
    A primitive  vector $v\in \Lambda_{A_n}$ is called 
    \begin{enumerate}
        \item $A_{n+1}$-type, if  
    $ v^2=-(n+1)(n+2),\ \ v^\ast=\zeta$
    so that $v^\perp  \cong \Lambda_{A_{n+1}}$. In this case, $[\Lambda_{A_n}: \bZ v \oplus v^\perp]=n+2$ and $\divi_{\Lambda_{A_n}}(v)=n+1$.

    \item  $D_{n+1}$-type, if 
    \begin{equation*}
        v^2= \begin{cases}
               -(n+1),  &  n\  \hbox{odd} \\
               -4(n+1),  & n \ \hbox{even}
        \end{cases}, \  \ \ 
        v^\ast=2\zeta
    \end{equation*}
    and $v^\perp  \cong \Lambda_{D_{n+1}}$. In this case, 
    \begin{equation*}
      [\Lambda_{A_n}: \bZ v \oplus v^\perp]=  \begin{cases}
             2,   & \  n\  \hbox{odd} \\
              4 , &\  n \ \hbox{even}
        \end{cases}
    \end{equation*}
    \item  $A_{n}'$-type, if for  $n=5,6$, $v^2=-150n+746,\ v^\ast= \xi +3\zeta$
    %\[ v^2= \begin{cases}  -4,  &  n=5\  \hbox{odd} \\ -154,  & n=6 \ \hbox{even}\end{cases},\ \ v^\ast= \xi +3\zeta \]
    so that $v^\perp  \cong \Lambda_{A_n'}$ and 
    $[\Lambda_{A_n}: \bZ v \oplus v^\perp]=9n-43,\ \ n=5,6$ and $\divi(v)=2$ for $n=5$. 
  \item  unigonal type, if $v^2=-2$ and $v^\ast =\xi$ so that sublattice $v^\perp$  has discriminant $(-1)^n\cdot (n+1)$ and $[\Lambda_{A_n}: \bZ v \oplus v^\perp]=2$, $\divi(v)=2$.
    
    \item  $E_6$-type, if for $n=5$ and $v^2=-2,\ v^\ast=3\zeta$ so that $v^\perp= \Lambda_{E_6}$.  In particular, $[\Lambda_{A_5}: \bZ v \oplus v^\perp]=2$ and $\divi(v)=2$.
    \end{enumerate}  
\end{defn}

Denote the  divisors
\begin{equation*}
    \begin{split}
       & H_{\rm {n}}(\Lambda_{A_n}):=(H_{v})^{\rm red} \ \hbox{for}\ v\ \hbox{nodal vector},\ \  H_{\rm u}(\Lambda_{A_n}):= (H_{\xi,-\frac{1}{4}})^{\rm red}\  ( \hbox{unigonal divisor}) \\
       & H_{A_{n+1}}(\Lambda_{A_n}):= (H_{\zeta,-\frac{n+2}{2(n+1)}})^{\rm red}, \  \ \ \ \ \  H_{D_{n+1}}(\Lambda_{A_n}):= (H_{2\zeta,-\frac{2}{n+1}})^{\rm red} \\
        & H_{A_{n}'}(\Lambda_{A_n}):= (H_{\xi+3\zeta,-\frac{17-n}{4(n+1)}})^{\rm red}, \ n=5,6,  \   \ \ \ \
          H_{E_6}:=(H_{3\zeta,-\frac{1}{4}})^{\rm red} \text{~for~} \sh(A_5)
    \end{split}
\end{equation*}
on $\sh(A_n)$  where $()^{\rm red}$ means taking reduced part. Note that $H_{D_{n+1}}(\Lambda_{A_n})=0$ for $n=1$ and \[H_h=H_{A_{2}}(\Lambda_{A_1}),\ \  H_{D_{3}}(\Lambda_{A_2})=H_{A_{3}}(\Lambda_{A_2}).\] 
By the triangular relation (\ref{trirel}), we have
\begin{equation*}
    H_{0,-1}(\Lambda_{A_n})= 2(H_{\rm {n}}(\Lambda_{A_n})+H_{\rm u}(\Lambda_{A_n}))
\end{equation*}
and $H_{\rm u}(\Lambda_{A_n})$ is irreducible.

Before we step into the deeper modifications and deal with the more complicated discriminant groups, we unify the notation here.
Let \[G = (\bZ / n_1\bZ) \times  (\bZ/n_2\bZ) \times \cdots  \times (\bZ/n_k \bZ),
\]
we denote the $(\Bar{i_1},\Bar{i_2},\cdots,\Bar{i_k})$ the corresponding element in $G$.

The discriminant group of $\Lambda_{A_n}'$ is given by
\begin{equation*}
    G(\Lambda_{A_n}')= 
      \bZ/ (17-n)\bZ ,\  \  n=5, 6.
\end{equation*}
\begin{defn} \label{1stmv}
    We call 
   \begin{enumerate}
       \item a primitive vector $v \in \Lambda_{A_5}'$ is $A_6'$-type if $v^2=-132$ and so that $v^\perp = \Lambda_{A_6}'$ and $ [\Lambda_{A_5}': \bZ v \oplus v^\perp]=11$.  
       \item  a primitive vector $v \in \Lambda_{A_6}'$ is $A_7'$-type if $v^2=-110$ and so that $v^\perp = \Lambda_{A_7}'$ and $ [\Lambda_{A_6}': \bZ v \oplus v^\perp]=10$.
       \item  a primitive vector $v \in \Lambda_{A_6}'$ is $A_6''$-type if $v^2=-88$ and so that $v^\perp = \Lambda_{A_6}''$ and $ [\Lambda_{A_6}': \bZ v \oplus v^\perp]=8$.
   \end{enumerate} 
 
\end{defn}

Define the Heegner divisors 
\begin{align*}
   H_n &:= (H_{\bar{0},-1})^{\rm red}, &
   H_{A_6'}(\Lambda_{A_5'})&:=H_{\bar{1},-\frac{11}{24}}, &
   H_{\rm u'}(\Lambda_{A_5'})&:=H_{\Bar{3},-\frac{1}{8}}, \\
   H_{A_6''}(\Lambda_{A_6'})&:=H_{\bar{1},\,-\frac{4}{11}},&
   H_{\rm u'(\Lambda_{A_6'})}&:=H_{\bar{3},-\frac{3}{11}},&
   H_{A_7'}(\Lambda_{A_6'})&:=H_{\bar{2},\,-\frac{5}{11}}
\end{align*}

on  $\sh(A_5')$ and $\sh(A_6')$.

\begin{rem} Indeed $H_{A_6'}$ is a movable divisor. This can be shown from the K-moduli side by \cite{psw2} but not purely from arithmetic computation.  We only give a heuristic reason from an arithmetic perspective. 
Indeed, the pullbacks $\iota'^{\ast}H_{\Bar{4},-\frac{1}{3}}, \iota'^{\ast}H_{\Bar{3},-\frac{1}{8}}$ can also be computed explicitly as Heegner divisors on $\sh(A_5')$. But as their discriminants are less than  that of the divisor $H_{A_6'}$, then they should be birationally contractable.  Thus  these  Heegner divisors will be viewed as irrelevant. For simplicity of computations, we will denote them by $Z_{\sh}$ on $\sh$.  
\end{rem}

The discriminant group of $\Lambda_{A_n}''$ is given by  $G(\Lambda_{A_n}'')=
       \bZ/ (26-3n) \bZ  $ for $ n=6, 7,8$.
%\begin{equation*} G(\Lambda_{A_n}'')=  \bZ/ (26-3n) \bZ,  \ \     n=6, 7,8 \end{equation*}
\begin{defn}
   We call 
   \begin{enumerate}
       \item a primitive vector $v \in \Lambda_{A_6}''$ is $A_7''$-type if $v^2=-40$ and so that $v^\perp = \Lambda_{A_7}''$ and $ [\Lambda_{A_6}'': \bZ v \oplus v^\perp]=5$.  
       \item  a primitive vector $v \in \Lambda_{A_7}''$ is $A_8''$-type if $v^2=-10$ and so that $v^\perp = \Lambda_{A_8}''$ and $ [\Lambda_{A_7}': \bZ v \oplus v^\perp]=2$.
   \end{enumerate} 
   Let $H_{A''_7} \defeq H_{\bar{1},-\frac{5}{16}}$ on $\sh(A''_6)$ and $H_{A''_8} \defeq H_{\bar{1},-\frac{1}{5}}$ on $\sh(A''_7)$.
\end{defn}

\subsubsection{D-type tower and their modification}
\begin{defn} \label{vectoronD}
    A primitive vector $v\in \Lambda_{D_n}$ is called 
    \begin{enumerate}

    \item unigonal type if $
    v^2=-2 , \divi(v)=2, v^* =\xi \in G(\Lambda_{D_{n}})$ so that $(v^\perp)^\perp_{\II_{2,26}} \cong E_8\oplus D_n $.
        \item $D_{n+1}$-type if  $v^2=-4$, $v^\ast=2\eta$ for $n$ odd and $v^\ast=\eta_1+\eta_2$ for $n$ even
   % \[ \ \ v^\ast=\tilde{\eta} \defeq  \begin{cases} 2\eta  &\  n\  \hbox{odd} \\ \eta_1+\eta_2  &\  n \ \hbox{even} \end{cases}  \]
    so that $v^\perp  \cong \Lambda_{D_{n+1}}$ and $[\Lambda_{D_n}:\bZ v \oplus v^\perp]=2$.
    \item  $D_{n}'$ type if  
    \[ v^2=\begin{cases}
        -4(10-n), &\  n =7,9\\
        -2 , & n=8
        %-4,  & n=9
    \end{cases},\ \ v^\ast=\eta' \defeq
        \begin{cases}
               \xi+\eta  &\  n\  \hbox{odd} \\
              \xi+\eta_1  &\  n \ \hbox{even}
        \end{cases}  \]
    so that $v^\perp  \cong \Lambda_{D_{n}}'$.
    \item  $E_{n+1}$ type if  $n \ge 5$ and 
    \[ v^2=\begin{cases}
               -4(8-n), &\  n =5, 7  \\
             -2 , &\  n =6, 
               %-4 , &\  n =7 
        \end{cases} ,  \ \ v^\ast=\eta_u \defeq
        \begin{cases}
               \eta  &\  n\  \hbox{odd} \\
              \eta_1 \text{~or~}\eta_2  &\  n \ \hbox{even}
        \end{cases}    \]
    so that $v^\perp  \cong \Lambda_{E_{n+1}}$ and 
    \[ [\Lambda_{D_n}:\bZ v \oplus v^\perp]= 
        \begin{cases}
               8-n,   &   \  n\  \hbox{odd} \\
                 1,   & \  n \ \hbox{even}
        \end{cases} . \]
   \end{enumerate}
\end{defn}
Denote the  divisors
    \begin{alignat*}{2}
        H_{\rm {n}}(\Lambda_{D_n})&:= (H_v)^{\rm red}\textbf{} \text{~for nodal vector~} v ,\quad \ \  
        H_{D_{n+1}}(\Lambda_{D_n}):= (H_{\tilde{\eta},-\frac{1}{2}})^{\rm red} \\
        H_{\rm u}(\Lambda_{D_n})&:= (H_{\xi,-\frac{1}{4}})^{\rm red}, \ \ 
        H_{E_{n+1}}(\Lambda_{D_n}):=  (H_{v_{E_{n+1}}^*})^{\rm red},\ \   
        H_{D_{n}'}(\Lambda_{D_n}):= 
            (H_{\eta',\frac{n-10}{8}})^{\rm red}, \  %n=7,8,9
    \end{alignat*}
on $\sh(\Lambda_{D_n})$ and $H_{E_{n+1}}(\Lambda_{D_n})=0$ for $n \le 3$.
Here note that \[
    (H_{\eta_u,- \frac{8-n}{8}})^{\rm red}= \begin{cases}
    H_{\eta}, & n \text{~odd} \\
    (H_{\eta_1})^{\rm red} + (H_{\eta_2})^{\rm red}, & n \text{~even}
    \end{cases}
    \;= \begin{cases}
       H_{E_{n+1}}(\Lambda_{D_n}), & n\neq 7 \\
       H_{E_8} + H_{D_8} , & n=7
   \end{cases}
\] and we make the convention that $H_{E_5} \defeq (H_{\eta_1})^{\rm red} + (H_{\eta_2})^{\rm red}$.

The  discriminant group of $\Lambda_{D_n}'$ is given by  
\begin{equation*}
    G(\Lambda_{D_n}')=
        \bZ/2  \bZ  \times \bZ / (10-n)\bZ ,\ 6 \le n\le 9. %\bZ/ 2\bZ \times 
\end{equation*}
\begin{defn}
    We call a primitive vector 
    \begin{enumerate}
        \item $v\in \Lambda_{D_7}'$ is a $\Lambda_{D_8}'$-type vector if $v^2=-6$ and $\divi(v) = 3$ so that $v^\perp \cong \Lambda_{D_8}'$ and   $[\Lambda_{D_7}':\bZ v \oplus v^\perp]=2$.
        \item $v\in \Lambda_{D_8}'$ is a $\Lambda_{D_9}'$-type vector if $v^2=-2$ and $\divi(v)=2$ so that $v^\perp=\Lambda_{D_9}'$ and $[ \Lambda_{D_8}': \bZ v \oplus v^\perp ]=2$. 
    \end{enumerate}
\end{defn}

Define the divisor $H_{D_{n+1}'}:= (H_{v_{D_{n+1}'}})^{\rm red}$ and $ 
   H_{\rm n}:=(H_v)^{\rm red}  \hbox{for}\   v \  \hbox{is nodal vector}  \ %H_{\rm u'}(\Lambda_{A_6'}):=H_{,-\frac{3}{11}}
   $
on  $\sh(D_n)'$ for $n=6,7,8$.

\subsubsection{E-type tower}

\begin{defn}
    A vector $v\in \Lambda_{E_n}$ is called 
     $E_{n+1}$-type if  
    $ v^2=4n-30,\  \ n=6,7,\ \   v^\ast=\eta_E $
    so that $v^\perp  \cong \Lambda_{E_{n+1}}$ and 
    $ [\Lambda_{E_n}:\bZ v \oplus v^\perp]=9-n,\ \ n=6,7$. 
\end{defn}
Denote the integral divisors
\begin{equation*}
    H_{\rm {n}}(\Lambda_{E_n}):= (H_{v})^{\rm red}, v \text{~is nodal ~}, \   H_{E_{n+1}}(\Lambda_{E_n}):= (H_{\eta,-\frac{8-n}{2(9-n)}})^{\rm red},\ H_{\rm u}(\Lambda_{E_n}):= H_{\xi,-\frac{1}{4}} 
\end{equation*}
on $\sh(E_n)$.  We have 
\[ 
    (H_{0,-1})^{\rm red}= H_{\rm {n}}(\Lambda_{E_m})+ H_{\rm u}(\Lambda_{E_m}) + \delta_{m,7}H_{E_{m+1}}(\Lambda_{E_m})  
\]
and $H_{\rm {n}}(\Lambda_{E_n}),H_{\rm u}(\Lambda_{E_n})$ are integral divisors on $\sh(E_n)$.

\subsection{Arithmetic stratification on unigonal divisor}\label{subsec:u arithstra}

Similar to the hyperelliptic case we will provide an arithmetic stratification on the unigonal divisor $H_u \subset \cF$.  The general point in $H_u\cap H_h$ comes from the period of the double covering of pairs $(Bl_p \bP(1,1,4), C)$ where $C$ is tangent to the exceptional divisor $E$ over $p$. Using the Hodge index theorem on can see that $H_u$ has no self-intersections in the sense of Looijenga, i.e., $H_u^{(2)}$ is empty.
    
\begin{defn}[Towers on unigonal divisor] 
 
Let
\[ U_1 \supset U_2' \supset  U_3' \supset U_4''\]  %\supset  U_5 \supset  U_6 
be the complement of Neron-Severi lattice  from  (\ref{nsU1}), (\ref{nsU2'}), (\ref{nsU3'}) and (\ref{nsU4''}) in K3 lattice.  

\end{defn}

By Nikulin's results \cite[Theorem 1.14.2]{Nikulin}, there are unique primitive embeddings $$U_1,U_2', U_3', U_4''\hookrightarrow \II $$ with complements  $E_8\oplus A_1 , E_8\oplus A_1 \oplus A_2 , E_8 \oplus A_2^2, E_8 \oplus A_2 \oplus D_4 $. 
This motivates us to consider the following towers \begin{equation*} E_8\oplus A_1 \subset E_8\oplus A_1^2 \subset E_8\oplus A_1 \oplus A_2 \subset E_8 \oplus A_2^2 \subset E_8 \oplus A_2 \oplus A_3 \subset E_8 \oplus A_2 \oplus D_4 \subset \II. \end{equation*}
Let $\Lambda_{R} = (E_8\oplus R)^\perp_{\II} $ for some root lattice $R$ and we denote the discriminant groups \begin{itemize}
    \item $G(\Lambda_{A_n\oplus A_m}) \cong \bZ \langle\zeta_{n}  \rangle \times \bZ \langle \zeta_{m} \rangle $ where $\ord{}{\zeta_i} = i+1$ with $\zeta_i^2 = -\frac{i+2}{i+1}$.

    \item $G(\Lambda_{A_n\oplus D_4}) \cong \bZ \langle\zeta_{n}  \rangle \times \bZ \langle \eta \rangle \times \bZ \langle \eta' \rangle $ where $\eta,\eta'$ are similar to those in the case of the hyperelliptic $D$-type tower.
\end{itemize}

%Then as in the study of towers for hyperelliptic divisor, we define 
\begin{defn} \label{uAn}
   We say a primitive vector 
   \begin{enumerate}
       \item $v \in  U_1= (E_8\oplus A_1)_{\II}^\perp$  is nodal if $v^2=-2$ and $\divi(v)=1$ so that $v^\perp \cong (E_8\oplus A_1^2)_{\II}^\perp$ and $[U_1: \bZ v \oplus v^\perp]=2$.  $v \in (E_8\oplus A_1)_{\II}^\perp$  is $A_2$-type if $v^2=-6$ and $\divi(v)=2$,  $v^\perp \cong (E_8\oplus A_2)_{\II}^\perp$. 
       
        \item  $v \in U_2=(E_8\oplus A_1^2 )_{\II}^\perp$  is $A_1 \oplus A_2$-type if $v^2= -6$ and $v^\ast =\zeta_1$ or $\zeta_1'$, so that $v^\perp \cong (E_8\oplus A_1 \oplus A_2 )_{\II}^\perp$ and $[U_2: \bZ v \oplus v^\perp]=3$. Denote 
        \[ H_{A_1 \oplus A_2}:=H_v ,\ \ v\  \hbox{is an}\  A_1 \oplus A_2 \hbox{-type vector} \] 
        
        \item  $v \in U_3=(E_8\oplus A_1 \oplus A_2 )_{\II}^\perp$  is  $A_2 \oplus A_2$-type if $v^2= -6$ and $v^\ast =\zeta_1$, so that $v^\perp \cong (E_8\oplus A_2^2)_{\II}^\perp$  and $[U_2: \bZ v \oplus v^\perp]=3$.
        Denote \[ H_{A_2 \oplus A_2}:=H_v ,\ \ v\  \hbox{is an}\  A_2 \oplus A_2 \hbox{-type vector} \] 
        
        \item  $v \in U_4=(E_8\oplus A_2 \oplus A_2 )_{\II}^\perp$  is $ A_2 \oplus A_3$-type if $v^2= -12$ and $v^\ast =\zeta_2$ or $\zeta_2'$, so that $v^\perp \cong (E_8\oplus A_2 \oplus A_3 )_{\II}^\perp$  and $[U_4: \bZ v \oplus v^\perp]=3$.
        Denote \[ H_{A_1 \oplus A_2}:=H_v ,\ \ v\  \hbox{is an}\  A_2 \oplus A_2 \hbox{-type vector} \] 
        
        \item  $v \in U_5=(E_8\oplus A_2 \oplus A_3 )_{\II}^\perp$ is $ A_2 \oplus D_4$-type if $v^2= -4$ and $v^\ast =2\zeta_3$ so that $v^\perp \cong (E_8\oplus A_2 \oplus D_4)_{\II}^\perp$   and $[U_5: \bZ v \oplus v^\perp]=2$.
        Denote \[ H_{A_2 \oplus A_4}:=H_v ,\ \ v\  \hbox{is an}\  A_2 \oplus D_4 \hbox{-type vector} \] 
   \end{enumerate}
\end{defn}

Let $\sh(U_n)$ be the locally symmetric variety associated to  $U_n$. Thus, we have \[  \sh(U_3)= \sh(U_2'),\ \ \sh(U_4)= \sh(U_3'),\ \ \sh(U_6)= \sh(U_4''). \]

\subsection{Arithmetic prediction for wall-crossings on  \texorpdfstring{$\cF(s)$}{d= 8 case}}\label{arithpred}

%\subsection{Arithmetic predictions for HKL}

According to the arithmetic principle discussed in \S \ref{arith},  we need to compute the  restriction of the $\bQ$-line bundle $\Delta(s)= \lambda+s(25H_u+H_h)$ on the arithmetic stratification. Now we describe the  stratification. Denote 
%\subsubsection{Arithmetic stratifications} 
\begin{equation*}
    \begin{split}
        &\NL_h(A_n):=\im(\sh(A_n) \rightarrow \cF), \ 
 \ \NL_h(A_n'):=\im(\sh(A_n') \rightarrow \cF),\\
        &\NL_h(D_n):=\im(\sh(D_n) \rightarrow \cF), \ 
 \ \NL_h(D_n)':= \im(\sh(D_n') \rightarrow \cF),\\
 & \NL_h(E_n):=\im(\sh(E_n) \rightarrow \cF)
    \end{split}
\end{equation*}
 the NL-loci on the moduli space $\cF$. Then we get an arithmetic stratification by the NL-loci 
\begin{enumerate}
    \item $A_n$-type stratification 
    \begin{equation*}
        \NL_h(A_4) \subset \NL_h(A_3)   \subset \NL_h(A_2) \subset  \NL_h(A_1)=\cF^\ast
    \end{equation*}
    and its modification 
    \begin{equation*}
        \NL_h(A_7'') \subset \NL_h(A_7') \subset \NL_h(A_6')   \subset \NL_h(A_5') \subset  \NL_h(A_5)
    \end{equation*}
    \item $D_n$-type stratification 
    \begin{equation*}
         \NL_h(D_7) \subset  \cdots   \subset \NL_h(D_4) \subset  \NL_h(D_3)=\NL_h(A_3) \subset \cF^\ast
    \end{equation*}
    and its modification
    \begin{equation*}
        \NL_h(D_9') \subset  \NL_h(D_8') \subset    \NL_h(D_7') \subset  \NL_h(D_7)
    \end{equation*}
    \item $E_n$-type stratification
    \begin{equation*}
        \NL_h(E_8) \subset \NL_h(E_7) \subset \NL_h(E_6) \subset \cF^\ast
    \end{equation*}
\end{enumerate}

On the unigonal divisor $H_u$, we similarly  define these Noether-Lefschetz loci 
\begin{equation*}
    \begin{split}
        \NL_u(U_n):= \im (\sh(U_n) \rightarrow \cF).
    \end{split}
\end{equation*}

%\subsubsection{Prediction of walls} Recall the $\bQ$-line bundle $ \Delta(s)=\lambda+s(25H_u+H_h)$ on $\sh(L)$.
\begin{thm} \label{cbf} Denote the $\varphi: \sh(L) \rightarrow  \cF$ the composition of morphism of towers for $L$ being one of the lattices in the  previous subsection. %\[ \{ \Lambda_{A_n}, \Lambda_{A_n}',\Lambda_{D_n},\Lambda_{D_n}', \Lambda_{E_n},\ U_n \}. \] 
 \begin{enumerate}
     \item  If \  $ L=\Lambda_{A_{n+1} }$ and  $ n\le 3$, then  $$ \varphi^\ast \Delta(s)=(1- n s) \lambda+s(H_{A_{n+2}}+n(1-\delta_{1,n}) H_{D_{n+2}} +25\cdot H_{\rm u}(\Lambda_{A_{n+1}})),$$
     where $\delta_{a,b}$ is the Kronecker delta function.
     \item    If $L=\Lambda_{A_5 }'$,  then  $ \varphi^\ast \Delta(s)=(1- 4 s) \lambda+ 2sH_{A_{6}'}+ \widetilde{Z}_A\ $
    where $\widetilde{Z}_A$ is the pullback of Heegner divisor $Z$ in Proposition \ref{pullA'}.
      \item    If $L=\Lambda_{A_6 }'' $, then  $ \varphi^\ast \Delta(s)=(1- 6s) \lambda+ 2sH_{A_7''} +\widetilde{Z}_A' \ $ %will be $$ (1- 6s) \lambda+ 2sH_{A_7''} +\widetilde{Z}_A' \  $$
      where $\widetilde{Z}_A'$ is the pullback of Heegner divisor $Z'$ in Proposition \ref{pullA''}. 
       \item   If $L=\Lambda_{A_7 }'''$, then $ \varphi^\ast \Delta(s)= (1- 8s) \lambda+\widetilde{Z}_A''$. %will be $$ (1- 8s) \lambda+\widetilde{Z}_A''\  $$
     \item   If  $ L=\Lambda_{D_{n+1} }$  and $3 \le n \le 5$, then
     $$\varphi^\ast \Delta(s)=
     (1- 2(n-1) s) \lambda+s\cdot \left[ 2 H_{D_{n+2}}+ 2^{n-2} H_{E_{n+2}} + 25 H_{\rm u}(\Lambda_{D_{n+1}}) \right].
     $$ 
     \item   If $L=\Lambda_{D_{n} }'$, then \[
     \varphi^\ast \Delta(s) = 
     \begin{cases}
         (1-10s)\lambda + 2sH_{D_8'} +Z_D, & n=7 \\
         (1-12s)\lambda + 4sH_{D_9'}  ,    & n=8 \\
         (1-16s)\lambda                  , & n=9 
     \end{cases}. \]
     
     \item  If \ $L=\Lambda_{E_n }$ and $ n=6,7,8$,  then  $ \varphi^\ast \Delta(s)=(1- \mu_{E_n} s) \lambda+25s\cdot H_{\rm u}(\Lambda_{E_{n}}) $ where $\mu_{E_6}=10, \mu_{E_7}=16, \mu_{E_8}=28$.
     \item  If $L=U_n$, then 
     \begin{equation*}
         \begin{split}
             \varphi^\ast \Delta(s)=&(1-25s) \lambda+s H_{A_2}  \ \hbox{on} \  \ \sh(U_1)  \\
             \varphi^\ast \Delta(s)=&(1-25s) \lambda+sH_{A_1 \oplus A_2}+s \widetilde{Z}_{U_2}   \  \ \hbox{on} \ \  \sh(U_2) \\
              \varphi^\ast \Delta(s)=&(1-27s) \lambda+sH_{A_2 \oplus A_2}+s \widetilde{Z}_{U_3}    \ \    \ \hbox{on} \ \sh(E_8\oplus A_1 \oplus A_2) \\
              \varphi^\ast \Delta(s)=&(1-28s) \lambda+sH_{A_2 \oplus A_3}+s \widetilde{Z}_{U_5}        \ \hbox{on} \ \sh(U_5) \\
              \varphi^\ast \Delta(s)=&(1-31s) \lambda+s \widetilde{Z}_{U_6}        \ \hbox{on} \ \sh(U_6) \\
         \end{split}
     \end{equation*}
     where $\widetilde{Z}_{U_i}$ are linear combinations pullback of $Z_{U_{i-1}}$  and other Heegner divisors in relation (\ref{relationU}).
 \end{enumerate}

\end{thm}
\begin{proof}
We prove the theorem by inductive computations on the towers. If $L=\Lambda_{A_n }$, then using the pullback formula (\ref{adjuctionAn}) repeatedly  we will have \[\varphi^\ast(\lambda+s(25H_u+H_h))=(1- n s) \lambda+sH_{A_{n+1}}+s(n-1)H_{D_{n+2}} +25s\cdot  H_{\rm u}(\Lambda_{A_{n+1}}) .\]
This proves the item (1).  Consider the pullback on $\sh(\Lambda_{A_5})$, and we have 
\[ (1-4s)\lambda+sH_{A_6}+3sH_{D_6}+2sH_{E_6}+25sH_{\rm u}\]
Since Heegner divisors $H_{A_6}, H_{D_6}, H_{E_6}, H_{\rm u}$ on $\sh(A_5)$ are different from $H_{A_5'}$, the pullback formula in Proposition \ref{pullA'} shows item (2).

Then, combining (2) in Proposition \ref{pullA'}, we get the pullback 
 \[  \varphi^\ast \Delta(s)=(1-4s)\lambda+sH_{A_6'}+\widetilde{Z}_A \]
 on $\sh(A_5')$   by 
 computation by (3) in Proposition \ref{pullA'}.  This finishes the proof of item (3). 

The remaining cases follow from the similar computations by applying pullback formulas (\ref{pullD}), (\ref{pullD'}) and (\ref{pullE2}).

\end{proof}

Based on the above result,  we make the following predictions for wall crossings of $\cF(s)$.
\begin{Pred} \label{pre1} The walls  $W_h$ for $\cF(s)$ from hyperelliptic divisor $H_h$ can be  divided into 5 types
\begin{enumerate}
    \item $A_n$ type walls are given by  $W_A=\{  \frac{1}{n} \ |\ n=1,2,3\ \} $;
      \item $A_n'$ type walls are 
      given by $ W_{A'}=\{\  \frac{1}{n} \ |\ n= 4,6,8\ \} $; 
    \item $D_n$ type walls are given by $ W_D=\{\  \frac{1}{2n}\ | \ n=2,3,4\   \} $;
      \item $D_n'$ type walls  are given by $W_{D'}=\{\  \frac{1}{2n}\ | \ n=5,6,8\   \}$;
    \item $E_n$ type walls are given by
     $W_E=\{\ \frac{1}{\mu_{E_n}} \} $
 where $\mu_{E_6}=10, \mu_{E_7}=16, \mu_{E_8}=28$. 
    %$ \mu_{E_n}$ is half of the number of roots in $E_7 \oplus E_n$ that are simply incident to a generator of $A_1$ under the primitive embedding \[ E_7 \oplus A_1 \hookrightarrow E_7 \oplus E_n .\]
\end{enumerate}
\end{Pred}

Assume the divisors $\widetilde{Z}$ in Theorem \ref{cbf}  are birationally contractible. By the arithmetic principle discussed in \S \ref{arith}, we get the predictions that walls are exactly the values so that the coefficient  of $\lambda$ in $\varphi^\ast \Delta(s)$ vanishes. In this way, we get the Prediction \ref{pre1}. Now let us give an explanation why the modifications appear. 
\begin{itemize}
    \item The modification for the $A$-type tower appears since $H_{A_5}$ should be a movable divisor on $\sh(A_4)$  and thus the corresponding NL locus $\NL(A_5)$ will be contracted earlier. That is,  $\NL(A_5)$ will be contracted together with $\NL(A_4)$ and  $\NL(A_5')$ will be contracted after $\NL(A_4)$. This forces the wall-crossing to jump from the $A_n$-tower to the $A_n'$-tower.  There are similar reasons to modification where $A_n''$-tower  and $A_n'''$-tower appear. 
    \item The  modification for the  $D$-type tower appears since  $H_{D_7}$ on $\sh(D_6)$ should be movable.  As $\rho(\sh(D_9'))=1$, then so is $\rho(\NL_h(D_9'))=1$, this shows that wall crossing at this stratum will terminate. 
\end{itemize}

\begin{Pred}\label{pre2}
The walls for $\cF(s)$ from unigonal divisor $H_u$  are  given by the NL locus $\NL_u$. The unigonal  walls  are \begin{equation*}   W_u=\left\{ \frac{1}{25}, \frac{1}{27}, \frac{1}{28}, \frac{1}{31} \right\}. \end{equation*}
\end{Pred}
%As before, we give an explanation for the modification.
By remark \ref{Humodify},  $H_{A_2}$ is big on $H_u$ and thus it will be contracted early and the modification  appears.  The walls are  predicted by (8) in Theorem \ref{cbf}.

\begin{rem}
From the arithmetic side, we cannot explain the termination of wall-crossing at stratum $\NL_u(U_6)$. 
We suspect that  the morphism $\sh(U_4'') \rightarrow \cF $ is not an isomorphism in its image $\NL_u(U_6)$ but just a finite degree morphism $\ge 2$.
\end{rem}
%\addtocontents{toc}{\SkipTocEntry}
\subsubsection*{Proof of Theorem \ref{mthm2}}
Under the assumption of birational contractibility and finite generation of the section ring $R(s)$, Theorem \ref{mthm2} follows from Theorem \ref{cbf}.

\section*{Appendix A: Geometric input - ADE degenerations}
\setupappendixcounters{A}
\label{Appa}

\subsection{Hyperelliptic K3 surfaces and  ADE degenerations}\label{hk3} 
Let $Bl_p \bP^2\xrightarrow{\pi} \bP^2$ be the blowup of  $\bP^2$ at $p=[0,0,1]\in \bP^2$ and $E$ be the exceptional curve.  Denote $C\in |-2K_{Bl_p \bP^2}|$   a general curve tangent to $E$ at only one point.  
  Let $\phi:  X \rightarrow  Bl_p \bP^2 $ be the double cover branched along the curve $C$. Then  $X$ is a smooth K3 surface. 
  The tangency condition of $E$ implies \[\phi^{-1}(E)=E_1+E_2\]
  where $E_i \cong \bP^1$ and $E_1.E_2=1$. Let $L\sim \phi^\ast \pi^\ast \cO_{\bP^2}(1)$, then $\NS(X)=\bZ [L]\oplus \bZ [E_1] \oplus\bZ [E_2]$ with the following intersection numbers
\[ L^2=2,\ E_1^2= E_2^2=-2,\ L.E_1=L.E_2=0,\ E_1.E_2=1.  \]

The K3 surface $X$ with antisymplectic involution from such del Pezzo pair is a general member in the hyperelliptic divisor $H_h$. Now we compute the Neron-Severi groups of the K3 surface of  the  possible geometric degeneration of pairs $(Bl_p \bP^2, C)$ which  is  the geometric input for arithmetic stratification in \S \ref{subsec: h arith strata}.  
First note that if  $\phi': X' \longrightarrow \bP^2$ is a double cover branched along a plane sextic curve $B=\pi(C)$ and $\mu': \widetilde{X}  \longrightarrow X' $ is a minimal resolution, there is a commutative diagram as follows
\begin{equation*}
\begin{tikzcd}
      \widetilde{X} \arrow[d,"\mu"]  \arrow[r,"="]  &  \widetilde{X}\arrow[d,"\mu'"] \\
      X \arrow[d,"\phi"]     & X' \arrow[d,"\phi'"] \\
  Bl_p \bP^2=Bl_p\bP^2 \arrow[r,"\pi"] & \bP^2     
\end{tikzcd}
\end{equation*}
Then we can compute the Neron-Severi group $\NS(\widetilde{X})$ via a covering of $\bP^2$. In this viewpoint a general element in $H_h$ is obtained from a plane sextic curve $B$ with only an $A_2$-singularity. 

\subsubsection*{Degeneration with $A_n$-singularities.}
If $B$ degenerates to an irreducible plane sextic curve $B_0=\pi(C_0)$ with  only  $A_n$ singularity at $p$ which is locally of the form 
  \[x^2+y^n=0,\]
  then $X'$ is a singular K3 with $A_n$-surface singularity under the double cover $\phi:  X' \rightarrow \bP^2  $.  Thus the  Neron-Severi group $\NS(\widetilde{X})$ is generated by  $H\sim \mu^\ast \cO_{\bP^2}(1)$ and the exceptional divisors $E_1,\cdots,E_n$. The Gram matrix  of $\NS(\widetilde{X})$ is given by 
  \begin{equation*}\label{NSA} 
  \left( \begin{array}{c|c c}    
    & H  &  \vec{E} \\ \hline  
 H  & 2 & 0  \\ 
   \vec{E} &0  &A_n  \end{array}\right)
    \end{equation*}
where $\vec{E}=(E_1,\cdots,E_n)$ and  $\NS(\widetilde{X})$ has rank $n+1$. Then we consider the following  tangency conditions.
\begin{enumerate}
    \item \textit{$A_n'$ modification: }
If $B$ degenerates to the case where there is a line $L=\{ x=0 \}$ tangent to $B$ and $B$ is irreducible, then the preimage  $\phi^{-1}(L)$ of $L$ will split into two lines $L_1, L_2 \ \subset X$. Indeed the tangency condition shows that the plane curve $B$ is defined by the following equation 
  \begin{equation*}
z^4xl(x,y)+z^3xf_2(x,y)+z^2xf_3(x,y)+zxf_4(x,y)+f_6(x,y)=0.
\end{equation*}
with the coefficient of $y^6$ nonzero.
Then locally around $p$  the germ $(B,p)$ is isomorphic to $ x^2+xy^n=0, \ \hbox{for}\ n \le 5$.
Thus locally around $p$, $X$ is defined by equation
\[ w^2= x^2+xy^n, \ \hbox{for}\ n \le 5,  \]
which is a $A_{2n-1}$-singularity. Then under the minimal resolution $\mu: \widetilde{X} \longrightarrow X $,
\[ \mu^\ast L_i=l_i+ \frac{E_1+\cdots + (n-1)E_{n-1}+ nE_n+(n-1)E_{n+1}+\cdots+E_{2n-1}}{2} \]
where $l_i$ is the proper transform of $L_i$.  Then $\NS(\widetilde{X})$ is generated by $E_1,\cdots, E_{2n-1}, l_1$ and $l_2$
and the  Gram matrix of $\NS(\widetilde{X})$ is 
\begin{equation}\label{NSA'} 
  \left( \begin{array}{c|c c c}    
    &  H   & \vec{E} & l_2 \\ \hline  
H   & 2   & 0 &  1 \\ 
  \vec{E} &  0 & A_n  &  w \\ 
  l_2  &1 &  w^t  & -2
   \end{array}\right)
    \end{equation}
    where $w^t=(0,0,\cdots,0,1,0,0)$. Moreover, $\NS(\widetilde{X})$ has rank $2+n$ and discriminant $(-1)^{n-1}(17-n)$.

\item \textit{$A_n''$ modification:}    
If $B$ degenerates to the
union of a line $L$ and a plane quintic curve $D$ such that $L$ and $D$ intersects at two points with multiplicity $4$ and $1$ respectively, then there is an additional class $\widetilde{L}$ in $\NS(\widetilde{X})$ which comes from the preimage of $L$ under the covering map. Then the new Gram matrix of $\NS(\widetilde{X})$ is 
   \begin{equation}\label{NSA''} 
  \left( \begin{array}{c|c c c}    
    & H  &  \vec{E}  & \widetilde{L} \\ \hline  
 H  & 2 & 0 &  1  \\ 
   \vec{E} & 0 &A_{n+1}  & l \\ 
    \widetilde{L} & 1 &l^t  & -2
   \end{array}\right)
    \end{equation}
    where $l^t=(0,0,\cdots,1,0,0,0)$. Thus,  the discriminant group 
 of the transcendental lattice of $\widetilde{X}$ is a finite group of order $17-n$ and $\NS(\widetilde{X})$ has rank $2+n$.
    
\item \textit{$A_n'''$ modification: } If in addition $B$ is union of a line $L$ and a plane quintic curve $D$ such that $D$ and $L$ are tangent at $p$, then the $\NS(\widetilde{X})$ is given by     \begin{equation}\label{NSA'''} 
  \left( \begin{array}{c|c c c}    
    & H  &  \vec{E}  & \widetilde{L} \\ \hline  
 H  & 2 & 0 &  1  \\ 
   \vec{E} & 0 &D_{n+2}  & l' \\ 
    \widetilde{L} & 1& (l')^t & -2
   \end{array}\right)
    \end{equation}
    where $(l')^t=(1,0,\cdots,0,1,0,0,0)$ with discriminant group $\bZ /12\bZ$.
    %(compare to Subsection \ref{subsec: arithmetic strata}).
\end{enumerate}

\subsubsection*{Degeneration with $D_n$-singularities.}   
If $C\in |-2K_{Bl_p \bP^2}|$ has only  $D_n$ ($E_n$ ) singularity at $p$, the Gram matrix of $\NS(\widetilde{X})$ is 
    \begin{equation*}\label{NSD} 
  \left( \begin{array}{c|cc}    
    & H  &  \vec{E} \\  \hline 
 H  & 2 & 0  \\ 
   \vec{E} &0  &D_n  \end{array}\right).
    \end{equation*}
If there is an additional line,  then  the Gram matrix of $\NS(\widetilde{X})$ is 
 \begin{equation}\label{NSD'} 
  \left( \begin{array}{c|c c c}    
    & H  &  \vec{E}  & \widetilde{L} \\ \hline  
 H  & 2 & 0 &  1  \\ 
   \vec{E} & 0 &D_n  & v \\ 
    \widetilde{L} & 1& v^t & -2
   \end{array}\right)
    \end{equation}   
  where $v^t=(1,0,\cdots,0,1,0)$. 
  For such smooth K3 surface $\widetilde{X}$, 
 $\NS(\widetilde{X})$ has rank $n+2$ and the discriminant group of its transcendental lattice is a finite group of order $4$. 

\subsubsection*{Degeneration with $E_n$-singularities.}
    The  plane sextic curve $B$ with $A_2$-singularity at $p$ degenerates  to the plane sextic with only $E_n$ singularity at $p$, then 
    \begin{equation*}\label{NSE} 
  \left( \begin{array}{c|cc}    
    & H  &  \vec{E} \\  \hline 
 H  & 2 & 0  \\ 
   \vec{E} &0  &E_n  \end{array}\right).
    \end{equation*}

\subsection{Unigonal K3 surfaces and  ADE type degenerations}\label{uk3}
We now carry out the analogous analysis for $H_u$.
A general point of $H_u$ corresponds to a K3 surface arising as the double cover of $Bl_p\bP(1,1,4)$, which is the surface model for the unigonal degeneration.
The surface $X:=Bl_p\bP(1,1,4)$ has only one singularity $v$. It is a quotient singularity of type $\frac{1}{4}(1,1)$. 
Let $\mu:\widetilde{X} \rightarrow X$ be the minimal resolution, then $\widetilde{X} \cong Bl_p \bbF_4$. 
Denote $E$ and $F$ the exceptional divisor over $p$ and  the quotient singularity.  Let $H_y\subset \widetilde{X}$ be the proper transform of the line $\{y=0\}$ passing $p$ and the quotient singularity $v$. By a double cover $\psi: S \rightarrow \widetilde{X}$ along a general smooth curve $B \in |-2K_{\widetilde{X}}|$, we get a smooth K3 surface $S$. 
Then the classes\[  E'=\psi^\ast E,\ 2F':=\psi^\ast F,\  H_y':=\psi^\ast H_y \]  form a basis of $\NS(S)$ and  by computation we have the following intersection numbers
\begin{equation}\label{nsU1}
  (E'){^2}= (F'){^2}= (H_y') {^2}=-2, \ \  E'. F'=0, \ \ E'.H_y'=2,\ \ F'.H_y'=1.
\end{equation}
\begin{comment}
  \begin{equation} \left(\begin{array}{c|ccc} 
      & E' & F' & H_y' \\\hline
    E' & -2 & 0 & 2 \\
    F' & 0 & -2 & 1 \\
    H_y' & 2 & 1 &-2  \\
    \end{array}\right) . \end{equation} % disc=2  
\end{comment}
The linear system $|-2K_{\widetilde{X}}|$ has a fixed component $F$ and a general element has the form $C+F$ where $C\in |\mu^*(-2K_X)|$ a smooth curve cannot pass through $v$ with $C\cdot F =0$.
   Hence unlike the hyperelliptic case, for $X = Bl_p\bP(1,1,4)$ we can not switch the order of taking the double cover and  the minimal resolution as the Picard number of the resulting K3 surfaces are different. 
Here we provide several degenerations for curve $B$ and compute the Neron-Severi groups of K3 surfaces associated with these degenerate curves.
\begin{enumerate}
    \item \textit{1st-level degeneration}: 
    $C$ passes through $E\cap H_y$ and $[0,1,0]$ and there are two irreducible singular (nodal) curves  $\Sigma_1$ and $\Sigma_2$  in $|-2K_{\widetilde{X}}-F|$ tangent to $C$ at the nodes lying on $E\cap H_y$ and   $[0,1,0]$ respectively.
      Then we have $\psi^*\Sigma_i= \Sigma'_i + \Sigma_i^{''}$.
    Since $\Sigma_i'$ and $\Sigma^{''}_i$ are both resolution of $\Sigma_i$, we have ${\Sigma'_{i}}^2={\Sigma^{''}_{i}}^2= 14 $.
    And $\Sigma'_1 \cdot \Sigma'_2 =15$ by Hodge index theorem. Then the Gram matrix of Neron-Severi group is given by
   % By the symmetry, we may assume that
   % \[\Sigma'_1 \cdot \Sigma'_2 = \Sigma'_1 \cdot \Sigma^{''}_2 = \frac{1}{2}\Sigma'_1\cdot \psi^*\Sigma_2 = \frac{1}{2}\Sigma_1 \cdot \Sigma_2 =16.
   % \]
\begin{center}
    \begin{equation} \label{nsU2'} \left(\begin{array}{c|ccccc} 
      &  E' & F'  & H_y'  & \Sigma_1' & \Sigma_2' \\\hline
    E' & -2 & 0 & 2  &2  &2 \\
    F' & 0 & -2 & 1  & 0 & 0 \\
   H_y'  & 2 &1 &-2 & 1 & 1 \\
      \Sigma_1 & 2 & 0 & 1  & 14 & 15 \\
     T_1   & 2 & 0 & 1  & 15 & 14\\
    \end{array}\right)  \end{equation} % disc=8
\end{center}

\item \textit{2nd-level degeneration}:
In addition to case (1), the exceptional divisor $E$ is tangent to the irreducible component $C$ of the branching locus and $E\cap F =\emptyset$.
Then $\psi^*E= E'_1 + E'_2$ are two rational curves and $-2=2E^2 = (\psi^*E)^2 = (E'_1 + E'_2)^2$. 
Hence $E_1'\cdot E_2'=1$ and $E'_i \cdot D = \frac{1}{2}\psi^*E \cdot D$ for $D = H_y' $ or $\Sigma'_j$ . Then the Gram matrix of Neron-Severi group is given by
\begin{center}
  \begin{equation}\label{nsU3'} \left(\begin{array}{c|cccccc} 
      &  E'_1 & E'_2 & F'  & H_y'  & \Sigma_1' & \Sigma_2' \\\hline
    E'_1  &  -2 & 1 & 0 & 1 & 1 & 1 \\
 E'_2  & 1 & -2 & 0 & 1 & 1 & 1 \\
F'  & 0 & 0 & -2 & 1 & 0 & 0 \\
H_y'  & 1 & 1 & 1 & -2 & 1 & 1 \\
\Sigma_1' & 1 & 1 & 0 & 1 & 14 & 15 \\
\Sigma_2' & 1 & 1 & 0 & 1 & 15 & 14 \\
    \end{array}\right)  \end{equation} % disc =-12
\end{center}

\item  \textit{3rd-level degeneration}: the branched curve is $B=C\sqcup F$ where $C=C_0\cup E$ has two irreducible components and $C_0$ intersects with $ E$ at three distinct points. 
Then the associated K3 surface $S$ has three $A_1$ singularities. 
Let $\widetilde{S}$ be its minimal resolution, then there are three exceptional curves $L_1, L_2$ and $L_3$.  
Clearly $L_i\cdot E'=1 ,\,L_i^2 =-2$ and $L_i$ does not intersect with $F',H_y'$ and $\Sigma'_j$. Then the Gram matrix of Neron-Severi group is given by    
\begin{center}
      \begin{equation} \label{nsU4''} \left(\begin{array}{c|cccccccc} 
      &  E' & F'  & H_y'  & \Sigma_1' & \Sigma_2' & L_1 & L_2 & L_3
      \\\hline
    E' & -2 & 0 & 1 & 1 & 1 & 1 & 1 & 1 \\
    F' & 0 & -2 & 1 & 0 & 0 & 0 & 0 & 0 \\
   H_y'  & 1 & 1 & -2 & 1 & 1 & 0 & 0 & 0 \\
   \Sigma_1' & 1 & 0 & 1 & 14 & 15 & 0 & 0 & 0 \\
   \Sigma_2' & 1 & 0 & 1 & 15 & 14 & 0 & 0 & 0 \\
  L_1  & 1 & 0 & 0 & 0 & 0 & -2 & 0 & 0 \\
    L_2 & 1 & 0 & 0 & 0 & 0 & 0 & -2 & 0 \\
   L_3 & 1 & 0 & 0 & 0 & 0 & 0 & 0 & -2
    \end{array}\right)  \end{equation} %disc = -16
\end{center}
\end{enumerate}

\vspace{0.3cm}

\section*{Appendix B: Arithmetic computation}
\setupappendixcounters{B}
\label{Appb}

\subsection{Pullback on hyperelliptic divisors towers}
We denote $\delta_{i,j} = 1$ for $i=j$ and $0$ for $i\neq j$ as the Kronecker symbol.
\subsubsection{Pullback of $A_n$ towers}
\begin{pro}The Picard group  $\pic(\sh(A_n))_\bQ$ is generated by $H_{\rm {n}}(\Lambda_{A_n})$, $H_{\rm u}(\Lambda_{A_n})$, $H_{A_{n+1}}(\Lambda_{A_n})$ and $H_{D_{n+1}}(\Lambda_{A_n})$. 
The Picard number of $\sh(A_n)$ is given in the following  Table:
\begin{table}[ht]
    \centering
    \begin{tabular}{ |c|*{9}{c|} }
    \hline
     $n$ &   $2$ & $3$ & $4$  & $5$  & $6$  & $7$ & $8$ & $9$   & $10$ \\ \hline
    $\rho(\sh(A_n)^\ast)$  &$3$ & $4$ & $4$  &  $5$ & $5$ & $6$  & $5$  &  $5$ & $4$ \\  \hline 
     \end{tabular}
    %\caption{Picard numbers of $A_n$-towers}\label{Picard of Atower}
\end{table}

Moreover if $n \le 4$,  there is a relation:
\begin{equation}\label{hodgerel1}
\begin{split}
   (75+\frac{n(n+1)}{2}) \lambda_{\sh(A_n)}=&H_{\rm {n}}+57 H_{\rm u}+(n+1) H_{A_{n+1}}+(1-\delta_{2,n})\frac{n(n+1)}{2}H_{D_{n+1}}. %(\Lambda_{A_n})
\end{split}
\end{equation}

If $n=5,6$, in addition to the relation:
 \begin{equation}
 \begin{split}
     90\lambda_{\sh(A_5)}=H_{\rm {n}}(\Lambda_{A_5})+57H_{\rm u}+6H_{A_6}+15H_{D_6}+20H_{E_6},\\
     96\lambda_{\sh(A_6)}=H_{\rm {n}}(\Lambda_{A_6})+57H_{\rm u}+7H_{A_7}+21H_{D_7}+35H_{E_7},
 \end{split}
\end{equation}
there is a second relation: 
\begin{equation*}
\begin{split}
    H_{A_6}=2\lambda_{\sh(A_5)}-4H_{\rm u}+2H_{D_6}+2H_{E_6}+2H_{A_5'}.
    %\\H_{A_7}=3\lambda_{\sh(A_6)}-2H_{\rm u}+H_{D_6}+2H_{A_6'}
\end{split}
%2\lambda_{\sh(A_5)}=H_{A_6}-2H_{E_6}+2H_{A_5'}-2H_{D_6}-H_{\xi^{(5)}_1+ 3\xi^{(5)}_2}
\end{equation*}
In particular, $H_{A_6}$ is big on $\sh(A_5)$ since the divisors $H_{\rm u},H_{D_6}, H_{E_6},H_{A_5'}$ are birationally contractible.
\end{pro}
\begin{proof}
By putting the lattice  $L=\Lambda_{A_n}$ into the software $\mathbf{Sage}$ package \href{https://github.com/btw-47/weilrep}{weilrep},%\footnote{This is a Sage package developed by Brandon Williams, for computing with vector-valued modular forms, Jacobi forms and theta lifts. See \cite{Williamsa}.}
 we get the Picard numbers of $\sh(\Lambda_{A_n})$ directly by Theorem \ref{modularform1}.  
For the relation (\ref{hodgerel1}), we explain how to deduce it through an example. 
Suppose  $L=\Lambda_{A_1}$ and its discriminant group is isomorphic to 
 $\bZ \eta \oplus \bZ \xi$. We may denote   \[ \vec{e}_1=e_0, \  \vec{e}_2=e_\zeta,\ \vec{e}_3=e_\xi,\ \vec{e}_4=e_{\zeta+\xi}  \]
as the basis for the $4$-dimensional vector space $\bC[A_L]$ associated with the discriminant group $A_L$.%. Here $e_v$ means the basis associated with an element in the discriminant group. 
The output of $\mathbf{Sage}$ \footnote{The data of modular form computation is \href{ https://changfeng1992.github.io/SiFei/Sage notebook for Hodge relations.pdf}{here}.}  is a basis of vector-valued modular forms with $q$ expansion as follows:
\begin{equation*}\label{qexan}
    \begin{split}
        f_0=& (1-152q+o(q))\cdot \vec{e}_1+ o(q) \cdot \vec{e}_2+o(q) \cdot \vec{e}_3+(-112q+o(q)) \cdot \vec{e}_4, \\
        f_1=&(-56q+o(q))\cdot \vec{e}_1+  o(q) \cdot  \vec{e}_2+(q^{\frac{1}{4}}+o(q) )\cdot \vec{e}_3+(56q+o(q)) \cdot \vec{e}_4 ,\\
        f_2=&(-2q+o(q))\cdot \vec{e}_1+(q^{\frac{3}{4}}+o(q))   \cdot  \vec{e}_2+ o(q) \cdot \vec{e}_3+(2q+o(q)) \cdot \vec{e}_4.
    \end{split}
\end{equation*}
By theorem  \ref{modularform1} the above vector-valued modular form gives the relation 
\[ 152 \lambda_{\sh(A_1)}=H_{0,-1}+56H_{\xi,-\frac{1}{4}}+2H_{\zeta,-\frac{3}{4}}.\]
Note that $H_{0,-1} = 2H_n +2H_u, H_{\xi,-\frac{1}{4}}=2H_u,  H_{\zeta,-\frac{3}{4}}=2H_{A_2}$.
This proves the relation (\ref{hodgerel1}) for $n=1$. The same arguments apply to the remaining cases and we omit the repetitive computation.
\end{proof}
\begin{rem}
    The coefficient of $\lambda$ in relation (\ref{hodgerel1}) is just the weight of the Borcherds' product modular form $\Phi_\Lambda$  in \cite{Bor} associated with the lattice $\Lambda=( E_7 \oplus A_n)^\perp_{\II}$, i.e., \[ 
   12+ \frac{1}{2}R(E_7\oplus A_n))=12+ \frac{126+n(n+1)}{2}=75+\frac{n(n+1)}{2}.\] 
    %Here $12$ comes from  the fact that $\Phi_\Lambda$ is the pullback of a weight $12$ modular form $\Phi_{12}$ on $\cD_{\II}$.
\end{rem}

\begin{pro}\label{proppullA}
 Let  $\iota  : \sh(\Lambda_{A_{n+1}}) \rightarrow  \sh(\Lambda_{A_n})$
be the natural morphism induced from lattice embedding, then for $n \le 4$,  we have 
\begin{enumerate}
    \item $ \iota^\ast H_{\rm {n}}(\Lambda_{A_n})=H_{\rm {n}}(\Lambda_{A_{n+1}})+H_{A_{n+2}}(\Lambda_{A_{n+1}})$;
    \item $ \iota^\ast H_{\rm u}(\Lambda_{A_n})=H_{\rm u}(\Lambda_{A_{n+1}})$ and $ \iota^\ast H_{D_{n+1}}(\Lambda_{D_n})=H_{D_{n+2}}(\Lambda_{D_{n+1}})$.
\end{enumerate}
In particular, we have 
\begin{equation} \label{adjuctionAn}
\iota^\ast H_{A_{n+1}}=
\begin{cases}
    %- \lambda_{\sh(\Lambda_{A_{n+1}})}+H_{A_{n+2}}(\Lambda_{A_{n+1}}), & n= 1\\
    - \lambda_{\sh(\Lambda_{A_{n+1}})}+H_{A_{n+2}}(\Lambda_{A_{n+1}})+2(n-1)H_{D_{n+2}}(\Lambda_{A_{n+1}}), & n=1, 2;\\
    - \lambda_{\sh(\Lambda_{A_{n+1}})}+H_{A_{n+2}}(\Lambda_{A_{n+1}})+H_{D_{n+2}}+4(n-3)H_{E_{n+2}},& n=3, 4.\\
    %-\lambda_{\sh(\Lambda_{A_5})}+H_{A_6}(\Lambda_{A_5})+H_{D_6}(\Lambda_{A_5}) +4(n-3)H_{E_6}(\Lambda_{A_5}),& n=3, 4.
\end{cases}
\end{equation}
\end{pro}
\begin{proof}
We first show the  pullback formula for the nodal divisor $ H_{\rm {n}}(\Lambda_{A_n})$. 
Fixed an $A_{n+1}$-type vector $v\in \Lambda_{A_n}$ defined in Definition \ref{Avector}.  
Let $w$ be a nodal vector ( i.e., $w^2=-2,\ \divi(w)=1$) in $\Lambda_{A_n}$ such that $\langle v,w\rangle$ is negative definite. 
We claim the projection $\pi_{v^\perp}(w)\in v^\perp=\Lambda_{A_{n+1}}$  of $w$ (see Definition \ref{proj_vector}) is either nodal or $A_{n+2}$-type. By the projection decomposition  $ m\cdot w=a_w\cdot v+ \pi_{v^\perp}(w)$ in  (\ref{vectorprojection}), 
   \[ -2m^2=-(n+1)(n+2)a_w^2+\pi_{v^\perp}(w)^2. \]
Since the negativity of the sublattice $\langle v,w\rangle$ is equivalent to $\pi_{v^\perp}(w)^2<0$, we have either $a_w^2=0$ or $a_w^2=1$ as $m^2 \le N^2=(n+2)^2$. If $a_w^2=0$, then $m\cdot w=\pi_{v^\perp}(w) \in v^\perp$. Actually, $m=1$ since $v^\perp \hookrightarrow \Lambda_{A_n}$ is a primitive embedding. So in this case we obtain a nodal vector $\pi_{v^\perp}(w)$ in $v^\perp$. If $a_w^2=1$, then the only possibility for $m$ is $m=N=n+2$ and thus 
\[ \pi_{v^\perp}(w)^2=-2(n+2)^2+(n+1)(n+2)=-(n+2)(n+3) \]
Moreover, the divisibility of $\pi_{v^\perp}(w)$ is $n+2$. Thus, $\pi_{v^\perp}(w)\in v^\perp=\Lambda_{A_{n+1}}$ is an $A_{n+2}$-type vector.  This finishes the proof of the claim. The claim implies 
\begin{equation*}
    \iota^\ast H_{\rm {n}}(\Lambda_{A_n})=a_{\rm n}\cdot H_{\rm {n}}(\Lambda_{A_{n+1}})+a_{A_{n+1}}\cdot H_{A_{n+2}}(\Lambda_{A_{n+1}})
\end{equation*}
To finish the proof of (1), it remains to show the coefficients $a_{\rm n}=1$ and $a_{A_{n+1}}=1$.  We follow the strategy of \cite[Section 5.2]{LO19} and  \cite[Section 8.4]{FLLST} by computing the intersection multiplicity of germs at a very general point in $\iota( H_{\rm {n}}(\Lambda_{A_n})$.  Let $[z] \in \cD_{\Lambda_{A_n}} \cap v^\perp$ be a very general point so that there is an isomorphism of local germs \[   \big (\widetilde{O}(\Lambda_{A_n}) \setminus \cD_{\Lambda_{A_n}} , \widetilde{O}(\Lambda_{A_n})[z] \big )\  \cong\  \big ( \cD_{V^\perp},[z] \big )   \times \big ( G \setminus V_\bC,0 \big )\]  
where $G=\stab([z])/\{\pm Id\}$ and $V=\spa_\bQ\{v,w\} \cap \Lambda_{A_n}$ is the rank $2$ sublattice.   Let $x_1,x_2$ be the coordinate  under the integral basis of $V$.  %If the projection of the nodal vector $w\in \Lambda_{A_n}$ onto $v^\perp=\Lambda_{A_{n+1}}$ is still nodal type, 
Then  the Gram matrix of $V$ is the following 
\[   \left( \begin{array}{cc c }    
-(n+1)(n+2)   & a    \\ 
  a &  -2   \\
   \end{array}\right).   \]
Then $a=0$ if the projection $\pi_{v^\perp}(w)$ of %the nodal vector w\in \Lambda_{A_n}$ 
onto $v^\perp=\Lambda_{A_{n+1}}$ is still nodal type, and $a=n+1$ if $\pi_{v^\perp}(w)$ is an $A_{n+2}$-type. Thus  by \cite[Example 8.13]{FLLST}, $a_{\rm n}=1$ and $a_{A_2}=1$, see also the similar computation in  \cite[Proposition 5.2.5.]{LO19}.

The similar argument works for the pullback of unigonal type divisors $H_{\rm u}(\Lambda_{A_n})$ and $D_{n+1}$-type divisors $H_{D_{n+1}}(\Lambda_{D_n})$ and thus proves (2).  We leave the details for these two cases to interested readers.
   
Finally combining the above pullback results in (1) and (2),  the pullback formula (\ref{adjuctionAn}) follows by taking  the difference of Hodge relations in (\ref{hodgerel1}).
\end{proof}
\begin{rem}
    If $ \iota$ is a closed embedding and the image $H_{A_{n}}$ is a  normal divisor on $\sh(\Lambda_{A_n})$, one can also compute the pullback $\iota^\ast H_{A_{n}}$ directly  via adjunction formula  
    \[ \iota^\ast (K_{\sh(\Lambda_{A_n})}+H_{A_{n}})=K_{\sh(\Lambda_{A_{n+1}})} \] combined with canonical bundle formula of $\sh(\Lambda_{A_n})$
    \[ K_{\sh(\Lambda_{A_n})}=(19-n) \lambda_{\sh(\Lambda_{A_n})}-\frac{1}{2}B_{\sh(\Lambda_{A_n})}.\]
    In this ideal situation, $\iota^\ast H_{A_{n}}=-\lambda_{\sh(\Lambda_{A_n})}+\frac{1}{2}(\iota^\ast B_{\sh(\Lambda_{A_n})}-B_{\sh(\Lambda_{A_{n+1}})} )$. 
    But in general it is hard to check the normality and it is possible that $\iota$ fails to be a closed embedding.
\end{rem}

\begin{pro}
The Picard number of $\sh(A_n')$ is given by $\rho(\sh(A_n')^\ast)=4$  for $n=5,6,7$.
Moreover, there are the first relations 
\begin{equation}\label{relaA'1}
    \begin{split}
        90 \lambda_{\sh(A_5')}=H_{\rm n}+14 H_{A_6'}+13 H_{\bar{4},-\frac{1}{3}}+ 78 H_{\Bar{3},-\frac{1}{8}} \\
    89 \lambda_{\sh(A_6')}=H_{\rm n}+91H_{\bar{5},-\frac{1}{11}}+28 H_{\rm u'}+15 H_{A_6''}
    \end{split}
\end{equation}
and the second relation
\begin{equation}\label{relaA'2}
    \begin{split}
       \lambda= H_{A_7'}-H_{\bar{3},-\frac{3}{11}}-H_{A_6''}.
     %  H_{A_7'}=&\lambda+H_{-\frac{3}{11}}+H_{A_6''}
    \end{split}
\end{equation}
In particular, $ H_{A_7'}$ is big on $\sh(\Lambda_{A_6'})$ if $H_{\bar{3},-\frac{3}{11}}$ and $H_{A_6''}$ are birationally contractible on $\sh(\Lambda_{A_6'})$.
\end{pro}

\begin{pro}
The Picard number of $\sh(A_n'')$ is given by $\rho(\sh(A_n'')^\ast)=9-n$ for $n=6,7$  and there are relations 
\begin{equation}\label{relaa'''}
    \begin{split}
       & 104 \lambda = H_{\rm n}+28 H_{\bar{2},-\frac{1}{4}}+ 16 H_{A_7''} \ \ \ \hbox{on} \  \ \sh(A_6''), \\
     &   120 \lambda=H_{\rm n}+45 H_{A_8''} \ \ \ \hbox{on} \  \ \ \sh(A_7'').
    \end{split}
\end{equation}
\end{pro}

%Denote \begin{equation}\label{morpull}  j: \ \sh(A_5') \longrightarrow \sh(A_5),\ \ \  \iota': \sh (A_6') \longrightarrow \sh(A'_{5})\end{equation}

\begin{pro}\label{pullA'}
Let $j:  \sh(A_5') \longrightarrow \sh(A_5)$ and $  \iota': \sh (A_6') \longrightarrow \sh(A'_{5})$ be the induced  morphism.  %$j$ and $ \iota'$ in (\ref{morpull}) will induce  the pullback 
 %\[j^\ast: \pic(\sh(A_5)) \rightarrow \pic(\sh(A_5')),\quad \iota'^\ast: \pic(\sh(A_{5}')) \rightarrow \pic(\sh(A_6')). \] 
Then the following pullback formulas hold
    \begin{enumerate}
          \item $j^\ast  H_{\rm n}= H_{\rm n}+ H_{\Bar{3},-\frac{1}{8}} +  \frac{1}{2}H_{\bar{6},\frac{1}{2}}$,  $ \ j^\ast H_{A_6}=2H_{A_6}'$;
          \item $j^\ast H_{D_6}=0,\ \ \ \  j^\ast H_{E_6}=H_{\bar{3},-\frac{1}{8}}$ ;
        \item $ \iota'^\ast H_{\rm n}=H_{\rm n}+H_{A_7'}$.
    \end{enumerate}
   In particular, the above will imply 
   \begin{equation} \label{pbA6'}
       \iota'^\ast H_{A_6'}=-\lambda+H_{A_7'}+Z
   \end{equation}
where $Z$ is linear combination Heegner divisors $\iota'^{\ast}H_{\bar{4},-\frac{1}{3}}, \iota'^{\ast}H_{\bar{3},-\frac{1}{8}}, H_{\bar{3},-\frac{3}{11}}, H_{\bar{5},-\frac{1}{11}}$ .
  
\end{pro}
\begin{proof}

We first deal with the pullback $j^\ast$. Fix an $A_5'$-type vector $v \in \Lambda_{A_5}$, then $v^\perp \cong  \Lambda_{A_5}'$ and $[\Lambda_{A_5}: \bZ v \oplus v^\perp]= 2$. Let $w$ be an $A_6$-type vector, i.e., $w^2=-42$ and $\divi(w)=6$. 
Let $x_w:=\pi_{v^\perp}(w)$ be the projection on $v^\perp=\Lambda_{A_5}'$. 
Then $mw=av+x_w$ for $m=1$ or $m=2$.  Note that $m \langle w,x_w \rangle =x_w^2$ implies that $x_w^2$ is also a multiple of $\divi(w)=6$. Then the only possibility is $m=2$  and $|a|=3$ by considering the divisibility of the vectors $v$ and $w$, and their corresponding elements in the discriminant group.
In this case, $x_w^2=-132$ and $x_w$ is an $A_6'$-type vector. This shows $j^\ast H_{A_6}$ is supported  on $H_{A_6'}$.  
If $w$ is a nodal vector, then the projection $x_w$  is either a nodal vector or a primitive vector with $x_w^2=-4$ ($x_w=2w\pm v$) and $\disc((x_w)_{v^\perp}^\perp)=2$ or $4$. 
This proves (1). 

Similarly, if $w$ is a $D_6$-type vector, by considering the projection $x_w$ of $w$ one can see that the intersection $H_{D_6} \cap H_{A_5'}$ is empty set-theoretically in $\sh(A_5)$.  
If $w$ is an $E_6$-type vector, then $x_w^2=-4$ with $\divi(x_w) =4$.  This proves $(2)$. 

Now we consider the  pullback $\iota'^\ast$.  Fixed an $A_6'$-type vector $v \in \Lambda_{A_5'}$. Let $w  \in \Lambda_{A_5'}$  be a nodal vector. Arguments as before will show that the projection $x_w$ of $w$ onto $v^\perp=\Lambda_{A_6'}$ is either a nodal vector ($x_w=w$)  or an $A_7'$-type vector ($x_w=-11w\pm v$) in $v^\perp=\Lambda_{A_6'}$.  The multiplicity of $H_{\rm n}$ and $H_{A_7'}$ on  $ \iota'^\ast H_{\rm n}$ are argued as in the previous case. This shows (3). 

For the last assertion, by the first relation in (\ref{relaA'1}) and the second relation (\ref{relaA'2}), we have 
\[  (89+15)\lambda= H_{\rm n} +15H_{A_7'}+13H_{\bar{3},-\frac{3}{11}} +91 H_{\bar{5},-\frac{1}{11}}\]
in order to replace $H_{A_6}''$ in (\ref{relaA'1}). Then again using the difference of the above relation with the pullback of the first identity in relation  (\ref{relaA'1}), one obtains (\ref{pbA6'}).
\end{proof}

\begin{pro} \label{pullA''}
  Let $j': \sh(A_6'') \longrightarrow \sh(A_6')$ and $\iota'': \sh(A_7'') \longrightarrow \sh(A_{6}'')$ be the natural morphism. 
  Then we have 
  \begin{enumerate}
     \item $j^{'\ast}  H_{\rm n}= H_{\rm n} + H_{\bar{2},-\frac{1}{4}}$, \ $j^{'\ast} H_{A_7'}= H_{\bar{2},-\frac{1}{4}} + H_{A_7''}$;
        \item $ \iota''^\ast H_{\rm n}= H_{\rm n} + H_{A_{8}''}$, \  $ \iota''^\ast H_{A_{7}''}=-\lambda+Z'$ where $Z'$  is the linear combination of Heegner divisors $\iota^{''\ast} H_{\bar{2},-\frac{1}{4}}$ and $H_{A_8''}$ .
    \end{enumerate}
\end{pro}

\begin{proof}
Fixed $v\in \Lambda_{A_6'}$ an $A_6''$-type vector as in Definition \ref{1stmv}.
If $w$ is a nodal vector, then the projection $\pi_{v^\perp}(w)$ can be of the following forms:
\begin{itemize}
    \item $\pi_{v^\perp}(w) = w$.
    \item $\pi_{v^\perp}(w) = 8w \pm v$, $\pi_{v^\perp}(w)^2 = -40, \divi(\projection{v}{w})=8$.
\end{itemize}

If $w$ is an $A_7'$-type vector, then the projection $\pi_{v^\perp}(w)$ can be of the following forms:
\begin{itemize}
    \item $11\pi_{v^\perp}(w) = 4w \pm 3v$, $\pi_{v^\perp}(w)^2 = -8, \divi(\projection{v}{w})=4$.
    \item $11\pi_{v^\perp}(w) = 8w \pm 5v$, $\pi_{v^\perp}(w)^2 = -40, \divi(\projection{v}{w})=8$.
\end{itemize}

Fixed $v\in \Lambda_{A_6''}$ an $A_7''$-type vector, and let $w\in  \Lambda_{A_6''}$ be a nodal vector, then it is easy to check that the projection of $w$ onto $v^\perp$ is either a nodal vector or an $A_8''$-type vector in $v^\perp$.  By taking the difference of the relation in (\ref{relaa'''}), we obtain 
    \begin{equation*}
        \begin{split}
            -32\lambda=&\iota^{''\ast} ( H_{\rm n}+56H_{\bar{2},-\frac{1}{4}}+ 32 H_{A_7''})-(H_{\rm n}+90 H_{A_8''})\\
            =&32\iota^{''\ast} ( H_{A_7''})+56 \iota^{''\ast} H_{\bar{2},-\frac{1}{4}}-90 H_{A_8''}.
        \end{split}
    \end{equation*}
    This finishes the proof of (2).
\end{proof}

\subsubsection{Pullback of $D_n$ towers}

\begin{pro}
The Picard numbers of $\sh(D_n)$ are given by the following  table:
\begin{table}[ht]
    \centering
    \begin{tabular}{ |c|*{9}{c|} }
    \hline
     $n$    & $4$  & $5$  & $6$  & $7$ & $8$ & $9$  & $10$  \\ \hline
    $\rho(\sh(D_n)^\ast)$   & $5$  &  $4$ & $5$ & $4$  & $4$  & $3$ & $2$ \\  \hline 
     \end{tabular}
    %\caption{Picard numbers of $D_n$-towers} \label{Picard of Dtower}
\end{table}

Moreover, there is a relation (called Hodge relation) for $n=4 ,5,6,7$ 
\begin{equation}\label{relaD}
\begin{split}
    (75+n(n-1)) \lambda_{\sh(D_n)} & =H_{\rm n}(\Lambda_{D_n})+57 H_{\rm u}(\Lambda_{D_n})+2nH_{D_{n+1}}(\Lambda_{D_n}) \\
    &\phantom{=} +(2^{n-1} + \delta_{n,6} + 2n \cdot\delta_{n,7}) H_{E_{n+1}}(\Lambda_{D_n}).
\end{split}
\end{equation}
\end{pro}

\begin{rem}
    $H_{D_{7}}(\Lambda_{D_6})$ is a  movable  divisor and thus is not contractible  on $\lambda_{\sh(D_6)}$. 
    This can not be seen directly from the arithmetic side.  
\end{rem}

\begin{pro}[Hodge relations on  $D_n'$-towers $\sh(D_n')$]
The Picard number of $\sh(D_n')$ is given by $\rho(\sh(D_n')^\ast)=10-n,  \ 6\le n \le 9$.
Moreover, there is a relation
\begin{equation}\label{relaD'}
    \begin{split}
       & 103 \lambda_{\sh(D_6')}=H_{\rm n}+14 H_{D_7'}+ 29 H_{(\bar{1},\bar{2}),-\frac{1}{4}} + 92 H_{(\bar{0},\bar{1}),-\frac{1}{8}} \\
      & 117 \lambda_{\sh(D_7')}=H_{\rm n}+ 15 H_{D_8'}+ 135 H_{(\bar{1},\bar{1}),-\frac{1}{12}} \\
       & 132 \lambda_{\sh(D_8')}=H_{\rm n}+33H_{D_9'}  \\
       & 165 \lambda_{\sh(D_9')} = H_{\rm n}
    \end{split}
\end{equation}
\end{pro}

\begin{proof}
    One only needs to notice that for $n=7$, we have $
    H_{(\bar{0},\bar{1}),-\frac{1}{3}} = H_{D_8'} + H_{(\bar{1},\bar{1}),-\frac{1}{12}} $.
\end{proof}

 \begin{pro}
    Let $\iota  : \sh(\Lambda_{D_{n+1}}) \rightarrow  \sh(\Lambda_{D_n})$  
be the natural morphism induced from lattice embedding, then   we have 
\begin{equation*}
    \begin{split}
      \iota^\ast H_{\rm {n}}(\Lambda_{D_n})=& \begin{cases}
    H_{\rm {n}}(\Lambda_{D_{n+1}})+H_{D_{n+2}}(\Lambda_{D_{n+1}}), & n=3; \\
        H_{\rm {n}}(\Lambda_{D_{n+1}})+2H_{D_{n+2}}(\Lambda_{D_{n+1}}) + \delta_{n,5} H_{E_{n+2}}(\Lambda_{D_{n+1}}), & \text{otherwise}.
    \end{cases} \\
    \iota^\ast H_{A_{4}}(\Lambda_{D_3})=& H_{D_4}(\Lambda_{D_{4}}) + 2 H_{E_{5}}(\Lambda_{D_{4}}), \quad  \iota^\ast H_{E_{n+1}}(\Lambda_{D_n})=2 H_{E_{n+2}}(\Lambda_{D_{n+1}}).
    \end{split}
\end{equation*}
Moreover, we get the pullback formula for $n \le 6$:
\begin{equation}\label{pullD}
    \iota^\ast H_{D_{n+1}}(\Lambda_{D_n})=\begin{cases}
        -\lambda_{\sh(D_{n+1})}+ \frac{1}{2} H_{D_{n+2}}, & n=3; \\
        -\lambda_{\sh(D_n)}+ H_{D_{n+2}}(\Lambda_{D_{n+1}}) +  H_{E_{n+2}} (\Lambda_{D_{n+1}}), & n=6; \\
        -\lambda_{\sh(D_n)}+ H_{D_{n+2}}(\Lambda_{D_{n+1}}) , & \text{otherwise}.
    \end{cases}
\end{equation}
\end{pro}
\begin{proof}  
The argument is parallel to the case of $A_n$ towers.
We give the proof of (1) here and leave the rest to interested readers.
Fix a $D_{n+1}$-type vector  $v\in \Lambda_{D_n}$ and let $w\in \Lambda_{D_n}$ be a nodal vector. Then the projection $x_w$ of $w$ onto $v^\perp=\Lambda_{D_{n+1}}$ is either $w$ itself or a $D_{n+1}$-type vector $x_w=2w\pm v$.
In the first case, by \cite[Proposition 1.4.6]{LO19}, we have that $w$ is a nodal vector and moreover an $E_7$-type vector when $n=5$.
The multiplicity computation is identical to that in the case of the $A_n$-tower in  Proposition \ref{proppullA}. This proves (1).  

The last pullback formula (\ref{pullD}) is obtained from relation (\ref{relaD}) by combining (1) and (2). 
\end{proof}

\begin{pro} 
Let $j: \sh(\Lambda_{D_{7}'}) \rightarrow  \sh(\Lambda_{D_7})$ and $\iota_m' : \sh(\Lambda_{D_{m+1}'}) \rightarrow  \sh(\Lambda_{D_m}')$
be the natural morphisms induced from lattice embedding, then
  \begin{enumerate}
      \item $j^\ast H_{\rm n}= H_{\rm n}+ H_{D_8'}$ ,\  \  \ $j^\ast H_{D_8}=H_{D_8'},\ \  \ j^\ast H_{E_8}=j^\ast H_{\rm u}= H_{(\Bar{0},\Bar{1}),-\frac{1}{12}} $.
      
      \item  $\iota_7'^\ast H_{\rm n}= H_{\rm n}+3 H_{D_8'}$, \  \  \  $  \iota_7'^\ast H_{(\bar{1},\bar{1}),-\frac{1}{12}}=0 $.
  \end{enumerate}  In particular, 
\begin{equation}\label{pullD'}
\iota'^\ast H_{D_{7}'}=-\lambda +H_{D_8'},\ \ \ \iota'^\ast H_{D_{8}'}=-\lambda +2 H_{D_9'},\ \ \ \iota'^\ast H_{D_{9}'}=-\lambda.
\end{equation}  
\end{pro}
\begin{proof}
(1). Fix a $D_{7}'$-type vector  $v\in \Lambda_{D_7}$ and let $w\in \Lambda_{D_7}$ be a $E_8$ vector. 
Take $w_1\in \Lambda_{D_7}$ a $E_8$-type vector, and $w_2\in \Lambda_{D_7}$ a unigonal vector then we have 
\[ 2\pi_{v^\perp}(w_1) = 3w_1 - v, \ \
        \pi_{v^\perp}(w_2) = 3w_2 - v,\ \ H_{\pi_{v^\perp}(w_1)}= H_{\pi_{v^\perp}(w_2)}=H_{(\Bar{0},\Bar{1}),-\frac{1}{12}}. \]

(2). Fix $v\in \Lambda_{D_7'}$ a $D_8'$-vector, and take $w\in\Lambda_{D_7'}$ with \[
    w^2=-6, \divi(w) =6, H_w = H_{(\bar{1},\bar{1}),-\frac{1}{12}}.
\]
We have $ mw=av+ l \cdot \pi_{v^\perp}(w) 
$ where $m\in\{1,2\}$.
Then the only possibility is $m=2$, $a=1$, $l^2 \pi^2 =-18$. 
Note that \[
    12 \,\big|\, (2w\cdot \pi = l\cdot \pi^2)\, \big|\, l^2 \pi^2 =-18
\] as $\divi(w) =6$.
This gives the contradiction and $  \iota_7'^\ast H_{(\bar{1},\bar{1}),-\frac{1}{12}}=0 $.

The remaining cases are similar and we omit here.
Combining relation (\ref{relaD'}),  the pullback formula (\ref{pullD'}) is obtained as before. 
\end{proof}

\subsubsection{Pullback of $E_n$-tower}
\begin{pro}%The Picard group of $\sh(E_n)$ is generated by the integral divisors $ H_{\rm {n}}(\Lambda_{E_n}), H_{\rm u}(\Lambda_{E_n}) $ and $H_{E_{n+1}}(\Lambda_{E_n})$.
The Picard number of $\sh(E_n)$ is given  by $\rho(\sh(E_n)^\ast)=6-\rounddown{\frac{n}{2}}, \  n=6,7, 8.$
%\begin{equation*}\rho(\sh(E_n)^\ast)=6-\rounddown{\frac{n}{2}},\ \  n=6,7, 8.\end{equation*}
Moreover, there are relations 
\begin{equation*}\label{pullE1}
\begin{split}
   & 111 \lambda_{\sh(E_6)}=H_{\rm {n}}(\Lambda_{E_n})+ 57H_{\rm u}(\Lambda_{E_n}) + 27 H_{E_{7}}(\Lambda_{E_6}), \\
   & 138 \lambda_{\sh(E_7)}=H_{\rm {n}}(\Lambda_{E_n})+57H_{\rm u}(\Lambda_{E_n})+57H_{E_{8}}(\Lambda_{E_7}), \\
    & 195  \lambda_{\sh(E_8)}=H_{\rm {n}}(\Lambda_{E_n})+57H_{\rm u}(\Lambda_{E_n}).
\end{split}  
\end{equation*}
\end{pro}

\begin{pro} 
    Let  $ \iota  : \sh(\Lambda_{E_{n+1}}) \rightarrow  \sh(\Lambda_{E_n}),\ \  \nu  : \sh(\Lambda_{E_{6}}) \rightarrow  \sh(\Lambda_{D_5})$ 
be the natural morphisms induced from lattice embedding, then   we have 
\begin{enumerate}
    \item $\nu^\ast  H_{\rm {n}}(\Lambda_{D_5})=H_{\rm {n}}+H_{E_7}$,\ \  \ \ $\nu^\ast H_{D_6}=H_{E_7}\ $  and $\ \nu^\ast H_{\rm u}(\Lambda_{D_5})=H_{\rm u}(\Lambda_{E_6})$
    
    \item $ \iota^\ast H_{\rm {n}}(\Lambda_{E_6})=H_{\rm {n}}(\Lambda_{E_{7}})+3 H_{E_{8}}(\Lambda_{E_{n+1}})$,\ \    \ \ $ \iota^\ast H_{\rm u}(\Lambda_{E_n})=H_{\rm u}(\Lambda_{E_{n+1}})$ 
\end{enumerate}
In particular, the following pullback formulas hold
\begin{equation}\label{pullE2}
\begin{split}
      &\nu^\ast H_{E_6}(\Lambda_{D_5}) =-\lambda_{\sh(E_6)}+ H_{E_7} ,\\  &\iota^\ast H_{E_{n+1}}(\Lambda_{E_n})=-\lambda_{\sh(E_{n+1})}+2(7-n)H_{E_{n+2}}(\Lambda_{E_{n+1}}), \ n=6,7. %\\ & \iota^\ast H_{E_{7}}(\Lambda_{E_6})=-\lambda_{\sh(E_{7})}+ 2 H_{E_{8}}(\Lambda_{E_{7}}),  \\ & \iota^\ast H_{E_{8}}(\Lambda_{E_7})=-\lambda.
\end{split}
\end{equation}
\end{pro}
\begin{proof}
The argument is similar to the previous one. 
We fix $v\in \Lambda_{D_5}$ be a vector of $E_6$-type as in Definition \ref{vectoronD} (3) and $v^\perp \cong \Lambda_{E_6}$.  
If $w\in  \Lambda_{D_5}$ is a nodal vector, i.e., $w^2=-2,\ \divi(w)=1$,  then 
\[ m w=a v+ \pi_{v^\perp}(w),\ m |3,\ \ a\in \bZ \ \]
since $[ \Lambda_{D_5}:\bZ v\oplus v^\perp]=3$. Clearly, $m=1$ implies $\pi_{v^\perp}(w)=w$. If  $m=3$, then $-18=-12a^2+(\pi_{v^\perp}(w))^2$. The negativity of $\pi_{v^\perp}(w)$ implies that the only solution is $a= \pm 1$  and $\pi_{v^\perp}(w)=-6$. One can also check $(\pi_{v^\perp}(w))^\ast =\eta_E$, i.e., $\pi_{v^\perp}(w) \in v^\perp \cong \Lambda_{E_6}$ is $E_7$-type vector. This proves that the support of $\nu^\ast  H_{\rm {n}}(\Lambda_{D_5})$ is $H_{\rm {n}}$ and $ H_{E_7}$. 
The multiplicity is $1$ by the same argument as in the proof of Proposition \ref{proppullA}. The remaining cases follow by same arguments.  This proves part (1).  
Part (2) follows from the same computation. We leave the details to the interested reader.
 
\end{proof}

\subsection{Pullback on unigonal divisor towers}
\begin{pro}\label{unigonalpicard}
The Picard group of $\sh(U_1)$ is generated by $H_{\rm {n}}(U_1)$  and $H_{A_{2}}(U_1)$, and  the Picard numbers is $2$.
Moreover, there  are  relations 
\begin{equation}\label{relationU}
\begin{split}
133 \lambda_{U_1} &= H_{\rm n}+ 2 H_{A_{2}}, \\
 134 \lambda_{ A_1 \oplus A_1 } &= H_n + 4 H_{\zeta_1+\zeta_1',-\frac{1}{2}}^{\mathbf{red}} + 2 H_{A_1\oplus A_2}, \\
  32 \lambda_{ A_1 \oplus A_1 } &=H_{A_1 \oplus A_2}- 8 H_{\zeta_1+\zeta_1',-\frac{1}{2}}, \\
 34 \lambda_{ A_1 \oplus A_2} &=H_{A_2 \oplus A_2}+ 2 H_{A_1\oplus A_3} - 14 H_{\zeta_1+\zeta_2,-\frac{5}{12}},\\
  36  \lambda_{ A_2 \oplus A_2} &= 2 H_{A_2 \oplus A_3}-12 (H_{\zeta_2+\zeta_2',-\frac{1}{3}} + H_{\zeta_2- \zeta_2',-\frac{1}{3}}),\\
  19\lambda_{A_2 \oplus A_3  }  &= H_{A_2 \oplus D_4}+H_{A_3+\oplus A_3}+H_{A_2\oplus A_4}-12H_{\zeta_2+2\zeta_3,-\frac{1}{6}}-5H_{\zeta_2\pm \zeta_3,-\frac{7}{24}}, \\
 21 \lambda_{ A_2 \oplus D_4}  &= H_{\zeta_2,-\frac{2}{3}}+3 H_{\eta,-\frac{1}{2}}-16(H_{\eta'+\zeta_2,-\frac{1}{6}}+H_{\eta+\eta'+\zeta_2,-\frac{1}{6}})+2 H_{\eta+\zeta_2,-\frac{1}{6}}.
\end{split}
\end{equation}
\end{pro}

\begin{rem} \label{Humodify}
 $H_{A_{2}}$ is big and thus not birationally contracted. Indeed in \cite{psw2} we prove K-moduli space is isomorphic to $\cF(s)$.  The bigness of  $H_{A_{2}}$ can be seen from the wall-crossing results on the K-moduli side. 
\end{rem}

\begin{pro}\label{unigonalpull}
Let $\nu: \sh(U_1) \longrightarrow \cF$ and $\iota:  \sh(U_{n+1}) \longrightarrow \sh(U_n)$ be the morphisms induced by the lattice embedding. 
Then 
\begin{enumerate}
    \item $\nu^\ast\lambda=\lambda,\  \nu^\ast H_{\rm n}= H_{\rm n}$, \ \   $\nu^\ast H_{h}= H_{A_2}$,\ \  $\nu^\ast (H_{\rm u})=-\lambda$;
    
    \item $\iota^\ast H_{A_2}= 2 H_{\zeta^h+\zeta,-\frac{1}{2}}^{\mathbf{red}} +H_{A_1 \oplus A_2}$ on $\  \sh(U_2)$;
    
    \item $\iota ^\ast H_{A_1 \oplus A_2}= -2\lambda+  H_{A_2 \oplus A_2}+Z_{U_3}\ $ on $\ \sh(U_3)$;
    
    \item  $\iota ^\ast H_{A_2 \oplus A_2}=-\lambda+H_{A_2 \oplus A_3}+Z_{U_4}\ $  on $\ \sh(U_4)$; 
    \item  $\iota ^\ast H_{A_2 \oplus A_3}=-\lambda+H_{A_2 \oplus D_4}+Z_{U_5}\  $ on $\ \sh(U_5)$;   
    \item   $\iota ^\ast H_{A_2 \oplus D_4}=-2\lambda+Z_{U_6}\ $ on   $\ \sh(U_6)$
\end{enumerate}
where each $Z_{U_i}$ is a linear combination of possible  Heegner divisors with smaller discriminant. 
\end{pro}
\begin{proof}
    Fixed a unigonal vector $v\in \Lambda_{A_1}=(E_7\oplus A_1)^\perp_{\II}$, i.e., $v^2=-2, \divib_{\Lambda_{A_1}}(v)=2$. Then $\Lambda_{A_1}=\bZ v \oplus v^\perp$ since  $v^\perp\cong (E_8\oplus A_1)^\perp_{\II}$. Thus a nodal vector $w$ in $\Lambda_{A_1}$ such that the sublattice $\langle v,w \rangle$ is negatively definite must be a nodal vector in $v^\perp$.  In other words, the projection of a nodal vector in $\Lambda_{A_1}$ onto $v^\perp$ must be nodal. This shows 
$ \nu^\ast H_{\rm n}= H_{\rm n}$.  Similarly, let $w$ be a hyperelliptic vector in $\Lambda_{A_1}$, i.e., $w^2=-6, \divib_{\Lambda_{A_1}}(w)=2$.  Then \[ v^2=-2,\ \ v.w=a,\ \ w^2=-6.\]
\begin{comment}
    \[ \left( 
    \begin{array}{c|cc}
         & v & w   \\ \hline
     v & -2 & a \\
       w & a  & -6
    \end{array}  
    \right)\]
\end{comment}
So $|a| \le 3$ as the sublattice $\langle v,w \rangle$ is negatively definite.  By the divisibility of $v$ and $w$, $a=0$ or $|a|=2$. The latter is impossible. Indeed, if $a=2$, then $w':=v+w\in v^\perp $ is the projection of $w$ onto $v^\perp$ and clearly $\divib_{v^\perp}(w')=2m,\ w'^2=-4$; therefore $(w'^\ast)^2=(\frac{w'}{2m})^2=-\frac{1}{m^2}$, a contradiction since the discriminant group of $v^\perp$ has only two elements $0$ and $\xi$ with $\xi^2=-\frac{1}{2}$. The same reasoning shows $a=-2$ is impossible.  This proves $a=0$ and thus $w\in v^\perp=U_1$ with  $w^2=-6, \divib_{v^\perp}(w)=2$, exactly the $A_2$-type vector in definition \ref{uAn}. Thus we finish the proof of $\nu^\ast H_h=H_{A_2}$. Recall we have the relation  $ 76\lambda=H_{\rm n} +57H_{\rm u}+2H_{h} $ on $\cF$.
%\begin{equation*}76\lambda=H_{\rm n} +57H_{\rm u}+2H_{h} \ \ \hbox{on}\ \cF.\end{equation*}
Combining the first identity in  (\ref{relationU}), we get $\nu^\ast (H_{\rm u})=-\lambda$. This proves (1).

Fixed a nodal vector $v\in U_1=(E_8\oplus A_1)^\perp_{\II}$ and so that $v^\perp =(E_8\oplus A_1 \oplus A_1 )^\perp_{\II}$, then $[U_1:\bZ v \oplus v^\perp]=2$. By similar arguments as in (1), we have $\iota^\ast H_{A_2}=H_{A_1 \oplus A_2}$.
Then combining second identity in (\ref{relationU}), we obtain the pullback formula in (2). The remaining cases follow by similar lattice arguments and by applying the identities in (\ref{relationU}). We omit the details here.
\end{proof}

\vspace{0.1cm}

\subsection*{Data availability} There is no data associated with this manuscript.

\subsection*{Conflict of interest} The author declares that there are no conflict of interest.

\bibliographystyle{alpha}
\bibliography{birational}
\end{document}